\documentstyle{article}

\newcommand{\inlabel}[1]{\refstepcounter{equation} \label{#1}
\hspace{1.5ex} (\theequation)}

\newcommand{\N}{{\rm I\kern-.5ex N}}
\newcommand{\Z}{{\sf \vrule height 1.55ex depth-1.2ex
width.03em\kern-.11em Z
\kern-.9ex Z\kern-.11em\vrule height 0.3ex depth0ex width.03em}}
\newcommand{\Q}{{\rm\kern.2ex\vrule height1.55ex depth-.05ex
width.03em\kern-.7ex Q}}
\newcommand{\R}{{\rm I\kern-.5ex R}}
\newcommand{\Rvar}{{\rm I\kern-.5ex R}}
\newcommand{\C}{{\rm\kern.3ex\vrule height1.55ex depth-.05ex
width.03em\kern-.7ex C}}
\newcommand{\Cvar}{{\, \rm\kern.3ex\vrule height1.1ex depth-.05ex
width.03em\kern-.7ex C}}

\newcommand{\nabp}{\nab \hspace{-1.05ex}
\rule[.5ex]{.2ex}{.8ex}   \hspace{1.05ex}}

\newcommand{\spat}{\hspace{4ex}}

\newcommand{\Ga}{\Gamma}

\newcommand{\nab}{\nabla}

\newcommand{\cL}{{\cal L}}
\newcommand{\cF}{{\cal F}}
\newcommand{\cG}{{\cal G}}

\newcommand{\cN}{{\cal N}}
\newcommand{\cM}{{\cal M}}

\newcommand{\od}{\odot}
\newcommand{\ot}{\otimes}

\newcommand{\la}{\Lambda}

\newcommand{\om}{\omega}
\newcommand{\io}{\iota}
\newcommand{\vfi}{\varphi}
\newcommand{\vep}{\varepsilon}

\newcommand{\al}{\alpha}
\newcommand{\be}{\beta}

\newcommand{\th}{\theta}

\newcommand{\si}{\sigma}
\newcommand{\Mfi}{{\cal M}_{\vfi}}
\newcommand{\Nfi}{{\cal N}_{\vfi}}

\newcommand{\Nps}{{\cal N}_{\psi}}

\newcommand{\Npsi}{{\cal N}_{\psi}}
\newcommand{\lafi}{\la_\vfi}
\newcommand{\laps}{\la_\psi}
\newcommand{\lapsi}{\la_\psi}
\newcommand{\pifi}{\pi_\vfi}
\newcommand{\pips}{\pi_\psi}

\newcommand{\cU}{{\cal U}}

\newcommand{\text}[1]{\mbox{#1}}

\newcommand{\cst}{\text{C}$\hspace{0.1mm}^*$}
\newcommand{\wst}{\text{W}$\hspace{0.1mm}^*$}

\newcommand{\qed}{\ \hfill \rule{2mm}{2mm}}
\newenvironment{demo}{\medskip\noindent\bf Proof :\ \
\rm}{\qed\bigskip\par }

\newtheorem{definition}{Definition}[section]
\newtheorem{proposition}[definition]{Proposition}
\newtheorem{lemma}[definition]{Lemma}
\newtheorem{corollary}[definition]{Corollary}
\newtheorem{remark}[definition]{Remark}
\newtheorem{theorem}[definition]{Theorem}

\newtheorem{notation}[definition]{Notation}
\newtheorem{result}[definition]{Result}

\begin{document}

\begin{center}
\LARGE\bf Weight theory for \cst-algebraic quantum groups
\end{center}

\bigskip\bigskip

\begin{center}
\rm Johan Kustermans

Department of Mathematics

University College Cork

Western Road

Cork

Ireland

\medskip

e-mail : johank@ucc.ie

\bigskip\bigskip

Stefaan Vaes\footnote{Research Assistant of the
Fund for Scientific Research -- Flanders (Belgium)\ (F.W.O.)}

Department of Mathematics

KU Leuven

Celestijnenlaan 200B

Leuven

Belgium

\medskip

e-mail : Stefaan.Vaes@wis.kuleuven.ac.be

\bigskip\bigskip

\bf January 1999 \rm
\end{center}

\bigskip

\begin{abstract}
\noindent In this paper, we collect some technical results about weights on \cst-algebras which are useful in de theory of locally compact quantum groups in the \cst-algebra framework (see \cite{JK-Va}).
We discuss the extension of a lower semi-continuous weight to a normal weight following S. Baaj
(see \cite{Baa3}), look into slice weights and their KSGNS-constructions and investigate the tensor product of weights together with a partial GNS-construction for such a tensor product.
\end{abstract}

\bigskip\medskip

\section*{Introduction}

The theory of lower semi-continuous weights on \cst-algebras is the non-commutative generalization of the theory of regular Borel measures on locally compact spaces. The first and most fundamental results were proven by F. Combes in \cite{Comb1} and \cite{Comb}. Where in \cite{Comb}, the focus lies on the general theory of lower semi-continuous weights, \cite{Comb1} revolves around an important subclass of these lower semi-continuous weights, the KMS weights. The KMS property guarantees that we have enough control over the non-commutativity of the \cst-algebra under the weight. Most of the results in this paper do not assume the weights to satisfy this extra KMS condition.

The theory of normal weights  on a von Neumann algebra (which has been investigated more intensively) on the other hand is the non-commutative version of measure theory in general. Some references here are \cite{Haa},\cite{Comb1},\cite{Con},\cite{Pe-Tak},\cite{Va2},\cite{Stramod},\cite{Stra},\cite{Haa1}.

\medskip

The main reason for writing this paper stems from the fact that we need a theory of weights which is powerful enough to deal with the Haar weights appearing in the theory of locally compact quantum groups in the \cst-algebra framework. As such, this paper accompanies our paper \lq Locally compact quantum groups\rq\ \cite{JK-Va} in which we propose a relatively simple definition of a locally compact quantum group in the reduced setting (i.e. in which the Haar weights are faithful). But it is also intended to be used in non-reduced approaches to locally compact quantum groups (e.g. the universal setting). The only aim of this paper is to provide a lot of relatively simple technical results which one does not want to deal with in a paper devoted to quantum group theory. Other references for these kind of technical results are \cite{Comb},\cite{Comb1},\cite{E},\cite{E-V},\cite{Q-V},\cite{JK1}.

\medskip

We will cover the following topics:
\begin{enumerate}
\item An overview of weight theory on \cst-algebras.
\item The GNS-construction of the normal extension of a lower semi-continuous weight (a result which is essentially due to S. Baaj \cite{Baa3}).
\item Elementary properties of slice weights and their KSGNS-constructions.
\item The tensor product of two lower semi-continuous weights and the partial GNS-construction of this tensor product.
\item A basic proof of a typical integration/strict convergence result.
\end{enumerate}

\bigskip\medskip

\section*{Notations and conventions}

For any subset $X$ of a Banach space $E$, we denote the linear span by $\langle X \rangle$, its closed linear span by $[X]$.

If $I$ is set, $F(I)$ will denote the set of finite subsets of $I$. We turn it into a directed set by inclusion.

\medskip

All tensor products in this paper are minimal ones. This implies that the tensor product functionals separate points of the tensor product (and also of its multiplier algebra). The completed tensor products will be denoted by $\ot$, the algebraic ones by $\odot$.

\medskip

The multiplier algebra of a  \cst-algebra $A$ will be denoted by $M(A)$. If $A$ and $B$ are two \cst-algebras, then the tensor product $M(A) \ot M(B)$ is naturally embedded in $M(A \ot B)$.

Consider two \cst-algebras $A$ and $B$ and a linear map $\rho : A \rightarrow M(B)$. We call $\rho$ strict if it is norm bounded and strictly continuous on bounded sets. If $\rho$ is strict, $\rho$ has a unique linear extension $\overline{\rho} : M(A) \rightarrow M(B)$ which is strictly continuous on bounded sets (see proposition 7.2 of \cite{JK3}). The resulting $\overline{\rho}$ is norm bounded and has the same norm as $\rho$. For $a \in M(A)$, we put $\rho(a) = \overline{\rho}(a)$.

Given two strict linear mappings $\rho : A \rightarrow M(B)$ and $\eta : B \rightarrow M(C)$, we define a new strict linear map $\eta \, \rho : A \rightarrow M(C)$ by  $\eta \, \rho = \overline{\eta} \circ \rho$.
The two basic examples of strict linear mappings are
\begin{itemize}
\item Continuous linear functionals on a \cst-algebra.
\item Non-degenerate $^*$-homomorphism. Recall that a $^*$-homomorphism  $\pi : A \rightarrow M(B)$ is called non-degenerate $\Leftrightarrow$ $ B = [ \, \pi(a)\,b \mid a \in A, b \in B\,]$.
\end{itemize}

\medskip

For $\om \in A^*$ and $a \in M(A)$, we define new elements $a \, \om$ and $\om \, a$ belonging to $A^*$ such that $(a\,\om)(x) = \om(x\,a)$ and
$(\om\,a)(x) = \om(a\,x)$ for $x \in A$.

We also define a  functional $\overline{\om} \in A^*$ such that
$\overline{\om}(x) = \overline{\om(x^*)}$ for all $x \in A$.

\medskip

In this paper, we will also use the notion of a Hilbert \cst-module over a \cst-algebra $A$. For an excellent treatment of Hilbert \cst-modules, we refer to \cite{Lan}.

If $E$ and $F$ are Hilbert \cst-modules over the same \cst-algebra, $\cL(E,F)$ denotes the set of adjointable operators from $E$ into $F$. When $A$ is a \cst-algebra and $H$ is a Hilbert space, $A \ot H$ will denote the Hilbert space over $A$, which is a Hilbert \cst-module
over $A$.

\medskip

Let $H$ be a Hilbert space. The space of bounded operators on $H$ will be denoted by $B(H)$, the space of compact operators on $H$ by $B_0(H)$. Notice that $M(B_0(H)) = B(H)$.

\medskip

Consider a \cst-algebra $A$ and a mapping $\al : \R \rightarrow \text{Aut}(A)$ (where $\text{Aut}(A)$ is the set of $^*$-automorphisms of $A$) such that
\begin{enumerate}
\item $\al_s \, \al_t = \al_{s+t}$ for all $t \in \R$.
\item We have for all $a \in A$ that the function $\R \rightarrow A : t
\rightarrow \al_t(a)$ is norm continuous.
\end{enumerate}
Then we call $\al$ a norm continuous one-parameter group on $A$.
It is then easy to prove that the mapping $\R \rightarrow M(A) : t \mapsto \al_t(a)$ is strictly continuous.

\smallskip

There is a standard way to define for every $z \in \C$ a closed
densely defined linear multiplicative operator $\al_z$ in $A$:
\begin{itemize}
\item The domain of $\al_z$ is by definition the set of elements
$ x \in A$ such that there exists a function $f$ from $S(z)$ into $A$
satisfying
\begin{enumerate}
\item $f$ is continuous on $S(z)$
\item $f$ is analytic on $S(z)^0$
\item We have  that $\al_t(x) = f(t)$ for every $t \in \R$
\end{enumerate}
\item Consider $x$ in the domain of $\al_z$ and $f$ the unique function from $S(z)$ into $A$ such that
\begin{enumerate}
\item $f$ is continuous on $S(z)$
\item $f$ is analytic on $S(z)^0$
\item We have  that $\al_t(x) = f(t)$ for every $t \in \R$
\end{enumerate}
Then we have by definition that $\al_z(x) = f(z)$.
\end{itemize}
where $S(z)$ denotes the strip $\{\,y \in \C \mid \text{Im}\,y \in [0,\text{Im}\,z] \, \}$

\medskip

The mapping $\al_z$ is closable for the strict topology in $M(A)$ and we define the strict closure of $\al_z$ in $M(A)$ by $\overline{\al}_z$. For $a \in D(\overline{\al}_z)$, we put $\al_z(a) := \overline{\al}_z(a)$.

Using the strict topology on $M(A)$, $\overline{\al}_z$ can be constructed from the mapping $\R \rightarrow \text{Aut}(M(A)) : t \rightarrow \overline{\al}_t$ in a similar way as  $\al_z$ is constructed from $\al$.
(See \cite{JK3} or \cite{Ant}, where they used the results in \cite{Zsido} to prove more general results.)

\medskip

If $M$ is a von Neumann algebra, then a strongly continuous one-parameter group (and its extension to the complex plane) is defined in a similar way as a norm continuous one-parameter group but you have to replace the norm topology  by the strong topology.

\bigskip\medskip

\section{Weight theory on \cst-algebras}
\label{weights}

In this section, we will collect some necessary information and conventions about weights. All weights in this paper will be assumed to be non-zero, densely defined and lower semi-continuous. These weights will be called proper weights.

\subsection{Weights on \cst-algebras}

In this first section, we give some information about weights. The standard reference for lower semi-continuous weights is \cite{Comb}. A substantial number of results are collected in \cite{JK1}. We start off with some standard notions concerning lower semi-continuous weights.

\medskip

\begin{definition} \label{weight.def1}
Consider a $C^*$-algebra $A$ and a function $\vfi : A^+ \rightarrow [0,\infty]$  such that
\begin{enumerate}
\item $\vfi(x + y) = \vfi(x) + \vfi(y)$ for all $x,y \in A^+$
\item $\vfi(r x) = r \vfi(x)$ for all $r \in \R^+$ and $x \in A^+$
\end{enumerate}
Then we call $\vfi$ a weight on $A$.
\end{definition}

Let $\vfi$ be a weight on a \cst-algebra $A$. We will use the following standard notations:
\begin{itemize}
\item $\Mfi^+  = \{\, a \in A^+ \mid \varphi(a) < \infty  \,\} $
\item $\Nfi = \{\, a \in A \mid \varphi(a^*a) < \infty \,\} $
\item $\Mfi = \text{span\ } \Mfi^+ = \Nfi^* \Nfi$ .
\end{itemize}
where $\Nfi^* \Nfi = \text{span\ }\{\,y^* x \mid x,y \in \Nfi \, \}$.

\medskip

Condition 1 in definition \ref{weight.def1} implies that $\vfi(x) \leq
\vfi(y)$ for all $x,y \in A^+$ such that $x \leq y$. Therefore
$\Mfi^+$ is a hereditary cone in $A^+$ and $\Nfi$ is a left ideal in
$M(A)$. So we see that $\Mfi \subseteq \Nfi$. We also have that $\Mfi$
is a sub-$^*$-algebra of $A$ and $\Mfi^+ = \Mfi \cap A^+$.

\medskip

It is not so difficult to see that there exists a unique linear map $\psi : \Mfi  \rightarrow \C$ such that $\psi(x) = \vfi(x)$ for all $x \in \Mfi^+$. For every $x \in \Mfi$, we put $\vfi(x) = \psi(x)$.

\medskip

Also the following terminology is standard:
\begin{itemize}
\item We say that $\vfi$ is densely defined $\Leftrightarrow$ $\Mfi^+$ is dense in $A^+$ $\Leftrightarrow$ $\Mfi$ is dense in $A$ $\Leftrightarrow$ $\Nfi$ is dense in $A$.
\item The weight $\vfi$ is called faithful $\Leftrightarrow$ $\bigl( \, \forall a \in A^+ : \vfi(a) = 0 \Rightarrow a = 0 \, \bigr)$.
\end{itemize}

\bigskip

The role of the L$^2$-space of a measure is taken over by the GNS-construction for a weight:

\begin{definition}
Consider a weight $\vfi$ on a \cst-algebra $A$. A GNS-construction for $\vfi$ is by definition a triple
$(H_\vfi,\pifi,\lafi)$ such that
\begin{itemize}
\item $H_\vfi$ is a Hilbert space
\item $\lafi$ is a linear map from $\Nfi$ into
$H_\vfi$ such that
\begin{enumerate}
\item  $\lafi(\Nfi)$ is dense in $H_\vfi$
\item   We have for  every $a,b \in \cN_\vfi$, that
$\langle \lafi(a),\lafi(b) \rangle
= \vfi(b^*a) $
\end{enumerate}
\item $\pifi$ is a representation of $A$ on
$H_\vfi$ such that $\pifi(a)\,\lafi(b) = \lafi(ab)$ for
every $a \in A$ and $b \in \Nfi$.
\end{itemize}
\end{definition}

It is not difficult to construct such a GNS-construction for any weight (cfr. the GNS-construction for a positive functional) and it is unique up to a unitary transformation.

When we use one of the notations $H_\vfi$, $\pifi$ or $\lafi$ without further comment, we implicitly fixed a GNS-construction for $\vfi$.

\bigskip

The following sets play a central role in the theory of lower
semi-continuous weights.

\begin{definition} \label{sweight.def1}
Consider a weight $\vfi$ on a \cst-algebra $A$. Then we define the sets $$\cF_\vfi = \{\, \omega \in A^*_+ \mid
\omega(x) \leq \varphi(x) \text{ for } x \in A^+ \,\}$$ and
$$\cG_\vfi = \{\, \alpha\,\omega \mid
\omega \in \cF_\vfi \, , \, \alpha \in  \,\, ]0,1[ \,\,\} \ \subseteq \cF_\vfi \  . $$
\end{definition}

On $\cF_\vfi$, we use the order inherited from the natural order on $A^*_+$. The advantage of $\cG_\vfi$ over $\cF_\vfi$ lies in the fact that $\cG_\vfi$
is a directed subset of $\cF_\vfi$: for every $\om_1,\om_2 \in \cG_\vfi$, there exists an element $\om \in \cG_\vfi$ such that $\om_1,\om_2 \leq \om$.
This implies that $\cG_\vfi$ can be used as the index set of a net.

 A proof of this fact can be found in \cite{Q-V}. We also included a proof in section 3 of \cite{JK1}.

\medskip

\begin{notation} \label{weight.not1}
Consider a weight $\vfi$ on a \cst-algebra $A$ and a GNS-construction $(H_\vfi,\pifi,\lafi)$ of $\vfi$. Let $\om \in \cF_\vfi$.
\begin{itemize}
\item We define $T_\om$ as the element in $B(H_\vfi) \cap \pifi(A)'$ with $0 \leq T_\om \leq 1$ such that $\langle T_\om \lafi(a) , \lafi(b) \rangle = \om(b^* a)$ for $a,b \in \Nfi$.
\item There exists a unique element $\xi_\om \in H_\vfi$ such that
$T_\om^{\frac{1}{2}} \lafi(a) = \pi_\vfi(a) \xi_\om$ for $a \in \Nfi$.
\end{itemize}
\end{notation}

\bigskip

In order to make weights manageable, we have to impose a continuity condition on them. It turns out that the usual lower semi-continuity is a useful continuity
condition.

\begin{definition}
Consider a  weight $\vfi$ on a \cst-algebra $A$. Then $\vfi$ is lower semi-continuous
\begin{description}
\item[$\Leftrightarrow$] We have for every $\lambda \in \R^+$ that the set
$\{ \, a \in A^+ \mid \vfi(a) \leq \lambda \, \}$ is closed.
\item[$\Leftrightarrow$] If $(x_i)_{i \in I}$ is a net in $A^+$ and $x \in A^+$ such that $(x_i)_{i \in I} \rightarrow x$, then
 $\vfi(x) \leq \liminf \bigl(\vfi(x_i)\bigr)_{i \in I}$.
\end{description}
\end{definition}

Notice that the last condition resembles the result in the classical lemma of Fatou. It implies also easily the next dominated convergence property:

\medskip

Consider $x \in A^+$ and $(x_i)_{i \in I}$ a net in $A^+$ such that $x_i \leq x$ for $i \in I$ and $(x_i)_{i \in I} \rightarrow x$. Then the net $\bigl(\vfi(x_i)\bigr)_{i \in I}$ converges to $\vfi(x)$.

\bigskip

The most important result concerning lower semi-continuous weights is the following one (proven in \cite{Comb} by F. Combes).

\begin{theorem} \label{weight.thm1}
Consider a lower semi-continuous weight $\vfi$ on a \cst-algebra $A$. Then
we have for every $x \in A^+$ that
$$ \vfi(x) = \sup \{\, \omega(x) \mid \om \in \cF_\vfi \,\} \ . $$
\end{theorem}

By writing any element of $\Mfi$ as a sum of elements in $\Mfi^+$,
 we get immediately that the net $\bigl(\om(x)\bigr)_{\om \in \cG_\vfi}$ converges to $\vfi(x)$ for every $x \in \Mfi$.

\medskip

Using this theorem, it is not hard to prove the following properties about a GNS-construction for a lower semi-continuous weight.

\begin{proposition} \label{weight.prop5}
Consider a lower semi-continuous weight $\vfi$ on a \cst-algebra $A$ and a GNS-construction $(H_\vfi,\pifi,\lafi)$ for $\vfi$. Then
\begin{itemize}
\item The mapping $\lafi : \Nfi \rightarrow H_\vfi$ is closed.
\item The $^*$-homomorphism $\pifi : A \rightarrow B(H_\vfi)$ is  non-degenerate.
\item The net $(T_\om)_{\om \in \cG_\vfi}$ converges strongly to 1.
\end{itemize}
\end{proposition}

See e.g. result 2.3 of \cite{JK1} for a proof of the closedness of $\la$.

\medskip

From now on, we will only work with proper weights, i.e. weights which
are non-zero, densely defined and lower semi-continuous.

\bigskip

\subsection{Extensions of lower semi-continuous weights to the multiplier algebra}

Consider a \cst-algebra $A$. Recall that every $\om \in A^*$ has a unique extension $\overline{\om}$ to $M(A)$ which
is strictly continuous and we put $\om(x) = \overline{\om}(x)$ for every $x \in M(A)$.

This implies immediately that any  proper weight has a natural extension to a weight on $M(A)$.

\begin{definition} \label{weight1.def3}
Consider a proper weight $\vfi$ on a \cst-algebra $A$. Then we define the weight $\overline{\vfi}$ on $M(A)$ such that $$\overline{\vfi}(x)=
 \sup \{\, \omega(x) \mid \om \in \cF_\vfi \,\} \ . $$
for every $x \in M(A)^+$. Then $\overline{\vfi}$ is an extension of $\vfi$ and we put $\vfi(x) = \overline{\vfi}(x)$ for all $x \in M(A)^+$.
\end{definition}

We will use the following notations: ${\bar{\cM}}_\vfi^+ ={\cM}_{\overline{\vfi}}^+$, \ ${\bar{\cM}}_\vfi = {\cM}_{\overline{\vfi}}$ and
${\bar{\cN}}_\vfi = {\cN}_{\overline{\vfi}}$.

For any
$x \in {\bar{\cM}}_\vfi$, we put $\vfi(x) = \overline{\vfi}(x)$. It is then clear that the net
$\bigl(\om(x)\bigr)_{\om \in \cG_\vfi}$converges to $\vfi(x)$.

\bigskip

The GNS-construction for a proper weight has a natural extension to a GNS-construction for its extension to the multiplier algebra (see e.g. definition 2.5 and proposition 2.6 of \cite{JK1}).

\begin{proposition}  \label{weight.prop11}
Consider a proper weight $\vfi$ on a \cst-algebra $A$ and a GNS-construction $(H_\vfi,\pifi,\lafi)$ for $\vfi$. Then the mapping $\lafi : \Nfi \rightarrow H_\vfi$ is closable for the strict topology on $M(A)$ and
the norm topology on $H_\vfi$. We denote its closure by $\overline{\la}_\vfi$. Then
$(H_\vfi,\overline{\la}_\vfi,\overline{\pi}_\vfi)$ is a GNS-construction for
$\overline{\vfi}$.
\end{proposition}

In particular, we have that $D(\overline{\la}_\vfi) = \bar{\cN}_\vfi$ and we put
$\lafi(a) = \overline{\la}_\vfi(a)$ for every $a \in \bar{\cN}_\vfi$.

\medskip

Consider $a \in \bar{\cN}_\vfi$. Then there exists a net $(a_i)_{i \in I}$ in $\Nfi$ such that
\begin{itemize}
\item $\|a_i\| \leq \|a\|$ for $i \in I$.
\item $(a_i)_{i \in I}$ converges strictly to $a$
\item $(\lafi(a_i))_{i \in I}$ converges to $\lafi(a)$.
\end{itemize}
This follows immediately by taking an approximate unit in $A$ and multiplying each element of the approximate unit by $a$ from the right.

\bigskip

Let $\om$ be a functional in $\cF_\vfi$. Then it is easy to check that, using the definitions of notation \ref{weight.not1}, the following holds:
\begin{itemize}
\item $\langle T_\om \lafi(a) , \lafi(b) \rangle = \om(b^* a)$ for $a,b \in \bar{\cN}_{\vfi}$.
\item $T_\om^{\frac{1}{2}} \lafi(a) = \pifi(a) \xi_\om$ for $a \in \bar{\cN}_\vfi$.
\end{itemize}

\bigskip

\subsection{KMS weights on a \cst-algebra}

Although a \cst-algebra is generally non-commutative, we would like to
have some control over the non-commutativity under the weight.
Therefore we will introduce the class of KMS weights. For full
details, we refer to \cite{JK1}.

\begin{definition} \label{sweight.def2}
Consider a \cst-algebra $A$ and a weight $\vfi$ on $A$. We say that $\vfi$ is a KMS weight
on A  $\Leftrightarrow$ $\vfi$ is a proper weight on $A$ and  there exists  a norm continuous one-parameter group $\si$ on $A$ satisfying the
following properties:
\begin{enumerate}
\item $\vfi$ is invariant under $\si$: $\vfi \, \si_t = \vfi$ for every $t \in \R$.
\item We have for every $a \in D(\si_\frac{i}{2})$ that
      $\vfi(a^* a) = \vfi(\si_\frac{i}{2}(a) \si_\frac{i}{2}(a)^*)$.
\end{enumerate}
The one-parameter group $\si$ is called a modular group for $\vfi$.
\end{definition}

If  the weight $\vfi$ is faithful, then the one-parameter group $\si$ is uniquely determined and is called the modular group of $\vfi$.

This is not the usual definition of a KMS weight on a \cst-algebra (see \cite{Comb1}), but we prove in \cite{JK1} that this definition is
equivalent with the usual one. More precisely,

\begin{proposition}
Consider a proper weight $\vfi$ on a \cst-algebra $A$ with GNS-construction $(H_\vfi,\pifi,\lafi)$. Let $\si$ be a norm continuous one-parameter group on $A$ such that $\vfi \, \si_t = \vfi$ for all $t \in \R$. Then the following conditions are equivalent
\begin{enumerate}
\item We have that $\vfi(a^* a) = \vfi(\si_{\frac{i}{2}}(a) \, \si_{\frac{i}{2}}(a)^*)$ for all $a \in D(\si_{\frac{i}{2}})$.
\item There exists a non-degenerate $^*$-antihomomorphism $\th : A \rightarrow B(H_\vfi)$ such that we have for all $x \in \Nfi$ and $a \in D(\si_{\frac{i}{2}})$ that $x a$ belongs to $\Nfi$ and
$\lafi(x a) = \th(\si_{\frac{i}{2}}(a)) \, \lafi(x)$.
\item For all $a,b \in \Nfi \cap \Nfi^*$,  there exists a function $f : S(i) \rightarrow \C$ such that
\begin{itemize}
\item $f$ is continuous and bounded on $S(i)$
\item $f$ is analytic on $S(i)^{\circ}$
\item $f(t) = \vfi(\si_t(b)\, a)$ and $f(t+i) = \vfi(a\,\si_t(b))$ for $t \in \R$.
\end{itemize}
\end{enumerate}
\end{proposition}

\bigskip

\subsection{Basic properties of KMS weights}

In this subsection, we collect the most common properties of KMS weights.
To this end, we will fix a \cst-algebra $A$ and a KMS weight $\vfi$ on $A$ with modular group $\si$. Let $(H_\vfi,\pifi,\lafi)$ be a GNS-construction
for $\vfi$.

\bigskip

In the next proposition,  we formulate some basic properties of KMS weights.

\begin{proposition} \label{weight.prop6}
Let $\vfi$ be a KMS weight on a \cst-algebra $A$, with GNS-construction $(H_\vfi,\pifi,\lafi)$.
Then the following properties hold:
\begin{enumerate}
\item There exists a unique anti-unitary operator $J$ on $H_\vfi$ such
      that $J \lafi(x) = \lafi(\si_\frac{i}{2}(x)^*)$ for every $x \in
      \Nfi \cap D(\si_\frac{i}{2})$.
\item Let $a \in D(\si_\frac{i}{2})$ and $x \in \Nfi$. Then $x a$
belongs to $\Nfi$ and $\lafi(x a) = J \pi_\vfi(\si_\frac{i}{2}(a))^* J \, \lafi(x)$.
\item Let $a \in D(\si_{-i})$ and $x \in \Mfi$. Then $a x$ and $x \si_{-i}(a)$
belong to $\Mfi$ and $\vfi(a x) = \vfi(x \si_{-i}(a))$.
\item Consider $x \in \Nfi \cap \Nfi^*$ and $a \in \Nfi^* \cap D(\si_{-i})$ such that $\si_{-i}(a) \in \Nfi$. Then  $\vfi(a x) = \vfi(x \si_{-i}(a))$.
\end{enumerate}
\end{proposition}

\medskip

The anti-unitary operator $J$ will be called the modular conjugation  of $\vfi$ in the GNS-construction $(H_\vfi,\pifi,\lafi)$. We also have  a strictly positive operator $\nab$ in $H_\vfi$ such that
$\nab^{it} \lafi(a) = \lafi(\si_t(a))$ for $t \in \R$ and $a \in \Nfi$.
The operator $\nab$ will be  called the modular operator of $\vfi$ in the GNS-construction $(H_\vfi,\pifi,\lafi)$.

\bigskip

Although the definition of the modular conjugation and the modular operator depend on $\si$, they only depend on the weight $\vfi$:

There exists a densely defined closed operator $T$ from within $H_\vfi$ into $H_\vfi$ such that $\lafi(\Nfi \cap \Nfi^*)$ is a core for $T$ and
$T \lafi(a) = \lafi(a^*)$ for $a \in \Nfi \cap \Nfi^*$.

Then $\nab = T^* T$ and $T = J \nab^{\frac{1}{2}} = \nab^{-\frac{1}{2}} J$. Also notice that $J \nab^t J = \nab^{-t}$ and  $J \nab^{it} J = \nab^{it}$ for $t \in \R$.

\bigskip

The above proposition can be easily extended to elements in the multiplier algebra by using the extensions $\overline{\vfi}$ and $\overline{\si}$. It boils down to approximate elements in the multiplier algebra by elements in the \cst-algebra $A$ in a nice way.

\begin{proposition} \label{extmult}
The following properties hold
\begin{enumerate}
\item Let $a \in D(\overline{\si}_\frac{i}{2})$ and $x \in \bar{\cN}_\vfi$. Then $x a$
belongs to $\bar{\cN}_\vfi$ and $\lafi(x a) = J \pi_\vfi(\si_\frac{i}{2}(a))^* J \, \lafi(x)$.
\item We have for all $x \in D(\overline{\si}_{\frac{i}{2}})$ that
$\vfi(x^* x) = \vfi(\si_{\frac{i}{2}}(x) \, \si_{\frac{i}{2}}(x)^*)$.
\item We have that $J \lafi(x) = \lafi(\si_\frac{i}{2}(x)^*)$ for every $x \in \bar{\cN}_\vfi \cap D(\overline{\si}_\frac{i}{2})$.
\item Let $a \in D(\overline{\si}_{-i})$ and $x \in \bar{\cM}_\vfi$. Then $a x$ and $x \si_{-i}(a)$
belong to $\bar{\cM}_\vfi$ and $\vfi(a x) = \vfi(x \si_{-i}(a))$.
\item Consider $x \in \bar{\cN}_\vfi \cap \bar{\cN}_\vfi^*$ and $a \in \bar{\cN}_\vfi^* \cap D(\overline{\si}_{-i})$ such that $\si_{-i}(a) \in \bar{\cN}_\vfi$. Then  $\vfi(a x) = \vfi(x \si_{-i}(a))$.
\end{enumerate}
\end{proposition}

\bigskip\medskip

If we have  a proper weight $\eta$ which agrees with $\vfi$  on the intersection $\Mfi^+ \cap \cM_\eta^+$ and such that the proper weight $\eta$ is invariant under a modular group of $\vfi$, then $\vfi = \eta$. This will follow easily once we have proven the next proposition.

\begin{proposition} \label{sweight.prop1}
Consider a proper weight $\eta$ on $A$ with GNS-construction $(H_\eta,\pi_\eta,\la_\eta)$ and such that there exists a strictly positive number $\lambda > 0$ such that
$\eta \, \si_t = \lambda^t \, \eta$ for $t \in \R$. Then $\Nfi \cap \cN_\eta$ is a core for both $\lafi$ and $\la_\eta$.
\end{proposition}
\begin{demo} We define  the injective positive operator $T$ in $H_\eta$ such that $T^{it} \la_\eta(a) = \lambda^{-\frac{t}{2}} \, \la_\eta(\si_t(a))$ for $t \in \R$.
\begin{enumerate}
\item Choose $x \in \cN_\eta$. For every $n \in \N$, we define $x_n \in A$ such that
$$ x_n = \frac{n}{\sqrt{\pi}} \int \exp(-n^2 t^2) \, \si_t(x) \, dt \ , $$
then $x_n$ belongs to $\cN_\eta \cap D(\si_{\frac{i}{2}})$ and
$$\la_\eta(x_n) = \frac{n}{\sqrt{\pi}} \int \exp(-n^2 t^2) \, \lambda^{\frac{t}{2}} \, T^{it} \la_\eta(x) \, dt \ .$$
So we get that $(x_n)_{n=1}^\infty$ converges to $x$ and  that
$(\la_\eta(x_n))_{n=1}^\infty$ converges to $\la_\eta(x)$.

It is well known that there exists a bounded net $(e_k)_{k \in K}$ in $\Nfi$ such
that $(e_k)_{k \in K}$ converges strictly to 1. We have for every $n \in \N$ and $k \in K$ that $e_k \, x_n$ belongs to $\Nfi \cap \cN_\eta$ (because $x_n$ belongs to $D(\si_{\frac{i}{2}})$ and $e_k$ belongs to $\Nfi$).

Clearly, $(e_k \, x_n)_{(n,k) \in \N \times K}$ converges to $x$. Since $\la_\eta(e_k \, x_n) = \pi_\eta(e_k) \la_\eta(x_n)$ for all $k \in K$ and $n \in \N$, we also see that the net
$(\la_\eta(e_k \, x_n))_{(n,k) \in \N \times K}$ converges to $\la_\eta(x)$.

So we have proven that $\Nfi \cap \cN_\eta$ is a core for $\la_\eta$.

\item Now choose $x \in \Nfi$. There exists a bounded net $(u_l)_{l \in L}$ in $\cN_\eta$ which converges strictly to 1.
For every $l \in L$, we define the element $v_l \in A$ such that
$$v_l = \frac{n}{\sqrt{\pi}} \int \exp(-n^2 t^2) \, \si_t(u_l) \, dt$$
which is, as above, an element in $\cN_\eta \cap D(\si_{\frac{i}{2}})$ such that
$$\si_{\frac{i}{2}}(v_l) = \frac{n}{\sqrt{\pi}} \int \exp(-n^2 (t-\frac{i}{2})^2\,) \, \si_t(u_l) \, dt \ .$$
These two formulas imply that the nets $(v_l)_{l \in L}$ and  $(\si_{\frac{i}{2}}(v_l))_{l \in L}$ are bounded and converge strictly to 1 (see result \ref{int.res1}).

We have for all $l \in L$ that $x\, v_l \in \Nfi \cap \cN_\eta$ and
$\lafi(x \,v_l) = J \pifi(\si_{\frac{i}{2}}(v_l))^* J \lafi(x)$. So we get that $(x \, v_l)_{l \in L}$ converges to $x$ and that $(\lafi(x \, v_l))_{l \in L}$ converges to $\lafi(x)$.

So we have proven that $\Nfi \cap \cN_\eta$ is a core for $\lafi$.
\end{enumerate}
\end{demo}

\medskip

This implies easily the following desired result.

\begin{corollary} \label{sweight.cor1}
Consider a  proper weight $\eta$ on $A$ such that
$\eta \, \si_t = \eta$ for $t \in \R$ and $\eta(x) = \vfi(x)$ for all $x \in \Mfi^+ \cap \cM_\eta^+$. Then $\eta = \vfi$.
\end{corollary}
\begin{demo}
Take a GNS-construction $(H_\eta,\pi_\eta,\la_\eta)$ for $\eta$.
Then we have clearly for all $x \in \Nfi \cap \cN_\eta$ that
$$\|\lafi(x)\|^2 = \vfi(x^* x) = \eta(x^* x) = \|\la_\eta(x)\|^2 \ .$$
This implies for every net $(x_i)_{i \in I}$ in $\Nfi \cap \cN_\eta$ that $(\lafi(x_i))_{i \in I}$ is Cauchy $\Leftrightarrow$ the net $(\la_\eta(x_i))_{i \in I}$ is Cauchy and hence that $(\lafi(x_i))_{i \in I}$ is convergent $\Leftrightarrow$ the net $(\la_\eta(x_i))_{i \in I}$ is convergent.

Combining this with the fact that $\Nfi \cap \cN_\eta$ is a core for both $\lafi$ and $\la_\eta$ and the closedness of $\lafi$,$\la_\eta$, we get that
$\Nfi = \cN_\eta$ and that $\|\lafi(x)\|^2 = \|\la_\eta(x)\|^2 $ for all $x \in \Nfi$.
\end{demo}

\bigskip\medskip

\section{Extending lower semi-continuous weights to normal weights}

Even in a \cst-algebraic approach to quantum groups it is (probably) unavoidable to get into the \wst-framework at some stages. The main reason for this is the existence of a Radon-Nikodym theorem for well behaved weights on von Neumann algebras. Another reason is the fact that a normal semi-finite faithful weight satisfies automatically some sort of KMS condition. In order to use the von Neumann algebra setting, we will need a well-behaved procedure to  extend a lower semi-continuous weight $\vfi$ on a \cst-algebra to the von Neumann algebra generated by this \cst-algebra
in the GNS-construction of $\vfi$.

\bigskip

\subsection{Normal weights on  a von Neumann algebra}

In this subsection, we will recall some notions concerning normal weights and state the main properties. If we consider weights on a  von Neumann algebra, it is more convenient to replace the norm topology by the $\si$-weak or $\si$-strong$^*$ topology.

\medskip

Let us first introduce the following terminology. Consider a  von Neumann algebra $M$ and a weight $\vfi$ on $M$. We call $\vfi$ semi-finite $\Leftrightarrow$ $\Mfi^+$ is $\si$-weak dense in $M^+$ $\Leftrightarrow$ $\Mfi$ is $\si$-weak dense in $M$ $\Leftrightarrow$ $\Nfi$ is $\si$-weak dense in $M$.

\medskip

Because the $\si$-strong$^*$ and the $\si$-weak topology have the same continuous linear functionals, the closure of a convex subset of $M$ agrees for both topologies. Suppose now that $\vfi$ is semi-finite. Then this argument shows that $\Mfi^+$ is $\si$-strong dense in $M^+$ and that $\Mfi$ and $\Nfi$ are $\si$-strong$^*$ dense in $M$.

\medskip

It is easy to improve this a little bit. Recall that the $\si$-strong$^*$ topology and the strong$^*$-topology agree on bounded subsets of $M$.

Because $\Mfi^+ = \Mfi \cap M^+$, Kaplansky's density theorem implies for every $x \in M^+$ the existence of a net $(x_i)_{i \in I}$ in $\Mfi^+$ such that $(x_i)_{i \in I}$ converges strongly to $x$ and $\|x_i\| \leq \|x\|$ for $i \in I$.

\medskip

Combining this with the polar decomposition, we get for all $x \in M$ the existence of a net $(x_i)_{i \in I}$ in $\Nfi$ such that $(x_i)_{i \in I}$ converges strongly$^*$ to $x$ and $\|x_i\| \leq \|x\|$ for $i \in I$.

\medskip

We will always work with semi-finite weights.

\bigskip\medskip

Let us define the obvious von Neumann algebra variant of the sets $\cF_\vfi$ and $\cG_\vfi$ for weights on a \cst-algebra (see definition \ref{sweight.def1}).

\begin{definition}
Consider a von Neumann algebra $M$ and a semi-finite weight $\vfi$ on $M$. Then we define the sets
$$\tilde{\cF}_\vfi = \{\, \om \in M_*^+ \mid \om(x) \leq \vfi(x) \text{ for all } x \in M^+ \, \} $$
and
$$\tilde{\cG}_\vfi = \{\, \al \, \om  \mid  \om \in \tilde{\cF}_\vfi, \al \in [0,1[ \, \} \ .$$
\end{definition}

\medskip

Again, the set $\tilde{\cG}_\vfi$ is upwardly directed with respect to the natural order. The proof of this fact for $\cG_\vfi$ in \cite{Q-V} remains valid (see also section 3 of \cite{JK1}).
If $\om_1, \om_2 \in \tilde{\cG}_\vfi$,  a functional $\om \in \cG_\vfi$ is constructed such that $\om_1,\om_2 \leq \om$. But the construction in \cite{Q-V} is such that there exist $\lambda_1, \lambda_2 \in \R^+$ such that
$\om \leq \lambda_1 \, \om_1 + \lambda_2 \, \om_2$. Because $\lambda_1 \, \om_1 + \lambda_2 \, \om_2$ is normal, this implies that $\om$ is normal (use Corollary~III.3.11 in \cite{Tak} for instance). Hence $\om \in \tilde{\cG}_\vfi$.

\bigskip

The more  natural condition in the setting of von Neumann algebras is the lower semi-continuity with respect to the $\si$-weak topology:

\begin{definition}
Consider a von Neumann algebra $M$ and a semi-finite weight $\vfi$ on $M$. Then we call $\vfi$ normal $\Leftrightarrow$ $\vfi$ is lower semi-continuous with respect to the $\si$-weak topology.
\end{definition}

\medskip

In \cite{Haa} U. Haagerup formulates some other conditions which he proves to be equivalent to the normality of the weight $\vfi$. One of these equivalences is the following rather appealing one.

\medskip

The weight $\vfi$ is normal $\Leftrightarrow$ We have for every increasing norm bounded net $(x_i)_{i \in I}$ in $M^+$  that \newline $\vfi\bigl(\sup\,(x_i)_{i \in I}\bigr)$ $= \sup\,(\vfi(x_i))_{i \in I}$.

\medskip\medskip

Another equivalent condition in \cite{Haa} is the content of the following fundamental result in the theory of normal weights. It is the normal variant of theorem \ref{weight.thm1}.

\begin{theorem}
Consider a von Neumann algebra $M$ and a normal semi-finite weight $\vfi$ on $M$. Then
$$\vfi(x) = \sup \, \{\,\om(x) \mid \om \in \tilde{\cF}_\vfi \, \} $$
for all $x \in M^+$.
\end{theorem}

\medskip

Because the $\si$-weakly continuous and $\si$-strong$^*$ continuous linear functionals are the same, this result implies that a normal weight is also $\si$-strong lower semi-continuous.

\medskip

Using the previous theorem,  theorem 7.2 of \cite{Pe-Tak} implies the following proposition.

\begin{proposition}
Consider a von Neumann algebra $M$ acting on a Hilbert space $H$ and a normal semi-finite weight $\vfi$ on $M$. Then there exists a family of vectors $(v_i)_{i \in I}$ in $H$ such that
$$\vfi(x) = \sum_{i \in I} \, \om_{v_i,v_i}(x) $$
for all $x \in M^+$.
\end{proposition}

\medskip

For normal weights, proposition \ref{weight.prop5} can be modified in the following way.

\begin{proposition}
Consider a von Neumann algebra $M$ and a normal semi-finite weight $\vfi$ on $M$. Let $(H_\vfi,\pifi,\lafi)$ be a GNS-construction for $\vfi$. Then
\begin{enumerate}
\item $\pifi : M \rightarrow B(H_\vfi)$ is a unital normal $^*$-homomorphism and $\pifi(M)$ is a von Neumann algebra acting on $H_\vfi$.
\item The mapping $\lafi : \Nfi \rightarrow H_\vfi$ is closed for both the $\si$-weak and the $\si$-strong$^*$ topology.
\item The net $(T_\om)_{\om \in \tilde{\cG}_\vfi}$ converges strongly to $1$.
\end{enumerate}
\end{proposition}

For a proof of the first result, we refer to Theorem~2.2 in \cite{Stramod}. The proof of the second one  is then similar to the proof of result 2.3 in \cite{JK1}.

\bigskip\medskip

A central role in the theory of normal weights is played by the faithful ones. Using left Hilbert algebra techniques, it is possible to show that a faithful normal weight automatically satisfies some sort of KMS property, provided we relax the continuity requirements on the modular group (see e.g. \cite{Stra}). In contrast with definition \ref{sweight.def2}, we will stick to the usual definition of the modular group in this case.

\medskip

A normal faithful semi-finite weight is also abbreviated to a n.f.s
weight.

\begin{theorem}
Consider a von Neumann algebra $M$ and a faithful semi-finite normal weight $\vfi$ on $M$. Then there exists a unique strongly continuous one-parameter group $\si$ on $M$ such that
\begin{itemize}
\item $\vfi \, \si_t = \vfi$ for all $t \in \R$.
\item For every $x,y \in \Nfi \cap \Nfi^*$ , there exists a function $f : S(i) \rightarrow \C$ such that
\begin{enumerate}
\item $f$ is continuous and bounded on $S(i)$
\item $f$ is analytic on $S(i)^{\circ}$
\item $f(t) = \vfi(\si_t(b)\, a)$ and $f(t+i) = \vfi(a\,\si_t(b))$ for $t \in \R$.
\end{enumerate}
\end{itemize}
We call $\si$ the modular group of $\vfi$.
\end{theorem}

\medskip

It is possible to proof the following equivalent defining properties for a modular group. The proof can be easily extracted from the results in \cite{JK1} by a translation from  the \cst-algebra to the von Neumann algebra setting.

\begin{proposition}
Consider a von Neumann algebra $M$,  a faithful semi-finite normal weight  $\vfi$ on $M$ with GNS-construction $(H_\vfi,\lafi,\pifi)$. Let $\si$ be a strongly continuous one-parameter group on $M$ such that $\vfi \, \si_t = \vfi$ for $t \in R$. Then the following conditions are equivalent.
\begin{enumerate}
\item $\si$ is the modular group of $\vfi$
\item $\vfi(a^* a) = \vfi(\,\si_{\frac{i}{2}}(a)\,\si_{\frac{i}{2}}(a)^*\,)$ for $a \in D(\si_{\frac{i}{2}})$
\item There exists a normal $^*$-antihomomorphism $\th : M \rightarrow B(H_\vfi)$ such that we have for all $x \in \Nfi$ and $a \in D(\si_{\frac{i}{2}})$ that
$x\,a$ belongs to $\Nfi$ and $\lafi(x\,a) = \th(\si_{\frac{i}{2}}(a)) \lafi(x)$.
\end{enumerate}
\end{proposition}

\bigskip

\subsection{The extension procedure}

We look at an extension of a lower semi-continuous weight to a normal weight on the von Neumann algebra generated in its GNS-space. We follow the discussion of S. Baaj in \cite{Baa3} and combine it with a technique of J. Verding and A. Van Daele to get to an extremely useful description of the GNS-construction of the normal extension. Most of the results (and their proofs) are taken from \cite{Baa3}. However, the most intricate results have been proven by U. Haagerup (in \cite{Haa}) and F. Combes (in \cite{Comb}).

\begin{definition} \label{weight.def5}
Consider a topological space $X$ with dense subspace $Y$. Let $\psi : Y \rightarrow [0,\infty]$ be a lower semi-continuous function. Then we define a new function $\hat{\psi} : X \rightarrow [0,\infty]$ as follows. Consider an element $x \in X$, then
$$\hat{\psi}(x)  =  \inf \, \{ \, L \in [0,\infty] \mid \exists \text{ a net } (x_i)_{i \in I} \text{ in } Y \text{ such that } (x_i)_{i \in I} \rightarrow x
\text{ and } \psi(x_i) \leq L \text{ for all } i \in I \, \} \ .$$
\end{definition}

\medskip

Then we can easily prove the next proposition.

\begin{proposition} \label{weight.prop18}
Consider a topological space $X$ with dense subspace $Y$ and let $\psi : Y \rightarrow [0,\infty]$ be a lower semi-continuous function. Then the following properties hold:
\begin{enumerate}
\item $\psi(x) = \hat{\psi}(x)$ for all $x \in Y$.
\item $\hat{\psi}$ is lower semi-continuous.
\item Consider a lower semi-continuous function $\eta : X \rightarrow [0,\infty]$ such that $\eta(x) = \psi(x)$ for all $x \in Y$. Then
$\eta(x) \leq \hat{\psi}(x)$ for all $x \in X$.
\end{enumerate}
\end{proposition}
\begin{demo}
\begin{enumerate}
\item Choose $x \in Y$. By taking the constant sequence with element $x$, we get immediately that $\hat{\psi}(x) \leq \psi(x)$. Take any element $L \in [0,\infty]$ such that there exists a net $(x_i)_{i \in I}$ in $Y$ such that
$(x_i)_{i \in I}$ converges to $x$ and $\psi(x_i) \leq L$ for $i \in I$.
By lower semi-continuity of $\psi$, we get that $\psi(x) \leq L$.

Hence $\psi(x) \leq \hat{\psi}(x)$, which gives that $\psi(x) = \hat{\psi}(x)$.
\item Choose $\lambda \in \R^+$. We have to prove that the set $\{ \, y \in X \mid \hat{\psi}(y) \leq \lambda \, \}$ is closed in $X$.

Take $x \in X$ such that there exists a net $(x_i)_{i \in I}$ in $X$ such that $(x_i)_{i \in I}$ converges to $x$ and $\hat{\psi}(x_i) \leq \lambda$ for $i \in I$.
Pick $\vep > 0$. Define $J$ to be the set of open neighbourhoods of $x$ in $X$. We order $J$ by reverse inclusion. This turns $J$ into a directed set.

Choose $U \in J$. Then there exists an element $i_U \in I$ such that
$x_{i_U} \in U$. Because $U$ is an open neighbourhood of $x_{i_U}$ and $\hat{\psi}(x_{i_U}) < \lambda + \vep$, the definition of $\hat{\psi}$ implies the existence of an element $y_U \in U \cap Y$ such that
$\psi(y_U) \leq \lambda + \vep$.

It is now clear that $(y_U)_{U \in J}$ converges to $x$ so the definition of $\hat{\psi}$ implies that $\hat{\psi}(x) \leq \lambda + \vep$.

From this all, we get that $\hat{\psi}(x) \leq \lambda$.
\item Take an element $x \in X$. Let $L$ be an element in $[0,\infty]$ such that there exists a net $(x_i)_{i \in I}$ in $Y$ such that
$\psi(x_i) \leq L$ for all $i \in I$ and $(x_i)_{i \in I} \rightarrow x$.
Then $\eta(x_i) = \psi(x_i) \leq L$ for $i \in I$, so the lower semi-continuity of $\eta$ implies that $\eta(x) \leq L$.
Hence $\eta(x) \leq \hat{\psi}(x)$.
\end{enumerate}
\end{demo}

\bigskip\medskip

For the rest of this section, we fix a \cst-algebra $A$ and a proper weight
$\vfi$ on $A$. Let $(H_\vfi,\pifi,\lafi)$ be a GNS-construction for $\vfi$.

\medskip

By notation \ref{weight.not1}, we have for all $\om \in \cF_\vfi$ that
$\om(x) = \langle \pifi(a) \, \xi_\om , \xi_\om \rangle$ for $a \in A$.
So there exists a unique element $\tilde{\om} \in (\pifi(A)'')_*^+$ such that $\tilde{\om} \, \pifi = \om$, i.e. $\tilde{\om}(x) = \langle x \, \xi_\om , \xi_\om \rangle$ for $x \in \pifi(A)''$.

\medskip

Because $\cG_\vfi$ is upwardly directed, the same is true for the set
$\{ \, \tilde{\om} \mid \om \in \cG_\vfi \, \}$. This will imply that we really get a weight in the next definition.

\begin{definition} \label{weight.def6}
We define the function $\tilde{\vfi} : (\pifi(A)'')^+ \rightarrow [0,\infty]$ such that $\tilde{\vfi}(x) = \sup \{\, \tilde{\om}(x) \mid \om \in \cF_\vfi \, \}$ for $x \in (\pifi(A)'')^+$. Then $\tilde{\vfi}$ is a normal semi-finite weight on $\pifi(A)''$ such that $\tilde{\vfi} \, \pifi = \vfi$.
\end{definition}

Notice that the last statement follows from theorem \ref{weight.thm1}.
We call $\tilde{\vfi}$ the \wst-lift of $\vfi$ in the GNS-construction $(H_\vfi,\pifi,\lafi)$.

\medskip

Using the fact that $\tilde{\vfi} \, \pifi = \vfi$, it is easy to check that $\tilde{\cF}_{\tilde{\vfi}} = \{\, \tilde{\om} \mid \om \in \cF_\vfi \, \}$ and that $\tilde{\cG}_{\tilde{\vfi}} = \{\, \tilde{\om} \mid \om \in \cG_\vfi \, \}$.

\bigskip

Although this definition is conceptually appealing it is not that practical in many applications. It will turn out that it is more useful to get a GNS-construction for $\tilde{\vfi}$ by closing $\lafi$ with respect to the $\si$-strong$^*$-topology. For this, the material in \cite{Baa3} gives us most of the necessary results. We will repeat the reasoning by S. Baaj in this section.

\medskip

Apply definition \ref{weight.def5} with $X = (\pifi(A)'')^+$ with the $\si$-strong topology. We define moreover $Y = \pifi(A)^+$ and $\psi =$ the restriction of $\tilde{\vfi}$ to $Y$.

Then definition \ref{weight.def5} gives us a $\si$-strong lower semi-continuous function $\hat{\psi} :  (\pifi(A)'')^+ \rightarrow [0,\infty]$. In a while, we will see that $\hat{\psi} = \tilde{\vfi}$.

\begin{lemma}  \label{weight.lem21}
Consider $x \in (\pifi(A)'')^+$. Then $\hat{\psi}(x)$ is equal to the infimum of the set
$$\bigl\{ \, L \in [0,\infty] \mid \exists \text{ a net } (x_i)_{i \in I} \text{ in } \pifi(A)^+ \text{ s.t. } (x_i)_{i \in I} \stackrel{\text{{\footnotesize strongly}}}{\longrightarrow} x  \ ,
 \|x_i\| \leq \|x\| \text{ and }\psi(x_i) \leq L \text{ for all } i \in I \, \bigr\}
$$
\end{lemma}
\begin{demo}
Take $L \in [0,\infty]$ such that there exists a net $(y_i)_{i \in I}$ in $\pifi(A)^+$  such that $(y_i)_{i \in I} \rightarrow x$ $\si$-strongly and
$\psi(y_i) \leq L$ for $i \in I$.

There exists clearly a continuous function $f : \R^+ \rightarrow \R^+$ with compact support such that
\begin{itemize}
\item $f(t) = t$ for $t \in [0,\|x\|]$
\item $f(t) \leq t$ for $t \in \R^+$
\item $f(t) \leq \|x\|$ for $t \in \R^+$
\end{itemize}
For every $i \in I$, we put $x_i = f(y_i) \in \pifi(A)^+$. Because $f$ has compact support and $(y_i)_{i \in I}$ converges strongly to $x$, the net
$(x_i)_{i \in I}$ converges strongly to $f(x)$ (see theorem 5.3.4 of \cite{Kad} for instance).
Because $f(t) = t$ for $t \in [0,\|x\|]$, we have that $f(x) = x$. So $(x_i)_{i \in I}$ converges strongly to $x$.
Since $f(t) \leq \|x\|$ for $t \in \R$, we have also immediately that $\|x_i\| \leq \|x\|$ for $i \in I$.

It is also clear for every $i \in I$ that $x_i = f(y_i) \leq y_i$, implying that  $\psi(x_i) \leq \psi(y_i) \leq L$.

\medskip

From this discussion and definition \ref{weight.def5}, the lemma follows now immediately.
\end{demo}

\medskip

\begin{lemma} \label{weight.lem20}
The mapping $\hat{\psi} : (\pifi(A)'')^+ \rightarrow [0,\infty]$ satisfies the following properties:
\begin{enumerate}
\item $\hat{\psi}(r \,x) = r \, \hat{\psi}(x)$ for $x \in (\pifi(A)'')^+$ and $r \in \R^+$.
\item $\hat{\psi}(x+y) \leq \hat{\psi}(x) + \hat{\psi}(y)$ for all $x,y \in
(\pifi(A)'')^+$.
\item $\hat{\psi}(x) \leq  \hat{\psi}(y)$ for all $x,y \in (\pifi(A)'')^+$.
\end{enumerate}
\end{lemma}
\begin{demo}
The first properties follow easily from definition \ref{weight.def5}. We will look into the third one.

So choose $x,y \in (\pifi(A)'')^+$ such that $x \leq y$ and such that $\hat{\psi}(y) < \infty$. Take $\vep > 0$.

By the previous lemma, there exists a bounded net $(y_i)_{i \in I}$ in $\pifi(A)^+$ such that $(y_i)_{i \in I}$ converges strongly to $y$ and
$\psi(y_i) \leq \hat{\psi}(y) + \vep$ for $i \in I$. Then $(y_i^{\frac{1}{2}})_{i \in I}$ is a bounded net in $\pifi(A)$ which converges strongly to $y^{\frac{1}{2}}$.

Because $x \leq y$, there exists an element $u \in \pifi(A)''$ such that
$x^{\frac{1}{2}} = u \, y^{\frac{1}{2}}$ and $\|u\| \leq 1$. By Kaplansky's density theorem, there exists a net $(v_k)_{k \in K}$ in $\pifi(A)^+$  such that $\|v_k\| \leq 1$ for $k \in K$ and $(v_k)_{k \in K}$ converges strongly to $u^* u$.

Then $(y_i^{\frac{1}{2}} v_k \, y_i^{\frac{1}{2}})_{(i,k) \in I \times K}$ is a bounded net in $\pifi(A)^+$ which converges strongly to $y^{\frac{1}{2}} u^* u y^{\frac{1}{2}} = x$.
We have moreover for all $i \in I$ and $k \in K$ that $y_i^{\frac{1}{2}} v_k \, y_i^{\frac{1}{2}} \leq y_i$, implying that $\psi(y_i^{\frac{1}{2}} v_k \, y_i^{\frac{1}{2}}) \leq \psi(y_i) \leq \hat{\psi}(y) + \vep$.

By definition \ref{weight.def5}, this gives  that $\hat{\psi}(x) \leq \hat{\psi}(y) + \vep$. Hence we conclude that $\hat{\psi}(x) \leq \hat{\psi}(y)$.
\end{demo}

\bigskip

In \cite{Haa}, U. Haagerup proved the following fundamental theorem about subadditive positive lower semi-continuous functions on a \wst-algebra.

\begin{theorem} \label{weight.thm2}
Consider a \wst-algebra $B$ and a function $\eta : B^+ \rightarrow [0,\infty]$ such that
\begin{enumerate}
\item $\eta$ is $\si$-strongly lower semi-continuous.
\item We have for all $x,y \in B^+$ that $\eta(x+y) \leq \eta(x) + \eta(y)$.
\item We have for all $x \in B^+$ and $r \in \R^+$ that $\eta(r \, x) = r \, \eta(x)$.
\item Let $x,y \in B^+$. If $x \leq y$, then $\eta(x) \leq \eta(y)$.
\end{enumerate}
and define the set $\cF = \{\,\om \in B_*^+ \mid \om(x) \leq \eta(x) \text{ for all } x \in B^+ \, \}$. Then
$$\eta(x) = \sup \, \{ \,\om(x) \mid \om \in \cF \,\} $$
for all $x \in B^+$.
\end{theorem}

In proposition 2.1 of \cite{Haa}, U. Haagerup shows two other conditions to be equivalent to the conclusions of this theorem. In theorem 2.2 of the same paper he then proves that the selfadjoint part of $B$ together with the $\si$-weak topology satisfies condition (1) of his proposition 2.1. But the $\si$-strong topology and the $\si$-weak topology have the same continuous linear functionals which implies for any convex subset $C$ of $B^+$ that the closures of $C$ in both topologies agree. It is then easy to see that
self adjoint part of $B$ together with the $\si$-strong topology also satisfies condition (1) of his proposition 2.1.

\bigskip

\begin{proposition}
We have that $\hat{\psi} = \tilde{\vfi}$.
\end{proposition}
\begin{demo}
By proposition \ref{weight.prop18}.3, we get immediately that $\tilde{\vfi}(x) \leq \hat{\psi}(x)$ for $x \in (\pifi(A)'')^+$.

By lemma \ref{weight.lem20} and theorem \ref{weight.thm2}, the following holds.
Define
$$\cF = \{\,\om \in (\pifi(A)'')^+_* \mid \om(x) \leq \hat{\psi}(x) \text{ for all } x \in (\pifi(A)'')^+ \, \} \ .$$
Then $\hat{\psi}(x) = \sup \, \{ \, \om(x) \mid \om \in \cF \, \}$ for all $x \in (\pifi(A)'')^+$.

Choose $\om \in \cF$. Put $\th = \om \pifi$. Then we have  for all $a \in A^+$ that
$$\th(a) = \om(\pifi(a)) \leq \hat{\psi}(\pifi(a)) = \psi(\pifi(a)) = \tilde{\vfi}(\pifi(a)) = \vfi(a) \ ,$$
hence $\th \in \cF_\vfi$. Because $\om$ is a normal positive functional  such that $\om \pifi = \th$, we have by definition that $\tilde{\th} = \om$.
So we see that $\cF \subseteq \{\,\tilde{\th} \mid \th \in \cF_\vfi \,\}$. Hence definition \ref{weight.def6} implies for every $x \in (\pifi(A)'')^+$ that
$$\hat{\psi}(x) = \sup \, \{ \, \om(x) \mid \om \in \cF \, \}
\leq \sup \, \{ \, \tilde{\th}(x) \mid \th \in \cF_\vfi \, \} = \tilde{\vfi}(x) \ .$$
\end{demo}

\medskip

In the next proposition, we use a technique of J. Verding and A. Van Daele (adapted in lemma \ref{app.lem2} to the strong topology framework) and combine it with the previous proposition. It is the crucial result to lift properties from the \cst-level to the \wst-level.

\begin{proposition} \label{weight.prop19}
Consider $x \in \cN_{\tilde{\vfi}}$. Then there exists a net $(a_i)_{i \in I}$ in $\Nfi$ such that
\begin{enumerate}
\item $\|a_i\| \leq \|x\|$ for every $i \in I$.
\item $(\pi_\vfi(a_i))_{i \in I}$ converges strongly$^*$ to $x$.
\item $(\,\lafi(a_i)\,\bigr)_{i \in I}$ is convergent in $H_\vfi$.
\end{enumerate}
\end{proposition}
\begin{demo}
Take $y \in \Nfi$ such  that $\pifi(y) = 0$. Then we have for all $\om \in \cG_\vfi$ that $T_\om^{\frac{1}{2}} \lafi(y) = \pifi(y) \xi_\om = 0$. Because $(T_\om)_{\om \in \cG_\vfi}$ converges strongly to 1, this implies that $\lafi(y) = 0$. So we can define a  linear map $\Gamma : \pifi(\Nfi) \rightarrow H_\vfi$ such that $\Gamma(\pifi(y)) = \lafi(y)$ for $y \in \Nfi$.

We can assume that $x \not= 0$. Because $x \in \cN_{\tilde{\vfi}}$, we have that $\hat{\psi}(x^* x) = \tilde{\vfi}(x^* x) < \infty$, so lemma \ref{weight.lem21} implies the existence of a net $(y_p)_{p \in P} \in \pifi(A)^+$ such that $\|y_p^2\| \leq \|x^* x\|$ for $p \in P$, $(y_p^2)_{p \in P}$ converges strongly to $x^* x$ and such that there exists $L \in \R^+$ such that $\psi(y_p^2) \leq L$ for $p \in P$. Then $(y_p)_{p \in P}$ converges strongly to $(x^* x)^{\frac{1}{2}}$ and $\|y_p\| \leq \|x\|$ for $p \in P$.

There exists clearly a partial isometry $u \in \pifi(A)''$ such that $x = u \,(x^* x)^{\frac{1}{2}}$. Using Kaplansky's density theorem, we find a net $(u_q)_{q \in Q}$ in $\pifi(A)$ which converges strongly$^*$ to $u$ and $\|u_q\| < 1$ for $q \in Q$.

Define $K = P \times Q$ and put the product order on $K$. In this way, $K$ becomes a directed set.

Choose $k=(p,q) \in K$. Then $u_q\,y_p \in \pifi(A)$ and $\|u_q \, y_p\| \leq \|u_q\|\,\|y_p\| < \|x\|$. Because  a $^*$-homomorphism sends the open unit ball onto the open unit ball of its image, we can find an element
$b_k \in A$ with $\|b_k\| < \|x\|$ such that $\pifi(b_k) = u_q \, y_p$.
Then we have also that
$$\vfi(b_k^* b_k) = \tilde{\vfi}(\pifi(b_k)^* \pifi(b_k))
= \tilde{\vfi}(y_p u_q^* u_q y_p)
\leq \tilde{\vfi}(y_p^2) = \psi(y_p^2) \leq L \ ,$$
which implies that $b_k \in \Nfi$ and $\|\lafi(b_k)\| \leq \sqrt{L}$.

It is also clear that the net $(\,\pifi(b_k)\,)_{k \in K}$ converges strongly$^*$ to $u\, (x^* x)^{\frac{1}{2}} = x$.

\medskip

Using lemma \ref{app.lem2}, we get a net $(z_i)_{i \in I}$ in the convex hull of $\{ \, \pifi(b_k) \mid k \in K \, \}$ such that $(z_i)_{i \in I}$ converges strongly$^*$ to $x$ and $(\Gamma(z_i))_{i \in I}$ is convergent in $H_\vfi$. It is clear that there exists for every $i \in I$ an element
$a_i$ in the convex hull of $\{\,b_k \mid k \in K\,\}$ such that
$\pifi(a_i) = z_i$. Then $(\pifi(a_i))_{i \in I}$ converges strongly$^*$ to $x$ and $(\lafi(a_i))_{i \in I}$ is convergent in $H_\vfi$.
Since $\|b_k\| \leq \|x\|$ for all $k \in K$, we have also that
$\|a_i\| \leq \|x\|$ for $i \in I$.
\end{demo}

\medskip

With this proposition in hand, it is now very easy to get the following natural and useful GNS-construction for $\tilde{\vfi}$.

\begin{proposition}
There exists a unique linear map $\la : \cN_{\tilde{\vfi}} \rightarrow H$ such that
\begin{itemize}
\item $(H_\vfi,\io,\la)$ is a GNS-construction for $\tilde{\vfi}$.
\item $\la(\pifi(a)) = \lafi(a)$ for all $a \in \Nfi$.
\end{itemize}
We call $(H_\vfi,\io,\la)$ the \wst-lift of $(H_\vfi,\pifi,\lafi)$
\end{proposition}
\begin{demo}
Uniqueness follows easily from the $\si$-strong$^*$ closedness of such a $\la$ and the previous proposition. We turn to the existence.
Therefore take any GNS-construction $(H_{\tilde{\vfi}},\pi_{\tilde{\vfi}},\la_{\tilde{\vfi}})$ for $\tilde{\vfi}$. Because $\tilde{\vfi} \, \pifi = \vfi$, we can define an isometry $U : H_\vfi \rightarrow H_{\tilde{\vfi}}$ such that $U \lafi(a) = \la_{\tilde{\vfi}}(\pifi(a))$ for $a \in \Nfi$.

Now take any $x \in \cN_{\tilde{\vfi}}$. By the previous proposition, there exists a bounded net $(a_i)_{i \in I}$ in $\Nfi$ such that
$(\pifi(a_i))_{i \in I}$ converges strongly$^*$ to $x$ and
$(\lafi(a_i))_{i \in I}$ converges to a vector $v \in H$.

So we get that the net $(\,\la_{\tilde{\vfi}}(\pifi(a_i))\,)_{i \in I}$ converges to $U v$. Hence the $\si$-strong$^*$ closedness of $\la_{\tilde{\vfi}}$ implies that $\la_{\tilde{\vfi}}(x) = U v$.

Consequently $U$ is a unitary. Now define the map $\la : \cN_{\tilde{\vfi}} \rightarrow H : x \mapsto U^* \la_{\tilde{\vfi}}(x)$. It is then easy to check that $\la$ satisfies all the requirements of the proposition.
\end{demo}

\medskip

\begin{remark} \rm \label{weight.rem1}
Remember that proposition \ref{weight.prop19} implies the following property. Consider $x \in \cN_{\tilde{\vfi}}$. Then there exists a net $(a_i)_{i \in I}$ in $\Nfi$ such that
\begin{enumerate}
\item $\|a_i\| \leq \|x\|$ for every $i \in I$.
\item $\bigl(\pi_\vfi(a_i)\bigr)_{i \in I}$ converges strongly$^*$ to $x$.
\item $\bigl(\lafi(a_i)\bigr)_{i \in I}$ converges to $\la(x)$.
\end{enumerate}
\end{remark}

\bigskip

Let us quickly extend all results to the multiplier algebra (which is very easy to do).

\begin{proposition}
We have the following properties.
\begin{enumerate}
\item $\tilde{\vfi}(\pifi(a)) = \vfi(a)$ for $a \in M(A)^+$.
\item We have for all $a \in \bar{\cN}_\vfi$ that $\pifi(a) \in \cN_{\tilde{\vfi}}$ and $\la(\pifi(a)) = \lafi(a)$.
\end{enumerate}
\end{proposition}
\begin{demo}
\begin{enumerate}
\item Take $a \in M(A)^+$. Take an approximate unit $(e_i)_{i \in I}$ for $A$. Then $(a^{\frac{1}{2}} e_i \, a^{\frac{1}{2}})_{i \in I}$ is an increasing bounded net in $A^+$ which converges strictly to $a$. So $(\pifi(a^{\frac{1}{2}} e_i \, a^{\frac{1}{2}}))_{i \in I}$ is an increasing bounded net in $\pifi(A)^+$ which converges strongly to $\pifi(a)$.

Consequently,
$$\tilde{\vfi}(\pifi(a)) = \sup \, \{ \, \tilde{\vfi}(\pifi(a^{\frac{1}{2}} e_i \, a^{\frac{1}{2}})) \mid i \in I \, \}
= \sup \, \{ \, \vfi(a^{\frac{1}{2}} e_i \, a^{\frac{1}{2}}) \mid i \in I \, \}
= \vfi(a) \ .$$
\item Take $a \in \bar{\cN}_\vfi$. Use the approximate unit from the first part of the proof. We have for every $i \in I$ that $e_i \, a \in \Nfi$, which implies that $\pifi(e_i \ a) \in \cN_{\tilde{\vfi}}$ and
$$\la(\pifi(e_i \ a)) = \lafi(e_i a) = \pifi(e_i) \, \lafi(a) \ .$$
So we have a bounded net $(\, \pifi(e_i a)\,)_{i \in I}$ in $\cN_{\tilde{\vfi}}$ which converges strongly$^*$ to $\pifi(a)$ and has the property that
the net $(\,\la(\pifi(e_i \ a))\,)_{i  \in I}$ converges to $\lafi(a)$.

Now the $\si$-strong$^*$ closedness of $\la$ implies that
$\pifi(a) \in \cN_{\tilde{\vfi}}$ and $\la(\pifi(a))= \lafi(a)$.
\end{enumerate}
\end{demo}

\bigskip\medskip

A first easy application can be found in the context of $^*$-automorphisms on $A$ which leave $\vfi$ relatively invariant.

\begin{proposition}
Consider a $^*$-automorphism $\al$ on $A$ such that there exists a strictly positive number $r > 0$ such that $\vfi \, \al = r \, \vfi$. Define the unitary operator $U$ on $H_\vfi$ such that $U \lafi(a) = r^{-\frac{1}{2}} \, \lafi(\al(a))$ for $a \in \Nfi$.
Then there exists a unique $^*$-automorphism $\tilde{\al}$ on $\pifi(A)''$ such that $\tilde{\al}\,\pifi = \pifi \, \al$. The $^*$-automorphism $\tilde{\al}$  is implemented by $U$, i.e. $\tilde{\al}(x) = U \,x \,U^*$ for all $x \in \pifi(A)''$.
We have moreover that $\tilde{\vfi} \, \tilde{\al} = r \, \tilde{\vfi}$ and that $U \la(x) = r^{-\frac{1}{2}} \, \la(\tilde{\al}(x))$ for all $x \in \cN_{\tilde{\vfi}}$.
\end{proposition}
\begin{demo}
Let $b$ be an element in $A$. Then we have for all $a \in \Nfi$ that
\begin{eqnarray*}
U \pifi(b) \lafi(a) & = & U \lafi(b a) = r^{-\frac{1}{2}} \, \lafi(\al(b a)) = r^{-\frac{1}{2}} \lafi(\al(b)\al(a)) \\
& = & r^{-\frac{1}{2}} \pifi(\al(b)) \lafi(\al(a)) = \pifi(\al(b)) U \lafi(a) \ ,
\end{eqnarray*}
implying that $U \pifi(b) = \pifi(\al(b)) U$, hence $U \pifi(b) U^* = \pifi(\al(b))$.
This implies first of all that $U \pifi(A)'' U^*$ $= \pifi(A)''$ so we can define a $^*$-automorphism $\tilde{\al}$ on $\pifi(A)''$ such that $\tilde{\al}(x) = U\, x \, U^*$ for $x \in \pifi(A)''$. Then we have immediately that $\tilde{\al}\, \pifi = \pifi \, \al$.

Choose $x \in \cN_{\tilde{\vfi}}$. By the remark before this proposition, there exists a bounded net $(a_i)_{i \in I}$ in $\Nfi$ such that
$(\pifi(a_i))_{i \in I}$ converges strongly$^*$ to $x$ and
$(\lafi(a_i))_{i \in I}$ converges to $\la(x)$.

Then $(\,\tilde{\al}(\pifi(a_i))\,)_{i \in I}$ is a bounded net which converges strongly$^*$ to $\tilde{\al}(x)$.

We have for $i \in I$ that $\tilde{\al}(\pifi(a_i)) = \pifi(\al(a_i))$ which implies that $\tilde{\al}(\pifi(a_i))$ belongs to $\cN_{\tilde{\vfi}}$ and
that
$$\la\bigl(\tilde{\al}(\pifi(a_i))\bigr) = \la\bigr(\pifi(\al(a_i))\bigr)
= \lafi(\al(a_i)) = r^{\frac{1}{2}} \, U \lafi(a_i) \ .$$
This implies that $\bigl(\,\la\bigl(\tilde{\al}(\pifi(a_i))\bigr)\,\bigr)_{i \in I}$ converges to $r^{\frac{1}{2}} \, U \la(x) \ .$
So the $\si$-strong$^*$ closedness of $\la$ implies that
$\tilde{\al}(x)$ belongs to $\cN_{\tilde{\vfi}}$ and
that $\la(\tilde{\al}(x)) = r^{\frac{1}{2}} \, U \la(x)$.
From this, we also conclude that $\tilde{\vfi}\,\tilde{\al} = r  \, \tilde{\vfi}$.
\end{demo}

\medskip

If we apply this result to a norm continuous one-parameter group under which $\vfi$ is relatively invariant, we arrive at the following conclusions.

\begin{corollary}
Consider a norm continuous one-parameter group $\al$ on $A$ such that there exists a strictly positive number $r > 0$ such that $\vfi \, \al_t = r^t \, \vfi$ for $t \in \R$. Define the strictly positive  operator $T$ in $H_\vfi$ such that $T^{it} \lafi(a) = r^{-\frac{t}{2}} \, \lafi(\al_t(a))$ for $a \in \Nfi$ and $t \in \R$.

Then there exists a unique strongly continuous one-parameter group  $\tilde{\al}$ on $\pifi(A)''$ such that $\tilde{\al}_t\,\pifi = \pifi \, \al_t$ for $t \in \R$. The one-parameter group $\tilde{\al}$ is implemented by $T$, i.e. $\tilde{\al}_t(x) = T^{it} \, x \, T^{-it}$ for all $t \in \R$ and $x \in \pifi(A)''$.
We have moreover that $\tilde{\vfi} \, \tilde{\al}_t = r^t \, \tilde{\vfi}$ for $t \in \R$ and that $T^{it} \la(x) = r^{-\frac{t}{2}} \, \la(\tilde{\al}_t(x))$ for all $t \in \R$ and $x \in \cN_{\tilde{\vfi}}$.
\end{corollary}

\bigskip\medskip

For the rest of this section, we assume that $\vfi$ is a KMS weight with modular group $\si$. We define moreover
\begin{itemize}
\item $J =$ the modular conjugation of $\vfi$ in the GNS-construction $(H_\vfi,\pifi,\lafi)$.
\item $\nab =$ the modular operator of $\vfi$ in the GNS-construction $(H_\vfi,\pifi,\lafi)$.
\end{itemize}

By the previous proposition, we know that there exists a unique strongly continuous one-parameter group $\tilde{\si}$ on $\pifi(A)''$ such that
$\tilde{\si}_t \, \pifi = \pifi \, \si_t$ for $t \in R$.
Recall also that $\tilde{\si}_t(x) = \nab^{it} x \nab^{-it}$ for all $t \in \R$ and $x \in \pifi(A)''$.

\smallskip

We know also that $\tilde{\vfi} \, \tilde{\si}_t = \tilde{\vfi}$ for $t \in \R$. Moreover, $\la(\tilde{\si}_t(x)) = \nab^{it} \la(x)$ for $t \in \R$ and
$x \in \cN_{\tilde{\vfi}}$.

\bigskip

We will quickly look into the left Hilbert algebra theory. This will allow us to use Tomita-Takesaki theory.

In section 6 of \cite{JK1}, we introduced the usual left Hilbert algebra $\cU$ on $H$. Just define the subspace $\cU = \lafi(\Nfi \cap \Nfi^*)$ of $H$. We turn
$\cU$ into a $^*$-algebra such that
\begin{itemize}
\item $\lafi(a) \, \lafi(b) = \lafi(a\, b)$ for $a,b \in \Nfi \cap \Nfi^*$.
\item $\lafi(a)^* = \lafi(a^*)$ for $a \in \Nfi \cap \Nfi^*$.
\end{itemize}
In section 6 of \cite{JK1}, we showed that $\cU$ is a left Hilbert algebra on $H$. It is also true  for every $a \in \Nfi$ that $\lafi(a)$ is left bounded w.r.t. $\cU$ and that $L_{\lafi(a)} = \pifi(a)$. Hence $\cL(U) = \pifi(A)''$.

We also proved that $\lafi(\Nfi \cap \Nfi^*)$ is a core for $J \nab^{\frac{1}{2}}$ and that
$J \nab^{\frac{1}{2}} \lafi(x) = \lafi(x^*)$ for all $x \in \Nfi \cap \Nfi^*$.

Hence $J$ is the modular conjugation of $\cU$ and $\nab$ is the
modular operator of $\cU$. \inlabel{ster4}

\medskip

Now Tomita-Takesaki theory implies immediately the following result.

\begin{proposition}
We have that $J \pifi(A)'' J = \pifi(A)'$.
\end{proposition}

\medskip

As to be expected, the weight $\tilde{\vfi}$ equals the weight which can be constructed from $\cU$.

\begin{proposition}
The weight $\tilde{\vfi}$ is the n.f.s weight on $\pifi(A)''$ canonically associated to $\cU$. We have moreover for all $a \in \cN_{\tilde{\vfi}}$ that $\la(a)$ is left bounded w.r.t. $\cU$ and $L_\la(a) = a$.
\end{proposition}
\begin{demo}
Denote the n.f.s. weight on $\pifi(A)''$ canonically associated to $\cU$ by
$\psi$. Let $(H,\io,\Gamma)$ be the canonical GNS-construction of $\psi$ (defined with the use of left bounded elements).

By the remarks above, we know that $\pifi(\cN_\vfi) \subseteq \cN_\psi$ and
that $\Gamma(\pifi(a)) = \lafi(a)$ for $a \in \cN_\vfi$.

Using remark \ref{weight.rem1} and the $\si$-strong$^*$ closedness of $\Gamma$, we get easily that $\cN_{\tilde{\vfi}} \subseteq \cN_\psi$ and $\Gamma(x) = \la(x)$ for all $x \in \cN_{\tilde{\vfi}}$. So $\psi$ is an extension $\tilde{\vfi}$. By statement~\ref{ster4} above, we know that $\tilde{\si}$ is the modular group of $\psi$. Since $\tilde{\vfi}$ is invariant under $\tilde{\si}$,  the \wst-version of corollary \ref{sweight.cor1} implies that $\tilde{\vfi} = \psi$ and $\la = \Gamma$.
\end{demo}

\medskip

Therefore we conclude that $\tilde{\vfi}$ is a n.f.s weight on $\pifi(A)''$
with modular group $\tilde{\si}$.

\bigskip\medskip

\section{Slicing with weights}

Fix two \cst-algebras $A$ and $B$ together with a proper weight $\vfi$
on $B$ and a GNS-construction $(H_\vfi,\pifi,\lafi)$ for $\vfi$. An
important tool in the theory of \cst-algebraic quantum groups is the
slice map $\io \ot \vfi$. If $A$ and $B$ arise from locally compact
spaces, $\vfi$ is implemented by a regular Borel measure $\mu$ and
$\io \ot \vfi$ will integrate out the second variable with respect to
$\mu$.

\smallskip

\medskip

In this section, we will define $\io \ot \vfi$ and prove all the necessary
properties of this slice map.

\begin{definition} \label{weight.def2}
We will use the following notations:
\begin{itemize}
\item We define the set  $${\bar{\cM}}_{\io \ot \vfi}^+ = \{\, x \in M(A \ot B)^+
\mid \text{\ the net \  } \bigl( \,(\io \ot \om)(x)\,\bigr)_{\om \in {\cal G}_\vfi}
\text{\  is strictly convergent in  } M(A) \, \} \ .$$
\item For $x \in {\bar{\cM}}_{\io \ot \vfi}^+$, we define $(\io \ot \vfi)(x)$ to be the element in $M(A)$ such that the net
$\bigl( \,(\io \ot \om)(x)\,\bigr)_{\om \in {\cal G}_\vfi}$
converges strictly to $(\io \ot \vfi)(x)$.
\end{itemize}
\end{definition}

\bigskip

It is possible to get strict convergence under weaker convergence conditions. We will look into two of these conditions. First we need a simple lemma.

\begin{lemma} \label{weight.lem1}
Let $b \in M(A)^+$. Then we have that $\|b\, a\|^2 \leq \|b\| \, \| a^* b\, a \|$ for $a \in A$.
\end{lemma}

The proof of this lemma is very simple because
$(b\,a)^* (b\, a) = a^* b^2 \, a \leq \|b\| \, a^* b \, a  , $
where we used the fact that $b^2 = b^{\frac{1}{2}} \, b \, b^{\frac{1}{2}} \leq \|b\| \, b$.
We will apply this small result in two situations.

\begin{lemma}  \label{weight.lem4}
Let $(b_i)_{i \in I}$ be a net in $M(A)^+$ and $b$ an element in $M(A)^+$ such that $b_i \leq b$ for every $i \in I$.
Then $(b_i)_{i \in I}$ converges strictly to $b$ $\Leftrightarrow$
$(\,a^* b_i  \,a)_{i \in I}$ converges to $a^* b \, a$ for every $a \in A$.
\end{lemma}
\begin{demo}
We have for every $i \in I$ that $\|b_i\| \leq \|b\|$.
Using the previous lemma, we get for every $i \in I$ and $a \in A$ that
$$\|b \, a - b_i \,a\|^2 \leq \|b - b_i\| \, \|a^* (b-b_i) \, a\| \leq 2 \|b\| \, \|a^* b \, a - a b_i \, a\|$$
and the lemma follows.
\end{demo}

\medskip

Using the Uniform Boundedness Principle, we can also prove easily the following result.

\begin{result}
Let $(b_i)_{i \in I}$ be an increasing net in $M(A)^+$.
Then $(b_i)_{i \in I}$ is strictly convergent in $M(A)^+$ $\Leftrightarrow$ the net $(\, a^* b_i \, a\,)_{i \in I}$ is convergent for every $a \in A$.
\end{result}
\begin{demo}
Suppose that $(\,a^* b_i \, a\, )_{i \in I}$ is convergent for every $a \in A$. First, we prove that $(b_i)_{i \in I}$ is bounded.

\begin{list}{}{\setlength{\leftmargin}{.4 cm}}

\item Choose $c \in A$. Because $(\,c^* b_i \, c\,)$ is convergent, there exists a positive number $N_c$ and an element $i_0$ such that $\| c^* b_i \, c\| \leq N_c$ for every $i \in I$ with $i \geq i_0$.

For every $j \in I$ there exists an element $i \in I$ with $i \geq i_0$
and $i \geq j$, implying that
$$ \| c^* b_j \,c \| \leq \| c^* b_i \, b \| \leq N_c .$$
So we get that the net $(\,c^* b_i \, c \,)_{i \in I}$ is bounded.

\smallskip

Let us fix $d \in A$. By polarization and the above result, we have for every $e \in A$ that the net $(\,e^* b_i \, d\,)_{i \in I}$ is bounded.
Using the uniform boundedness principle, we get that the net
$(\,b_i \,d\, )_{i \in I}$ is bounded.

\smallskip

Applying the uniform boundedness principle once again, we see that the net $(b_i)_{i \in I}$ is bounded.

\end{list}

Hence there exists a strictly positive
number $M$ such that $\|b_i\| \leq M$ for every $i \in I$.

Choose $a \in A$. Take $\vep > 0$. Then there exist an element $i_1 \in I$ such that $\| \, a^* b_i \, a - a^*  b_{i_1}\, a \|  \leq \frac{\vep}{2M}$ for every $i \in I$ with $i \geq i_1$. Using lemma \ref{weight.lem1}, we have for every $i \in I$ with $i \geq i_1$ that
$$\| b_i \,a - b_{i_1} \, a \|^2  \leq  \|b_i - b_{i_1}\| \, \| a^* (b_i - b_{i_1}) \, a\|
 \leq  2 M \, \frac{\vep}{2M} = \vep .$$

So we see that $(\,b_i \, a\,)_{i \in I}$ is Cauchy and hence convergent in $A$.

\medskip

From this all, we infer the existence of a mapping $T$ from $A$ into $A$ such that $(\,b_i\,a\,)_{i \in I}$ converges to $T(a)$ for every $a \in A$. It follows immediately that $c^* T(a) = T(c)^* a$ for every $a,c \in A$. This implies the existence of an element $b \in M(A)$ with $b^* = b$ such that $b \, a = T(a)$ for $a \in A$. It  is also clear that $a^* b \, a \geq 0$ for every $a \in A$, which implies that $b \geq 0$.
\end{demo}

\medskip

Applying this result to definition \ref{weight.def2},  we get immediately the following one.

\begin{result} \label{weight.res3}
Consider $x \in M(A \ot B)^+$. Then $x$ belongs to $\bar{\cM}_{\io \ot \vfi}^+$ $\Leftrightarrow$ we have for all $a \in A$ that the net
$\bigl( \,a^* (\io \ot \om)(x) \,a\,\bigr)_{\om \in {\cal G}_\vfi}$
is convergent in $A$.
\end{result}

\medskip

Combining this criterium with the Cauchy property, it becomes easy to prove the third property of $\io \ot \vfi$ in the next result. Just notice that
$$a^* (\io \ot \th)(x) \,a - a^* (\io \ot \om)(x) \,a
\leq a^* (\io \ot \th)(y) \,a - a^* (\io \ot \om)(y) \,a$$
for all $a \in A$ and $\om,\th \in \cG_\vfi$ such that $\om \leq \th$.
The proofs of the other properties are straightforward.

\begin{result}
The slice map $\io \ot \vfi$ satisfies the following algebraic properties:
\begin{enumerate}
\item We have for $x,y \in \bar{\cM}_{\io \ot \vfi}^+$ that $x+y \in \bar{\cM}_{\io \ot \vfi}^+$ and $(\io \ot \vfi)(x + y) = (\io \ot \vfi)(x) + (\io \ot \vfi)(y)$.
\item We have for $x \in \bar{\cM}_{\io \ot \vfi}^+$ and $\lambda \in \R^+$ that $\lambda \, x \in \bar{\cM}_{\io \ot \vfi}^+$ and
$(\io \ot \vfi)(\lambda \,x) = \lambda \, (\io \ot \vfi)(x)$.
\item Consider $y \in \bar{\cM}_{\io \ot \vfi}^+$ and $x \in M(A \ot B)^+$ such that $ x \leq y$. Then $x \in \bar{\cM}_{\io \ot \vfi}^+$ and $(\io \ot \vfi)(x) \leq (\io \ot \vfi)(y)$.
\item We have for $a \in M(A)^+$ and $b \in \bar{\cM}_\vfi^+$ that
$a \ot b \in \bar{\cM}_{\io \ot \vfi}^+$ and $(\io \ot \vfi)(a \ot b) = a \, \vfi(b)$.
\end{enumerate}
\end{result}

\medskip

As for ordinary weights, this allows us to extend $\io \ot \vfi$ to a linear mapping defined on a subalgebra of $M(A \ot B)$ :

\begin{notation}
The following notations will be used
\begin{itemize}
\item We define $\bar{\cM}_{\io \ot \vfi}$ as the linear span of $\bar{\cM}_{\io \ot \vfi}^+$ in $M(A \ot B)$.
Then $\bar{\cM}_{\io \ot \vfi}$ is a sub-$^*$-algebra of $M(A \ot B)$ such that $\bar{\cM}_{\io \ot \vfi}^+ = \bar{\cM}_{\io \ot \vfi}
\cap M(A \ot B)^+$.
\item There exists a unique linear map $F : \bar{\cM}_{\io \ot \vfi} \rightarrow M(A)$ such that $F(x) = (\io \ot \vfi)(x)$ for $x \in \bar{\cM}_{\io \ot \vfi}^+$.

For every $x \in \bar{\cM}_{\io \ot \vfi}$, we put $(\io \ot \vfi)(x) = F(x)$.
\item We define $\bar{\cN}_{\io \ot \vfi} = \{ x \in M(A \ot B) \mid x^* x \in \bar{\cM}_{\io \ot \vfi}^+ \}$. Then
$\bar{\cN}_{\io \ot \vfi}$ is a left ideal in $M(A \ot B)$
such that $\bar{\cM}_{\io \ot \vfi} =
\bar{\cN}_{\io \ot \vfi}^{\ *} \ \bar{\cN}_{\io \ot \vfi}$.
\end{itemize}
\end{notation}

\medskip

Then we have immediately the following natural properties concerning $\bar{\cM}_{\io \ot \vfi}$.

\begin{lemma}
The following properties hold:
\begin{itemize}
\item Let $x \in \bar{\cM}_{\io \ot \vfi}$. Then the net
$\bigl(\,(\io \ot \om)(x)\,\bigl)_{\om \in \cG_\vfi}$ converges strictly to
$(\io \ot \vfi)(x)$.
\item We have for $a \in M(A)$ and  $b \in \bar{\cM}_\vfi$ that
$a \ot b \in \bar{\cM}_{\io \ot \vfi}$ and $(\io \ot \vfi)(a \ot b) = a \, \vfi(b)$.
\item We have for $a \in M(A)$ and $b \in \bar{\cN}_\vfi$ that
$a \ot b \in \bar{\cN}_{\io \ot \vfi}$.
\end{itemize}
\end{lemma}

\bigskip

The next Fubini-like proposition is an easy consequence of the definitions above.

\begin{proposition} \label{weight.prop1}
Consider  $x \in \bar{\cM}_{\io \ot \vfi}$ and $\th \in A^*$. Then
$(\th \ot \io)(x)$ belongs to $\bar{\cM}_{\vfi}$ and
$$\vfi\bigl((\th \ot \io)(x)\bigr) = \th\bigl((\io \ot \vfi)(x)\bigr) \ .$$
\end{proposition}
\begin{demo}
Because any element in $A^*$ can be written as a linear combination of positive linear functionals and any element of $\bar{\cM}_{\io \ot \vfi}$ can be written as a linear combination of elements in $\bar{\cM}_{\io \ot \vfi}^+$, it is enough to prove the result if  $x \in \bar{\cM}_{\io \ot \vfi}^+$ and $\th \in A^*_+$.

But we have for every $\om \in \cG_\vfi$ that $\om\bigl((\th \ot \io)(x))
= \th((\io \ot \om)(x))$. Therefore definition \ref{weight.def2} and the strict continuity of $\th$ imply that the net
$\bigl(\,\om\bigl((\th \ot \io)(x)\bigr)\,\bigr)_{\om \in \cG_\vfi}$ converges to $\th((\io \ot \vfi)(x))$. Hence definition \ref{weight1.def3} implies that $(\th \ot \io)(x)$ belongs to
$\bar{\cM}_\vfi$ and $\vfi\bigl((\th \ot \io)(x)\bigr) = \th\bigl((\io \ot \vfi)(x))$.
\end{demo}

\medskip

Let $\th$ be a continuous linear functional on a \cst-algebra $A$. Proposition~III.4.6 of \cite{Tak} implies the existence of a unique positive linear functional
$|\th|$ on $A$ such that $\|\,|\th|\,\| = \|\th\|$ and
$|\th(a)|^2 \leq \|\th\| \, |\th|(a^* a)$ for $a \in A$.

\begin{lemma}
Let $x \in M(A \ot B)$. Then
$$(\th \ot \io)(x)^* (\th \ot \io)(x) \leq \|\th\| \, \, (|\th| \ot \io)(x^* x) \ .$$
\end{lemma}
\begin{demo}
Choose $y \in A \odot B$. So there exists $a_1,\ldots,a_n \in A$ and
$b_1,\ldots,b_n \in B$ such that $y = \sum_{i=1}^n a_i \ot b_i$. Define $M \in M_n(\C)$ by $M_{ij} = \|\th\| \, |\th|(a_i^* a_j) - \overline{\th(a_i)} \th(a_j)$ for $i,j = 1,\ldots,n$.

Then we have for $c_1,\ldots,c_n \in \C$ that
$$\sum_{i,j=1}^n \overline{c_i} M_{ij} c_j
= \|\th\| \, |\th|\bigl(\,[\,\sum_{i=1}^n c_i a_i ]^* \, [\,\sum_{i=1}^n c_i a_i ]\,\bigr) - \overline{\th(\sum_{i=1}^n c_i a_i)} \, \th(\sum_{i=1}^n c_i a_i) \ ,$$
which is positive. This implies that the matrix $M$ is positive.

It follows that $\sum_{i,j=1}^n b_i^* M_{ij} b_j$ is positive. Now
$$\|\th\| \, (|\th| \ot \io)(y^* y) - (\th \ot \io)(y)^* (\th \ot \io)(y)
= \sum_{i,j}^n \|\th\| \ b_i^* |\th|(a_i^* a_j) b_j - b_i^* \overline{\th(a_i)} \th(a_j) b_j = \sum b_i^* M_{ij} b_j \geq 0 \ .$$
So we see that $(\th \ot \io)(y)^* (\th \ot \io)(y)
\leq \|\th\| \, (|\th| \ot \io)(y^* y)$. The lemma follows now by continuity considerations.
\end{demo}

\medskip

Combining this lemma with proposition \ref{weight.prop1}, we get immediately the following one.

\begin{proposition}
Consider  $x \in \bar{\cN}_{\io \ot \vfi}$ and $\th \in A^*$. Then
$(\th \ot \io)(x)$ belongs to $\bar{\cN}_{\vfi}$ and
$$\vfi\bigl((\th \ot \io)(x)^* (\th \ot \io)(x)\bigr) \leq \|\th\| \ |\th|\bigl((\io \ot \vfi)(x^* x)\bigr) \ .$$
\end{proposition}

\medskip

Now we want to prove a converse of proposition \ref{weight.prop1}. This will be very easy once the following lemma is proven (a result we borrowed from \cite{Rie}).

\begin{lemma}
Consider an increasing net $(a_i)_{i \in I}$ in $A^+$ and an element $a \in A^+$ such that
$\om(a) = \sup \{\,\om(a_i) \mid i \in I \, \}$ for all $\om \in A^*_+$. Then the net $(a_i)_{i \in I}$ converges to $a$.
\end{lemma}
\begin{demo}
Define $S = \{ \, \om \in A^*_+ \mid  \|\om\| \leq 1 \, \}$ and equip $S$ with the weak$^*$-topology. Then $S$ becomes a compact Hausdorff space.

For every $i \in I$, we define the function $f_i \in C(S)^+$ such that
$f_i(\om)= \om(a_i)$ for $\om \in A^*_+$. So $(f_i)_{i \in I}$ is an increasing net in $C(S)^+$.
We also define the function $f \in C(S)^+$ such that $f(\om)= \om(a)$ for $\om \in A^*_+$.

By assumption, the net $(f_i)_{i \in I}$ converges pointwise to $f$. Therefore Dini's theorem implies that $(f_i)_{i \in I}$ converges uniformly to $f$.

Because $\|x\| = \sup \{ \, |\om(x)| \mid \om \in S  \, \}$ for $x \in A^+$, this implies that $(a_i)_{i \in I}$ converges in norm to $a$.
\end{demo}

\medskip

Combining this lemma with lemma \ref{weight.lem4}, it is not difficult to prove the following strict version.

\begin{lemma} \label{weight.lem5}
Consider  an increasing net $(a_i)_{i \in I}$ in $M(A)^+$ and an element $a \in M(A)^+$ such that
$\om(a) = \sup \{\,\om(a_i) \mid i \in I \, \}$ for all $\om \in A^*_+$. Then the net $(a_i)_{i \in I}$ converges strictly to $a$.
\end{lemma}

\medskip

Now we can easily prove the converse of proposition \ref{weight.prop1}.

\begin{proposition} \label{weight.prop5a}
Consider $x \in M(A \ot B)^+$ and $a \in M(A)^+$ such that $(\th \ot \io)(x)$ belongs to $\bar{\cM}_\vfi^+$ and $\vfi\bigl((\th \ot \io)(x)\bigr) = \th(a)$ for all
$\th \in A^*_+$. Then $x$ belongs to $\bar{\cM}_{\io \ot \vfi}^+$ and $(\io \ot \vfi)(x) = a$.
\end{proposition}
\begin{demo}
Take $\th \in A^*_+$. Then we have for every $\om \in \cG_\vfi$ that
$\th((\io \ot \om)(x)) = \om((\th \ot \io)(x))$ so that definition \ref{weight1.def3} implies that the net $\bigl(\,\th((\io \ot \om)(x))\,\bigr)_{\om \in \cG_\vfi}$ converges to $\vfi((\th \ot \io)(x)) = \th(a)$.

Hence the previous lemma implies that the net $\bigl(\,(\io \ot \om)(x)\,\bigr)_{\om \in \cG_\vfi}$ converges strictly to $a$. Therefore the proposition follows from definition \ref{weight.def2}.
\end{demo}

\bigskip

In this setting of \lq\cst-valued weights\rq , the Cauchy-Schwarz inequality is generalized in the following way.

\begin{proposition} \label{weight.prop3}
Let $x,y \in \bar{\cN}_{\io \ot \vfi}$. Then
$$[(\io \ot \vfi)(y^* x)]^* \,[(\io \ot \vfi)(y^* x)]  \leq \|(\io \ot \vfi)(y^* y)\| \,\, (\io \ot \vfi)(x^* x) \ .$$
\end{proposition}
\begin{demo}
Take $\om \in \cG_\vfi$. Let $(H,\pi,v)$ be a cyclic GNS-construction for $\om$ and define $\th_v \in B(\C,H)$ such that $\th_v(c) = c v$ for $c \in \C$.

Look now at the Hilbert \cst-module $\cL(A,A \ot H)$ over $M(A)$ : $\langle S , T \rangle = T^* S$ for all operators $S,T \in \cL(A,A \ot H)$.
Then $(\io \ot \om)(b^* a) = \bigl(\,(\io \ot \pi)(b)(1 \ot \th_v)\,\bigr)^*
\bigl(\,(\io \ot \pi)(a)(1 \ot \th_v)\,\bigr)$ for $a,b \in M(A \ot B)$
and
\begin{eqnarray*}
& & [(\io \ot \om)(y^* x)]^* \,[(\io \ot \om)(y^* x)] \\
& & \spat = [\langle (\io \ot \pi)(x)(1 \ot \th_v) , (\io \ot \pi)(y)(1 \ot \th_v) \rangle]^* \,
[\langle (\io \ot \pi)(x)(1 \ot \th_v) , (\io \ot \pi)(y)(1 \ot \th_v) \rangle] \\
& & \spat  \leq \| \langle (\io \ot \pi)(y)(1 \ot \th_v) , (\io \ot \pi)(y)(1 \ot \th_v) \rangle \| \,  \langle (\io \ot \pi)(x)(1 \ot \th_v) , (\io \ot \pi)(x)(1 \ot \th_v) \rangle \\
& & \spat =  \|(\io \ot \om)(y^* y)\| \,\, (\io \ot \om)(x^* x) \ .
\end{eqnarray*}

Because $\bigl(\,(\io \ot \om)(b^* a) \,\bigr)_{\om \in \cG_\vfi}$ is a bounded net which converges strictly to $(\io \ot \vfi)(b^* a)$ for every $a,b \in \bar{\cN}_{\io \ot \vfi}$, the proposition follows immediately from the above inequality.
\end{demo}

\bigskip

In the next part, we are going to construct a \lq KSGNS-construction\rq\ for the  slice map $\io \ot \vfi$ which is a special case of the KSGNS-constructions considered in \cite{JK2}.

\medskip

\begin{result}  \label{weight.res2}
Consider $x \in \bar{\cN}_{\io \ot \vfi}$ and $v \in H_\vfi$. Then there exists a unique element $q \in M(A)$ such that \newline $\th(q) = \langle \lafi\bigl((\th \ot \io)(x)\bigr) , v \rangle$ for $\th \in A^*$.
\end{result}
\begin{demo}
For $a \in \Nfi$, we put $q(a) = (\io \ot \vfi)((1 \ot a^*) x) \in M(A)$.

Using proposition \ref{weight.prop3}, we see
for all $a,b \in \Nfi$ that
\begin{eqnarray*}
& & \|q(a) - q(b) \|^2  =  \| q(a-b)^* q(a-b) \|   \\
& & \spat = \|(\io \ot \vfi)((1 \ot (a-b)^*)\,x)^*\, (\io \ot \vfi)((1 \ot (a-b)^*)\,x)\| \\
& & \spat \leq \|(\io \ot \vfi)([1 \ot (a-b)]^*[1 \ot (a-b)]) \| \, \| (\io \ot \vfi)(x^* x) \| \\
& & \spat = \| \lafi(a) - \lafi(b) \|^2 \, \| (\io \ot \vfi)(x^* x) \| \ .
\end{eqnarray*}
Now take a sequence $(a_n)_{n=1}^\infty$ in $\Nfi$ such that
$(\lafi(a_n))_{n=1}^\infty$ converges to $v$. Then the previous inequality implies that $(q(a_n))_{n=1}^\infty$ is Cauchy and hence convergent in $M(A)$. So there exists an element $q \in M(A)$ such that $(q(a_n))_{n=1}^\infty$ converges in norm to $q$.

We have for every $\th \in A^*$ and $n \in \N$ that
$\th(q(a_n)) =  \vfi( a_n^* (\th \ot \io)(x)) = \langle \lafi((\th \ot \io)(x)) , \lafi(a_n) \rangle \ ,$
which implies that
$(\,\th(q(a_n))\,)_{n=1}^\infty$ converges to $\langle \lafi((\th \ot \io)(x)) , v \rangle$. So $\th(q)$ must be equal to
$\langle \lafi((\th \ot \io)(x)) , v \rangle$.
\end{demo}

\medskip

\begin{lemma}
Consider a non-degenerate $^*$-representation
$\pi$ of $A$ on a Hilbert space $H$ and an orthonormal basis $(e_i)_{i \in I}$ for $H$.
Let $x,y \in M(A \ot B)$ and $v,w \in H$. Then the net $$\bigl(\,\sum_{i \in J} \, (\om_{w,e_i} \ot \io)(y)^* (\om_{v,e_i} \ot \io)(x)\,\bigr)_{J \in F(I)}$$ is bounded and converges strictly to $(\om_{v,w}  \ot \io)(y^* x)$.
\end{lemma}
\begin{demo}
By polarization, it is enough to prove the result for $y=x$ and $v=w$.

Choose $\mu \in B^*_+$ and take a cyclic GNS-construction $(K,\th,u)$ for $\mu$. Fix also an orthonormal basis $(f_l)_{l \in L}$ for $K$.
Then
\begin{eqnarray*}
& & \sum_{i \in I} \mu\bigl(\, (\om_{v,e_i} \ot \io)(x)^* (\om_{v,e_i} \ot \io)(x)\,\bigr) = \sum_{i \in I} \langle \th((\om_{v,e_i} \ot \io)(x)) u , \th((\om_{v,e_i} \ot \io)(x)) u \rangle \\
& & \spat = \sum_{i \in I} \sum_{l \in L} \langle \th((\om_{v,e_i} \ot \io)(x)) u , f_l \rangle \,\, \langle f_l , \th((\om_{v,e_i} \ot \io)(x)) u \rangle \\
& & \spat = \sum_{i \in I} \sum_{l \in L} \langle (\pi \ot \th)(x) (v \ot u) , e_i \ot  f_l \rangle \, \, \langle e_i \ot f_l , (\pi \ot \th)(x) (v \ot u) \rangle \\
& & \spat = \langle (\pi \ot \th)(x) (v \ot u) , (\pi \ot \th)(x) (v \ot u) \rangle
= (\om_{v,v} \ot \om_{u,u})(x^* x) = \mu((\om_{v,v} \ot \io)(x^* x)) \ .
\end{eqnarray*}
Consequently, the lemma follows from lemma \ref{weight.lem5}.
\end{demo}

\medskip

We will now apply these results to get a sort of KSGNS-construction for the \lq\cst-valued weight\rq \ $\io \ot \vfi$:

\begin{proposition} \label{weight.prop2}
There exists a unique linear map $\la : \bar{\cN}_{\io \ot \vfi} \rightarrow \cL(A,A \ot H_\vfi)$ such that
$$\la(x)^* (a \ot \lafi(b))  = (\io \ot \vfi)(x^* (a \ot b))$$
for $a \in A$ and $b \in \Nfi$.

\smallskip

For $x \in \bar{\cN}_{\io \ot \vfi}$, we put $(\io \ot \lafi)(x) =
\la(x)$. Then we have the following properties:
\begin{itemize}
\item We have for all $x,y \in \bar{\cN}_{\io \ot \vfi}$ that
$(\io \ot \lafi)(y)^* (\io \ot \lafi)(x) = (\io \ot \vfi)(y^* x)$.
\item Consider $a \in M(A)$ and $b \in \bar{\cN}_\vfi$, then
$(\io \ot \lafi)(a \ot b) = a \ot \lafi(b)$.
\item We have for $x \in M(A \ot B)$ and $y \in \bar{\cN}_{\io \ot \vfi}$ that $(\io \ot \lafi)(x\, y) = (\io \ot \pifi)(x) (\io \ot \lafi)(y)$.
\end{itemize}
\end{proposition}
\begin{demo}
Fix an orthonormal basis $(e_i)_{i \in I}$ for $H_\vfi$.

Take $x \in \bar{\cN}_{\io \ot \vfi}$. Using result \ref{weight.res2}, we define for every $i \in I$ the element $q_i \in M(A)$ such that
$\eta(q_i) = \langle \lafi((\eta \ot \io)(x)) , e_i \rangle$ for $\eta \in A^*$.

\medskip

Choose $\mu \in A^*_+$ and a cyclic GNS-construction $(K,\th,v)$ for $\mu$. Fix also an orthonormal basis $(f_l)_{l \in L}$ for $K$. Then we have that
\begin{eqnarray*}
\sum_{i \in I} \mu(q_i^* q_i)
& = & \sum_{i \in I} \langle \th(q_i) v , \th(q_i) v \rangle
= \sum_{i \in I} \sum_{l \in L} \langle \th(q_i) v , f_l \rangle \, \, \langle f_l , \th(q_i) v \rangle \\
& = & \sum_{i \in I} \sum_{l \in L} | \om_{v,f_l}(q_i) |^2
= \sum_{l \in L} \sum_{i \in I} | \langle \lafi((\om_{v,f_l} \ot \io)(x)) , e_i \rangle|^2  \\
& = & \sum_{l \in L} \vfi((\om_{v,f_l} \ot \io)(x)^*
(\om_{v,f_l} \ot \io)(x)) \ .
\end{eqnarray*}
Therefore the previous lemma implies that
$$\sum_{i \in I} \mu(q_i^* q_i) = \vfi((\mu \ot \io)(x^* x))= \mu((\io \ot \vfi)(x^* x)) \ .$$
Hence lemma \ref{weight.lem5} implies that the net $(\,\sum_{i \in J} q_i^* q_i\,)_{J \in F(I)}$ converges strictly to $(\io \ot \vfi)(x^* x)$.

\medskip\medskip

Take finite subsets $J$ and $K$ of $I$ such that $J \subseteq K$ and $a \in A$. Then we have that
$$\| \,  \sum_{i \in K} q_i a \ot  e_i  - \sum_{i \in J} q_i a \ot e_i  \, \|^2  = \| \, \sum_{i \in K \setminus J} a^* q_i^* q_i a \, \| \ ,$$
implying that
the net $$(\, \sum_{i \in J} q_i a \ot e_i \,)_{J \in F(I)}$$ is Cauchy and hence convergent in $A \ot H_\vfi$.

So we can define a linear operator $F_x : A \rightarrow A \ot H_\vfi$ such that $F_x(a) = \sum_{i \in I} q_i a \ot e_i$ for $a \in A$.

Because $\sum_{i \in I} q_i^* q_i = (\io \ot \vfi)(x^* x)$ in the
strict topology, it follows that $\langle F_x(a) , F_x(a) \rangle =
a^* (\io \ot \vfi)(x^* x) a$ for $a \in A$. \inlabel{ster1}

\medskip

Choose $a \in A$, $b \in A$ and $c \in \Nfi$. Then we have for $\om \in A^*$ that
\begin{eqnarray*}
& & \om(\langle F_x(a) , b \ot \lafi(c) \rangle) =  \sum_{i \in I} \om(\langle q_i a \ot e_i  , b \ot \lafi(c) \rangle)
 = \sum_{i \in I} \om(b^* q_i a) \langle e_i , \lafi(c) \rangle \\
& & \spat =  \sum_{i \in I} \langle \lafi((a \om b^* \ot \io)(x)) , e_i \rangle \, \langle e_i , \lafi(c) \rangle
=  \langle \lafi((a \om b^* \ot \io)(x)) , \lafi(c) \rangle
= \vfi( c^* (a \om b^* \ot \io)(x)) \\
& & \spat =  (a \om b^*)\bigl( (\io  \ot \vfi)( (1 \ot c^*)x )  \bigr)
= \om( (\io \ot \vfi)((b^* \ot c^*)x)\,  a ) \ .
\end{eqnarray*}
So we see that $\langle F_x(a) , b \ot \lafi(c) \rangle
= (\io \ot \vfi)((b^* \ot c^*) x) \, a$. \inlabel{ster2}

\medskip

The Cauchy-Schwarz inequality in proposition \ref{weight.prop3} implies for $y \in A \odot \Nfi$ that
$$\| (\io \ot \vfi)(x^* y) \|^2 \leq \|(\io \ot \vfi)(x^* x)\| \,
\|(\io \ot \vfi)(y^* y)\| = \|(\io \ot \vfi)(x^* x)\| \,
\|(\io \odot \lafi)(y)\|^2$$
This implies that there exists a unique bounded linear map $G_x : A \ot H_\vfi \rightarrow A$ such that $G_x (b \ot \lafi(c)) = (\io \ot \vfi)(x^*(b \ot c))$ for $b \in A$ and $c \in \Nfi$. Equation~\ref{ster2} then implies that
$\langle F_x(a) , v \rangle = G_x(v)^* a$ for $a \in A$ and $v \in A \ot H_\vfi$. So $F_x$ belongs to $\cL(A,A \ot H_\vfi)$ and $F_x = G_x^*$.

\medskip

Now define the map $\la : \bar{\cN}_{\io \ot \vfi} \rightarrow \cL(A,A \ot H_\vfi)$ such that $\la(x) = F_x$ for $x \in \bar{\cN}_{\io \ot \vfi}$. So we have that $\la(x)^* (b \ot \lafi(c))  =
(\io \ot \vfi)(x^* (b \ot c))$ for $b \in A$, $c \in \Nfi$ and $x \in \bar{\cN}_{\io \ot \vfi}$.

This implies immediately that $\la$ is linear and that $\la(a \ot b) = a \ot \lafi(b)$ for $a \in M(A)$ and $b \in \bar{\cN}_\vfi$.

Equation~\ref{ster1} implies that $\la(x)^*  \la(x) = (\io \ot \vfi)(x^* x)$ for $x \in \bar{\cN}_{\io \ot \vfi}$, so polarization yields that
$\la(y)^* \la(x) = (\io \ot \vfi)(y^* x)$ for $x,y \in \bar{\cN}_{\io \ot \vfi}$.

\medskip

Choose $b \in A$ and $c \in \Nfi$. Then we have for $z \in A \ot B$ that
\begin{eqnarray*}
\|(\io \ot \lafi)(z(b \ot c))\|^2 & = & \| (\io \ot \lafi)(z(b \ot c))^* (\io \ot \lafi)(z(b \ot c)) \| = \| (\io \ot \vfi)((b \ot c)^* z^* z (b \ot c)) \| \\
&  \leq & \|z^* z\| \, \|(\io \ot \vfi)((b \ot c)^* (b \ot c))\|
= \| z \|^2 \, \|b\|^2 \, \vfi(c^* c) \ .
\end{eqnarray*}
So the linear mapping $A \ot B \rightarrow \cL(A,A \ot H_\vfi) : z \mapsto
(\io \ot \lafi)(z(b \ot c))$ is bounded. This is of course also true for the linear mapping $A \ot B \rightarrow \cL(A,A \ot H_\vfi) : z \mapsto
(\io \ot \pifi)(z)(b \ot \lafi(c))$.

It is easy to see that both mappings above agree on $A \odot B$ so they agree on $A \ot B$, i.e. $(\io \ot \lafi)(z(b \ot c)) =
(\io \ot \pifi)(z)(b \ot \lafi(c))$ for $z \in A \ot B$.

\medskip

Now choose $x \in \bar{\cN}_{\io \ot \vfi}$ and $z \in M(A \ot B)$. Take an approximate unit $(u_j)_{j \in J}$ for $A \ot B$. Then we have for $j \in J$ and  $b \in A$, $c \in \Nfi$ that
\begin{eqnarray*}
& & [(\io \ot \pifi)(u_j) (\io \ot \pifi)(z) (\io \ot \lafi)(x)]^* (b \ot \lafi(c)) = [(\io \ot \pifi)(u_j z) (\io \ot \lafi)(x)]^* (b \ot \lafi(c)) \\
& & \spat = (\io \ot \lafi)(x)^* (\io \ot \pifi)(z^* u_j^*)(b \ot \lafi(c))
= (\io \ot \lafi)(x)^* (\io \ot \lafi)(z^* u_j^*  (b \ot c)) \\
& & \spat = (\io \ot \vfi)(x^* z^* u_j^* (b \ot c))
= (\io \ot \lafi)(z x)^* (\io \ot \lafi)(u_j^* (b \ot c)) \\
& & \spat = (\io \ot \lafi)(z x)^* (\io \ot \pifi)(u_j^*) (b \ot \lafi(c))
= [(\io \ot \pifi)(u_j) (\io \ot \lafi)(z x)]^* (b \ot \lafi(c)) \ .
\end{eqnarray*}
Hence $(\io \ot \pifi)(u_j) (\io \ot \pifi)(z) (\io \ot \lafi)(x)
= (\io \ot \pifi)(u_j) (\io \ot \lafi)(z x)$ for all $j \in J$,
so $(\io \ot \pifi)(z) (\io \ot \lafi)(x) = (\io \ot \lafi)(z x)$.
\end{demo}

\bigskip

\begin{remark} \rm \label{slice.rem1}
Take $a,b \in \bar{\cN}_{\io \ot \vfi}$. Then we have for all $x \in M(A \ot B)$ that $(\io \ot \vfi)(b^*\,x\,a)$ \newline $= (\io \ot \lafi)(b)^* (\io \ot \pifi)(x)\,(\io \ot \lafi)(a)$. So the linear mapping
$$M(A \ot B) \rightarrow M(A) : x \rightarrow (\io \ot \vfi)(b^*\,x\,a)$$
is bounded and strictly continuous on bounded sets.

\medskip

This implies immediately for all $c,d \in \bar{\cN}_\vfi$ and
$x \in M(A \ot B)$ that $$(\io \ot \vfi)((1 \ot d^*)\,x\,(1 \ot c))
= (\io \ot c \, \vfi \, d^*)(x) \ .$$
\end{remark}

\bigskip

In notation \ref{weight.not1}, we introduced the operators $T_\om \in B(H_\vfi)$ which where used to cut off $\lafi$ and get a continuous linear mapping in this way. We can also use them to cut off $\io \ot \lafi$.

First we formulate a simple lemma.

\begin{lemma} \label{weight.lem6}
We have for all $a,b \in A$, $v \in H_\vfi$, $\om \in A^*$ and $x \in \bar{\cN}_{\io \ot \vfi}$ that
$$\om(\langle (\io \ot \lafi)(x) \, a , b \ot v \rangle )
= \langle \lafi\bigl((a \om b^* \ot \io)(x)\bigr) , v \rangle \ .$$
\end{lemma}
\begin{demo}
It is enough to prove the lemma for a dense set of elements $v$. So take $c \in \Nfi$. Then we have that
\begin{eqnarray*}
 \om(\langle (\io \ot \lafi)(x) \, a , b \ot \lafi(c) \rangle )
& = & \om\bigl(\, (\io \ot \vfi)((b^* \ot c^*) x) \, a \, \bigr)
 = \vfi\bigl(\,(a \om \ot \io)((b^* \ot c^*) x)\,\bigr) \\
& = & \vfi\bigl(\,c^* (a \om b^* \ot \io)(x)\,\bigr)
= \langle \lafi\bigl((a \om b^* \ot \io)(x)\bigr) , \lafi(c) \rangle  \ .
\end{eqnarray*}
\end{demo}

\medskip

For every $v \in H_\vfi$, we define $\th_v \in  B(\C,H_\vfi)$ such that $\th_v(c) = c\,v$ for all $c \in \C$.

\begin{result} \label{cutoff}
Consider $\om \in \cF_\vfi$. Then we have for all $x \in \bar{\cN}_{\io \ot \vfi}$ that
$$(1 \ot T_\om^{\frac{1}{2}})(\io \ot \lafi)(x) = (\io \ot \pifi)(x) (1  \ot \th_{\xi_\om}).$$
\end{result}
\begin{demo}
Take $v \in H_\vfi$ , $b \in A$ and $\om \in A^*$.
Then the previous lemma implies that
\begin{eqnarray*}
& & \om(\langle (1 \ot T_\om^{\frac{1}{2}})(\io \ot \lafi)(x) \, a , b \ot  v \rangle) =
\om(\langle (\io \ot \lafi)(x) \, a , b \ot T_\om^{\frac{1}{2}} v \rangle) \\
& & \spat = \langle \lafi((a \om b^* \ot \io)(x)) , T_\om^{\frac{1}{2}} v  \rangle  = \langle T_\om^{\frac{1}{2}} \lafi((a \om b^* \ot \io)(x)) , v \rangle \\
& & \spat = \langle  \pifi((a \om b^* \ot \io)(x)) \, \xi_\om , v \rangle = \om(\langle (\io \ot \pifi)(x) (a \ot \xi_\om) , b \ot v \rangle)\ .
\end{eqnarray*}
So we see that $(1 \ot T_\om^{\frac{1}{2}})(\io \ot \lafi)(x) \, a =
(\io \ot \pifi)(x) (a \ot \xi_\om)$ for $a \in A$.
\end{demo}

\begin{corollary} \label{weight.cor1}
Consider  $\om \in \cF_\vfi$. Then we have for all $x,y \in \bar{\cN}_{\io \ot \vfi}$ that
$$(\io \ot \lafi)(y)^* (1 \ot T_\om) (\io \ot \lafi)(x)
= (\io \ot \om)(y^* x)\ .$$
\end{corollary}

\bigskip

In proposition \ref{weight.prop5}, we mentioned that the mapping $\lafi$ is closed if $\vfi$ is lower semi-continuous. A similar property holds for $\io \ot \lafi$.

\begin{proposition} \label{weight.prop7}
The linear mapping
$\bar{\cN}_{\io \ot \vfi}  \rightarrow \cL(A,A \ot H) : x \mapsto (\io \ot \lafi)(x)$ is closed for the strict topology on $M(A \ot B)$ and the strong topology on $\cL(A,A \ot H)$.
\end{proposition}
\begin{demo}
Choose a net $(a_i)_{i \in I}$ in $\bar{\cN}_{\io \ot \vfi}$, $a \in M(A \ot B)$ and $S \in \cL(A, A \ot H_\vfi)$ such that $(a_i)_{i \in I}$ converges to $a$ in the strict topology and such that
$\bigl((\io \ot \lafi)(a_i)\bigr)_{i \in I}$ converges strongly to $S$.

Choose $b \in \bar{\cN}_{\io \ot \vfi}$, $c,d \in A$ and $\om \in \cF_\vfi$. By corollary \ref{weight.cor1}, we have for all $i \in I$ that
$$\langle (1 \ot T_\om)(\io \ot \lafi)(a_i) \, c ,  (\io \ot \lafi)(b) \, d \rangle = d^* (\io \ot \om)(b^* a_i) \, c \, .$$
This implies that the net $(\, \langle (1 \ot T_\om)(\io \ot \lafi)(a_i) \, c ,  (\io \ot \lafi)(b) \, d \rangle  \,)_{i \in I}$ converges to $d^* (\io \ot \om)(b^* a) \, c$. It is also clear the net $(\, \langle (1 \ot T_\om)(\io \ot \lafi)(a_i) \, c ,  (\io \ot \lafi)(b) \, d \rangle\,)_{i \in I}$ converges to
$\langle (1 \ot T_\om) S c , (\io \ot \lafi)(b) \, d \rangle$.

Combining these two facts, we get that
\begin{equation} \label{vgl2}
\langle (1 \ot T_\om) \, S c , (\io \ot \lafi)(b)\, d \rangle
= d^* (\io \ot \om)(b^* a) \, c
\end{equation}
for all $b \in \bar{\cN}_{\io \ot \vfi}$, $c,d \in A$ and $\om \in \cF_\vfi$.
\begin{enumerate}
\item Choose $c \in A$. Fix $\om \in \cG_{\vfi}$ for the moment.
We have immediately that the net
$(c^* (\io \ot\om)(a_i^* a) \, c)_{i \in I}$ converges to $c^* (\io \ot \om)(a^* a) \, c$.
Equality~\ref{vgl2} implies on the other hand that the net
$(c^* (\io \ot\om)(a_i^* a) \, c)_{i \in I}$ converges to $\langle (1 \ot T_\om) S c , S c \rangle$.
So we get that
$c^* (\io \ot \om)(a^* a) \, c = \langle (1 \ot T_\om) S c , S c \rangle$.

This implies that the net $(c^* (\io \ot \om)(a^* a) \, c)_{\om \in \cG_\vfi}$ converges to $\langle S c , S c \rangle$.
Hence result \ref{weight.res3} implies that $a \in \bar{\cN}_{\io \ot \vfi}$.

\item Choose $b \in A \od \Nfi$ and $c,d \in A$. By equality~\ref{vgl2} and corollary \ref{weight.cor1}, we have for all
$\om \in \cG_\vfi$ that
$$\langle (1 \ot T_\om) S(c) , (\io \ot \lafi)(b) \, d \rangle
= d^* (\io \ot \om)(b^* a) \, c = \langle (1 \ot T_\om) (\io \ot \lafi)(a) \, c , (\io \ot \lafi)(b) \, d \rangle \ . $$
Because $(T_\om)_{\om \in \cG_\vfi}$ converges strongly to 1, this implies that
$\langle S(c) , (\io \ot \lafi)(b) \, d \rangle
=$ \newline $\langle (\io \ot \lafi)(a) \, c , (\io \ot \lafi)(b) \, d \rangle$.
Hence $S = (\io \ot \lafi)(a)$.
\end{enumerate}
\end{demo}

\bigskip

In a next step, we quickly prove some kind of dominated convergence property for the \lq\cst-valued weight\rq \ $\io \ot \vfi$.

\begin{proposition} \label{weight.prop8}
Consider  a net $(x_i)_{i \in I}$ in $\bar{\cM}_{\io \ot \vfi}^+$ and an element $x$ in $\bar{\cM}_{\io \ot \vfi}^+$. If $(x_i)_{i \in I}$ converges strictly to $x$ and $x_i \leq x$ for every $i \in I$,
then $\bigl((\io \ot \vfi)(x_i)\bigr)_{i \in I}$ converges strictly to $(\io \ot \vfi)(x)$.
\end{proposition}
\begin{demo}
It is clear that $0 \leq (\io \ot \vfi)(x_i) \leq (\io \ot \vfi)(x)$ for every $i \in I$. Choose $a \in A$.

Take $\vep > 0$. Then there exist $\om \in {\cal G}_\vfi$ such that
$\| a^* (\io \ot \vfi)(x) \, a - a^* (\io \ot \om)(x) \, a \| \leq \frac{\vep}{2}$.

Because $((\io \ot \om)(x_i))_{i \in I}$ converges strictly to $(\io \ot \om)(x)$, there
exists an element $i_0 \in I$ such that \newline
$\|a^* (\io \ot \om)(x_i)\, a - a^* (\io \ot \om)(x)\,  a\| \leq \frac{\vep}{2}$ for every $i \in I$ with $i \geq i_0$.

Choose $j \in I$ with $j \geq i_0$. Then
\begin{eqnarray*}
& & \| a^* (\io \ot \om)(x_j) \, a - a^* (\io \ot \vfi)(x) \, a \| \\
& & \spat \leq \| a^* (\io \ot \om)(x_j) \, a - a^* (\io \ot \om)(x) \, a \| + \| a^* (\io \ot \om)(x) \, a - a^* (\io \ot \vfi)(x)\, a \|
\leq \frac{\vep}{2} + \frac{\vep}{2} = \vep .
\end{eqnarray*}
Because $0 \leq a^* (\io \ot \om)(x_j) \, a \leq a^* (\io \ot \vfi)(x_j) \, a \leq a^* (\io \ot \vfi)(x) \, a $, we get that
$$\| a^* (\io \ot \vfi)(x_j) \, a - a^* (\io \ot \vfi)(x) \, a \|
\leq \| a^* (\io \ot \om)(x_j) \, a - a^* (\io \ot \vfi)(x) \, a \| \leq \vep \ .$$

Consequently, we see that $(a^* (\io \ot \vfi)(x_i) \,a)_{i \in I}$ converges to $a^* (\io \ot \vfi)(x)\,a$. Lemma \ref{weight.lem4} implies now that $\bigl((\io \ot \vfi)(x_i)\bigr)_{i \in I}$ converges strictly to $(\io \ot \vfi)(x)$.
\end{demo}

\medskip

\begin{proposition} \label{weight.prop9}
Consider   a net $(x_i)_{i \in I}$ in $\bar{\cM}_{\io \ot \vfi}^+$ and an element $x$ in $M(A \ot B)^+$. Suppose that $(x_i)_{i \in I}$ converges strictly to $x$ and $x_i \leq x$ for every $i \in I$.
Then $x$ belongs to $\bar{\cM}_{\io \ot \vfi}^+$ $\Leftrightarrow$ The net $(a^* (\io \ot \vfi)(x_i) \,a)_{i \in I}$ is convergent for every $a \in A$.
\end{proposition}
\begin{demo}
One implication follows from the previous result, we will turn to the other one. Therefore assume that the net $(a^* (\io \ot \vfi)(x_i) \, a)_{i \in I}$ is convergent for every $a \in A$. Choose $c \in A$.

\medskip

By assumption, there exists an element $d \in A^+$ such that
$(c^* (\io \ot \vfi)(x_i) \, c)_{i \in I}$ converges to $d$.
\inlabel{lettera}

First, we will prove that $c^* (\io \ot \om)(x) \, c \leq d $ for every $\om \in  \cG_\vfi$.
\inlabel{letterb}

\begin{list}{}{\setlength{\leftmargin}{.4 cm}}

\item Take $\om \in \cG_\vfi$. Choose $n \in \N$.
From statement~\ref{lettera}, we get the existence of an element $i_0 \in I$
such that $\| d - c^* (\io \ot \vfi)(x_i) \, c \| \leq \frac{1}{n}$ for every $i \in I$ with $i \geq i_0$.
Therefore, we have for every $i \in I$ with $i \geq i_0$
that $c^* (\io \ot \om)(x_i) \, c \leq c^* (\io \ot \vfi)(x_i) \, c \leq d + \frac{1}{n} \, 1$

Because $(c^* (\io \ot \om)(x_i) \, c)_{i \in I}$ converges to $c^* (\io \ot \om)(x) \, c$,
this implies that $c^* (\io \ot \om)(x) \, c \leq d + \frac{1}{n} \, 1$.

If we let $n$ tend to $\infty$, we get that $c^* (\io \ot \om)(x) \, c \leq d$.

\end{list}

Choose $\vep > 0$.
By statement~\ref{lettera}, there exists an element $j \in I$ such that
$\| c^* (\io \ot \vfi)(x_j) \, c - d \| \leq \frac{\vep}{2}$.

We have furthermore an element $\om_0 \in \cG_\vfi$ such that $\| c^* (\io \ot \vfi)(x_j) \, c - c^* (\io \ot \om)(x_j) \, c \| \leq \frac{\vep}{2}$ for every $\om \in \cG_\vfi$ with $\om \geq \om_0$.
So we get that $\|c^* (\io \ot \om)(x_j) \, c - d \| \leq \vep$ for every
$\om \in \cG_\vfi$ with $\om \geq \om_0$. \inlabel{letterc}

\smallskip

Take $\eta \in \cG_\vfi$ with $\eta \geq \om_0$.
Using inequality~\ref{letterb}, we have that $0 \leq c^* (\io \ot \eta)(x_j) \, c \leq c^* (\io \ot \eta)(x) \, c \leq d$. Consequently inequality~\ref{letterc} implies that
$$\| c^* (\io \ot \eta)(x) \, c - d \| \leq \| c^* (\io \ot \eta)(x_j) \, c - d \| \leq \vep \ . $$

\medskip

Hence, we see that $(c^* (\io \ot \om)(x) \, c)_{\om \in \cG_\vfi}$ converges to $d$. Therefore $x$ belongs to $\bar{\cM}_{\io \ot \vfi}^+$.
\end{demo}

\medskip

\begin{proposition}
Consider  a $^*$-isomorphism $\th : A \rightarrow A$ such that there exists a strictly positive number $r > 0$ such that $\vfi \, \th = r \, \vfi$.  Then the following holds:
\begin{enumerate}
\item We have for all $x \in  \bar{\cM}_{\io \ot \vfi}$ that
$(\io \ot \th)(x)$ belongs to $\bar{\cM}_{\io \ot \vfi}$ and
$(\io \ot \vfi)\bigl((\io \ot \th)(x)\bigr) = r \, (\io \ot \vfi)(x)$ \ .
\item Define the unitary operator $U$ on $H_\vfi$ such that $U \lafi(a) = r^{-\frac{1}{2}} \, \lafi(\th(a))$ for all $a \in \Nfi$. Then we have for all $x \in \bar{\cN}_{\io \ot \vfi}$ that $(\io \ot \th)(x)$ belongs to $\bar{\cN}_{\io \ot \vfi}$ and $(1 \ot U)(\io \ot \lafi)(x) =$ \newline $ r^{-\frac{1}{2}} \, (\io \ot \lafi)\bigl((\io \ot \th)(x)\bigr)
$\ .
\end{enumerate}
\end{proposition}
\begin{demo}
\begin{enumerate}
\item It is again sufficient to prove the result in the case that $x \geq 0$.  Take $\om \in A^*_+$. By proposition \ref{weight.prop1}, $(\om \ot \io)(x)$ belongs to $\bar{\cM}_\vfi^+$. So $\th\bigl((\om \ot \io)(x)\bigr)$ belongs to
$\bar{\cM}_\vfi^+$. Hence $(\om \ot \io)\bigl((\io \ot \th)(x)\bigr)$ belongs to $\bar{\cM}_\vfi^+$ and
$$\vfi\bigl(\,(\om \ot \io)\bigl((\io \ot \th)(x)\bigr)\,\bigr) = \vfi\bigl(\,\th((\om \ot \io)(x)) \, \bigr)
= r \, \vfi((\om \ot \io)(x)) = r \, \om((\io \ot \vfi)(x)) \ .$$
Hence proposition \ref{weight.prop5a} implies that  $(\io \ot \th)(x)$ belongs to $\bar{\cM}_{\io \ot \vfi}^+$ and
$(\io \ot \vfi)\bigl((\io \ot \th)(x)\bigr) = r \, (\io \ot \vfi)(x)$.
\item Choose $x \in \bar{\cN}_{\io \ot \vfi}$. Then we have for all $a \in A$ and $b \in \Nfi$ that
\begin{eqnarray*}
& & [(1 \ot U)(\io \ot \lafi)(x)]^* (a \ot \lafi(b))
= (\io \ot \lafi)(x)^* (a \ot U^* \lafi(b)) \\
& & \spat = r^{\frac{1}{2}} \, (\io \ot \lafi)(x)^* (\,a \ot \lafi(\th^{-1}(b))\,)
= r^{\frac{1}{2}} \, (\io \ot \vfi)(\,x^* (a \ot \th^{-1}(b))\,) \ .
\end{eqnarray*}
Therefore the result proven in the first part implies that
\begin{eqnarray*}
& & [(1 \ot U)(\io \ot \lafi)(x)]^* (a \ot \lafi(b))
= r^{-\frac{1}{2}} \, (\io \ot \vfi)\bigl(\,(\io \ot \th)(x^* (a \ot \th^{-1}(b)))\,\bigr) \\
& & \spat = r^{-\frac{1}{2}} \, (\io \ot \vfi)((\io \ot \th)(x)^* (a \ot b))
= r^{-\frac{1}{2}} \, (\io \ot \lafi)((\io \ot \th)(x))^* (a \ot \lafi(b)) \ .
\end{eqnarray*}
Consequently, $(1 \ot U)(\io \ot \lafi)(x) = r^{-\frac{1}{2}} \, (\io \ot \lafi)((\io \ot \th)(x))$.
\end{enumerate}
\end{demo}

\medskip

The next right $M(A)$-module property is also an easy consequence of the definition of the slice weight $\io \ot \vfi$.

\begin{proposition}
The following properties hold:
\begin{enumerate}
\item We have for all $x \in \bar{\cM}_{\io \ot \vfi}$ and $a,b \in M(A)$ that  $(b^* \ot 1) \,x\, (a \ot 1)$ belongs to $\bar{\cM}_{\io \ot \vfi}$ and \newline $(\io \ot \vfi)\bigl((b^* \ot 1)\, x\, (a \ot 1)\bigr) = b^*  (\io \ot \vfi)(x) \, a$.
\item We have for all $x \in \bar{\cN}_{\io \ot \vfi}$ and $a \in M(A)$ that $x\,(a \ot 1)$ belongs to $\bar{\cN}_{\io \ot \vfi}$ and
$(\io \ot \lafi)(x \,(a \ot 1)) = (\io \ot \lafi)(x) \, a$.
\end{enumerate}
\end{proposition}
\begin{demo}
\begin{enumerate}
\item As usual, we can assume that $x \geq 0$ and that $b = a$. By definition, we know that the net \newline $(\,(\io \ot \om)(x)\,)_{\om \in \cG_\vfi}$ converges strictly to $(\io \ot \vfi)(x)$. Hence $(\,a^* (\io \ot \om)(x)\,a\,)_{\om \in \cG_\vfi}$ converges strictly to $a^* (\io \ot \vfi)(x) \, a$.

Consequently, the net $\bigl(\,(\io \ot \om)((a^* \ot 1)\, x\,(a \ot 1))\,\bigr)_{\om \in \cG_\vfi}$ converges strictly to $a^* (\io \ot \vfi)(x) \, a$. By definition, this implies that $(a^* \ot 1)\, x \,(a \ot 1)$ belongs to $\bar{\cM}_{\io \ot \vfi}^+$ and
$(\io \ot \vfi)((a^* \ot 1) \,x\, (a \ot 1)) = a^* (\io \ot \vfi)(x)\, a$.
\item Take $b \in A$ and $c \in \Nfi$. Then
\begin{eqnarray*}
& &  (\io \ot \lafi)(x\,(a \ot 1))^* (b \ot \lafi(c))
= (\io \ot \vfi)((a^* \ot 1)\, x^* (b \ot c)) \\
& & \spat = a^*  (\io \ot \vfi)(x^* (b \ot c))
= a^* (\io \ot \lafi)(x)^* (b \ot \lafi(c)) \ .
\end{eqnarray*}
Hence  $(\io \ot \lafi)(x(a \ot 1))^* = a^* (\io \ot \lafi)(x)^*$.
\end{enumerate}
\end{demo}

\medskip

If $\vfi$ is a KMS weight, we can also prove some right $M(B)$-module properties and get a version of proposition \ref{weight.prop6}.

\begin{proposition}
Suppose that $\vfi$ is a KMS weight on $B$ with modular group $\si$. Let $J$ denote the modular conjugation of $\vfi$ in the GNS-construction $(H_\vfi,\pifi,\lafi)$. Then
\begin{enumerate}
\item We have for all $x \in \bar{\cN}_{\io \ot \vfi}$ and $b \in D(\overline{\si}_{\frac{i}{2}})$ that $x\, (1 \ot b)$ belongs to
$\bar{\cN}_{\io \ot \vfi}$ and $(\io \ot \lafi)(x\,(1 \ot b))
= (1 \ot J \pi(\si_{\frac{i}{2}}(b))^* J)(\io \ot \lafi)(x)$.
\item Let $x \in \bar{\cM}_{\io \ot \vfi}$ and $b \in D(\overline{\si}_i)$. Then $x \, (1 \ot b)$ belongs to $\bar{\cM}_{\io \ot \vfi}$, $(1 \ot \si_i(b))\, x$ belongs to $\bar{\cM}_{\io \ot \vfi}$ and
$(\io \ot \vfi)\bigl(x\,(1 \ot b)\bigr) = (\io \ot \vfi)\bigl((1 \ot \si_i(b))\, x\bigr)$.
\item Let $x \in \bar{\cN}_{\io \ot \vfi} \cap \bar{\cN}^*_{\io \ot \vfi}$ and $b \in D(\overline{\si}_i) \cap \bar{\cN}_{\vfi}$ such that $\si_i(b) \in \bar{\cN}^*_{\vfi}$. Then
$(\io \ot \vfi)\bigl(x\,(1 \ot b)\bigr) = (\io \ot \vfi)\bigl((1 \ot \si_i(b))\, x\bigr)$.
\end{enumerate}
\end{proposition}
\begin{demo}
\begin{enumerate}
\item Take $c \in \Nfi \cap D(\si_{\frac{i}{2}})$. Then $x(1 \ot c)$ belongs certainly to $\bar{\cN}_{\io \ot \vfi}$.

Choose $d,e \in A$, $v \in H_\vfi$ and $\om \in A^*$. Then lemma \ref{weight.lem6} and proposition \ref{extmult} imply that
\begin{eqnarray*}
& & \om(\langle (1 \ot J \pifi(\si_{\frac{i}{2}}(c))^* J) (\io \ot \lafi)(x) \, d , e \ot v \rangle)
= \om(\langle (\io \ot \lafi)(x) \, d , e \ot J \pifi(\si_{\frac{i}{2}}(c)) J v \rangle) \\
& & \spat = \langle \lafi\bigl((d \om e^* \ot \io)(x)\bigr) , J \pifi(\si_{\frac{i}{2}}(c)) J v \rangle
= \langle J \pifi(\si_{\frac{i}{2}}(c))^* J \lafi\bigl((d \om e^* \ot \io)(x)\bigr) , v \rangle \\
& & \spat = \langle  \lafi\bigl((d \om e^* \ot \io)(x) \, c \bigr) , v \rangle
= \langle  \lafi\bigl((d \om e^* \ot \io)(x(1 \ot c))  \bigr) , v \rangle \\
& & \spat = \om(\langle (\io \ot \lafi)(x(1 \ot c)) \, d , e \ot v \rangle)
\ .
\end{eqnarray*}
Hence $ (1 \ot J \pifi(\si_{\frac{i}{2}}(c))^* J) (\io \ot \lafi)(x)
= (\io \ot \lafi)(x(1 \ot c))$.

\smallskip

There exists a net $(e_j)_{j \in J}$ in $\Nfi \cap D(\si_{\frac{i}{2}})$ such that
\begin{itemize}
\item $(e_j)_{j \in J}$ is bounded and converges strictly to $1$
\item $(\si_{\frac{i}{2}}(e_j))_{j \in J}$ is bounded and converges strictly to $1$.
\end{itemize}
Cfr. the proof of proposition \ref{sweight.prop1}.

For $j \in J$, we define $c_j = b \, e_j \in \Nfi \cap D(\si_{\frac{i}{2}})$. Then
\begin{itemize}
\item $(c_j)_{j \in J}$ is bounded and converges strictly to $b$
\item $(\si_{\frac{i}{2}}(c_j))_{j \in J}$ is bounded and converges strictly to $\si_{\frac{i}{2}}(b)$.
\end{itemize}

By the first part, we have for every $j \in J$ that $x (1 \ot c_j)$ belongs to $\bar{\cN}_{\io \ot \vfi}$ and $(\io \ot \lafi)(x(1 \ot c_j))
= (1 \ot J \pifi(\si_{\frac{i}{2}}(c_j))^* J)(\io \ot \lafi)(x)$.

So $\bigl(\,x\,(1 \ot c_j )\,\bigr)_{j \in J}$ converges strictly to
$x\,(1 \ot b)$ and $\bigl(\,(\io \ot \lafi)(x\,(1 \ot c_j))\,\bigr)_{j \in J}$
converges strongly to $(1 \ot J \pifi(\si_{\frac{i}{2}}(b))^* J)(\io \ot \lafi)(x)$. Therefore proposition \ref{weight.prop7} implies that
$x\,(1 \ot b)$ belongs to $\bar{\cN}_{\io \ot \vfi}$ and
$(\io \ot \lafi)(x\,(1 \ot b)) = (1 \ot J \pifi(\si_{\frac{i}{2}}(b))^* J)(\io \ot \lafi)(x)$.

\item By the first part of this proof, it follows easily that $x \, (1 \ot b)$ and $(1 \ot \si_i(b))\, x$ belong to $\bar{\cM}_{\io \ot \vfi}$.
Using propositions  \ref{weight.prop1} and \ref{weight.prop6}.3, we get for all $\th \in A^*$ that
$$\th\bigl((\io \ot \vfi)((1 \ot \si_i(b))\, x)\bigr)
= \vfi(\si_i(b)\,(\th \ot \io)(x)) = \vfi((\th \ot \io)(x)\,b)
= \th\bigl((\io \ot \vfi)(x\,(1 \ot b))\bigr) \ .$$
Therefore $(\io \ot \vfi)((1 \ot \si_i(b))\, x)
= (\io \ot \vfi)(x\,(1 \ot b))$.
\item Similar to the proof of 2.
\end{enumerate}
\end{demo}

\bigskip\bigskip

Up to now, we mainly worked with the strict topology. In a last short part, we will focus on the norm topology. In connection with this, we will also look at the set of elements $x \in \bar{\cM}_{\io \ot \vfi}^+$ for which $(\io \ot \vfi)(x)$ belongs to $A \ot B$.

\bigskip

The proof of the following form of lower semi-continuity in the norm topology is completely analogous to the proof of proposition \ref{weight.prop8}.

\begin{proposition}
Consider a net $(x_i)_{i \in I}$ in $\bar{\cM}_{\io \ot \vfi}^+$ and an element $x$ in $\bar{\cM}_{\io \ot \vfi}^+$. If $(x_i)_{i \in I}$ converges to $x$ and $x_i \leq x$ for every $i \in I$,
then $\bigl((\io \ot \vfi)(x_i)\bigr)_{i \in I}$ converges  to $(\io \ot \vfi)(x)$.
\end{proposition}

\medskip

\begin{lemma}
Consider a Hilbert \cst-module $E$ over $A$ and  an element $T$ in $\cL(A,E)$. Then $T^* T$ belongs to $A$
$\Leftrightarrow$
There exists an element $v \in E$ such that $T a = v a $ for every $a \in A$
\end{lemma}
\begin{demo}
\begin{itemize}
\item Suppose that $T^* T$ belongs to $A$. Define $S \in \cL(A \oplus E , A \oplus E)$  with $S^* = S$ such that
$$ S = \left( \begin{array}{cc} 0 & T^* \\ T  & 0 \end{array} \right) .$$
We have that
$$ S^2 = \left( \begin{array}{cc} T^* T & 0 \\ 0 & T T^* \end{array} \right) .$$
By proposition 1.4.5 of \cite{Ped}, there exists an element $U \in \cL(A \oplus E, A \oplus E)$ such that $S = U (S^2)^\frac{1}{4}$.
This implies that
$$  \left( \begin{array}{cc} 0 & T^* \\ T  & 0 \end{array} \right)
 = U \, \left( \begin{array}{cc} (T^* T)^\frac{1}{4} & 0 \\ 0 &
(T T^*)^\frac{1}{4} \end{array} \right).$$
From this, we get the existence of $V \in \cL(A,E)$ such that $T = V (T^* T)^\frac{1}{4}$.
Because $T^* T$ belongs to $A$, we have that $(T^* T)^\frac{1}{4}$ belongs to $A$. Put $v = V\bigl((T^* T)^\frac{1}{4}\bigr) \in E$. Then $T(a) = v a$ for every $v \in A$.
\item Suppose there exists an element $v \in E$ such that $T(a) = v a$ for every $a \in A$. In this case, we have for every $a \in A$ that
$$ a^* (T^* T) a = \langle (T^* T) a ,  a \rangle = \langle T a , T a \rangle = \langle v a , v a \rangle = a^* \langle v , v \rangle \, a $$
which implies that $T^* T = \langle v , v \rangle \in A$.
\end{itemize}
\end{demo}

\medskip

\begin{proposition} \label{weight.prop10}
Let $x$ be an element in $\bar{\cN}_{\io \ot \vfi}$. Then $(\io \ot \vfi)(x^* x)$ belongs to $A$
\vspace{-0.2cm}
\begin{trivlist}
\item[\ $\Leftrightarrow$] We have for every $\om \in \cG_{\vfi}$ that $(\io \ot \om)(x^* x)$ belongs to $A$ and $\bigl(\,(\io \ot \om)(x^* x)\,\bigr)_{\om \in \cG_\vfi}$ converges to \newline  $\text{\ \ \ \ } (\io \ot \vfi)(x^* x)$
\item[\ $\Leftrightarrow$] There exists $v \in A \ot H_\vfi$ such that
$(\io \ot \lafi)(x) \, a = v a$ for $a \in A$.
\end{trivlist}
\end{proposition}
\begin{demo}
\begin{trivlist}
\item[\  (1) $\Rightarrow$ (3)] This follows from the previous lemma because
$(\io \ot \lafi)(x)^* (\io \ot \lafi)(x) = (\io \ot \vfi)(x^* x)$.
\item[\  (3) $\Rightarrow$ (2)] By corollary \ref{weight.cor1}, we have for all $\om \in \cG_\vfi$ \newline that
$$ (\io \ot \om)(x^* x) = (\io \ot \lafi)(x)^* (1 \ot T_\om) (\io \ot \lafi)(x) = \langle (1 \ot T_\om)\,v , v \rangle \ .$$
This implies immediately that $(\io \ot \om)(x^* x) \in A$ for $\om \in \cG_\vfi$ and that $\bigl(\,(\io \ot \om)(x^* x)\,\bigr)_{\om \in \cG_\vfi}$ converges to $\langle v , v \rangle = (\io \ot \vfi)(x^* x)$.
\item[\  (2) $\Rightarrow$ (1)] Trivial.
\end{trivlist}
\end{demo}

\medskip

Consider a \cst-algebra $A$. Let $a$ be an element in $A^+$ and $b$ an element in $M(A)^+$ such that
$b \leq a$. By proposition 1.4.5 of \cite{Ped}, there exists an element $U \in M(A)$ such that $b^\frac{1}{2} = U  a^\frac{1}{4}$ which implies that
$b^\frac{1}{2}$ belongs to $A$. Therefore $b \in A^+$.
This argument shows that $A^+$ is hereditary in $M(A)^+$.

\medskip

Let us introduce some natural sets.

\begin{definition}
We define the following sets :
\begin{itemize}
\item $\cM_{\io \ot \vfi}^+ = \{ \ x \in (A \ot B)^+  \mid x \in \bar{\cM}_{\io \ot \vfi}^+ \text{ and } (\io \ot \vfi)(x) \in A^+ \,  \}$\ , then $\cM_{\io \ot \vfi}^+$ is a hereditary cone in $(A \ot B)^+$.
\item $\cN_{\io \ot \vfi} = \{\, x \in A \ot B \mid x^* x \in \cM_{\io \ot \vfi}^+ \, \}$\ , so $\cN_{\io \ot \vfi}$ is a left ideal in $M(A \ot B)$.
\item $\cM_{\io \ot \vfi} = \text{span } \cM_\vfi^+ = \cN_{\io \ot \vfi}^* \, \cN_{\io \ot \vfi}$\ , so $\cM_{\io \ot \vfi}$ is a sub-$^*$-algebra of $A \ot B$.
\end{itemize}
\end{definition}

\medskip

There is of course also an alternative way to get to the set $\cM_{\io \ot \vfi}^+$. The next two results follow immediately from proposition \ref{weight.prop10}

\begin{result}
We have for all $x \in \cM_{\io \ot \vfi}^+$  that
the net $\bigl((\io \ot \om)(x)\bigr)_{\om \in \cG_\vfi}$ converges to
$(\io \ot \vfi)(x)$.
\end{result}

\begin{result}
Let $x \in (A \ot B)^+$. Then $x$ belongs to $\cM_{\io \ot \vfi}^+$ $\Leftrightarrow$
the net $\bigl((\io \ot \om)(x)\bigr)_{\om \in \cG_\vfi}$ is norm convergent in
$A$.
\end{result}

\medskip

Proposition \ref{weight.prop10} implies of course also the next two results.

\begin{result}
Consider  $x \in \cN_{\io \ot \vfi}$. Then there exists an element $v \in A \ot H_\vfi$ such that
$(\io \ot \lafi)(x) \, a = v \, a$ for $a \in A$ and $(\io \ot \vfi)(x^* x) = \langle v , v \rangle$.
\end{result}

\begin{result}
Consider  $x \in (A \ot B) \, \cap \, \bar{\cN}_{\io \ot \vfi}$. Then $x \in \cN_{\io \ot \vfi}$ $\Leftrightarrow$ There exists an element $v \in A \ot H_\vfi$ such that $(\io \ot \lafi)(x) \, a = v \, a$ for $a \in A$.
\end{result}

\medskip

Combining these last results with proposition \ref{weight.prop7}, it is easy to prove the following one.

\begin{corollary}
The mapping
$\cN_{\io \ot \vfi}  \rightarrow A \ot H_\vfi : x \mapsto (\io \ot \lafi)(x)$ is closed for the strict topology on $A \ot B$ and the norm topology on $A \ot H_\vfi$.
\end{corollary}

\bigskip\bigskip

We end this section with a natural generalization of proposition \ref{weight.prop1}.

\medskip

\begin{proposition}
Consider \cst-algebras $A$,$B$ and $C$ and a strict completely positive mapping
$\rho : A \rightarrow M(B)$. Let $\vfi$ be a proper weight on $C$ with GNS-construction $(H_\vfi,\pifi,\lafi)$. Then the following holds:
\begin{enumerate}
\item We have for all $x \in \bar{\cM}_{\io_A \ot \vfi}$ that
$(\rho \ot \io)(x)$ belongs to $\bar{\cM}_{\io_B \ot \vfi}$
and $(\io_B \ot \vfi)\bigl((\rho \ot \io)(x)\bigr)
=$ \newline $\rho\bigl((\io_A \ot \vfi)(x)\bigr)$ \ .
\item Consider $x \in \bar{\cN}_{\io_A \ot \vfi}$. Then $(\rho \ot \io)(x)$ belongs to $\bar{\cN}_{\io_B \ot \vfi}$ and
$(1 \ot \th_v^*)(\io_B \ot \lafi)\bigl((\rho \ot \io)(x)\bigr)
=$ \newline  $\rho\bigl(\, (1 \ot \th_v^*)(\io_A \ot \lafi)(x)\, \bigr)$
for every $v \in H_\vfi$.
\end{enumerate}
\end{proposition}
\begin{demo}
\begin{enumerate}
\item It is clear that it is enough to prove the result in the case that
$x \geq 0$. But then it is a straightforward application of definition \ref{weight.def2}. By this definition, we have that
$\bigl(\,(\io_A \ot \om)(x)\,\bigr)_{\om \in \cG_\vfi}$ is a bounded net which converges strictly to $(\io_A \ot \vfi)(x)$. Because $\rho$ is assumed to be strict, this implies that $\bigl(\,\rho((\io_A \ot \om)(x))\,\bigr)_{\om \in \cG_\vfi}$ converges strictly to $\rho((\io_A \ot \vfi)(x))$.

Hence $\bigl(\,(\io_B \ot \om)((\rho \ot \io)(x))\,\bigl)_{\om \in \cG_\vfi}$ converges strictly to $\rho((\io_A \ot \vfi)(x))$.
By definition, this implies that $(\rho \ot \io)(x)$ belongs to $\bar{\cM}_{\io_B \ot \vfi}^+$ and
$(\io_A \ot \vfi)\bigl((\rho \ot \io)(x)\bigr) = \rho((\io_A \ot \vfi)(x))$.

\item By the first result, we know that $(\rho \ot \io)(x^* x)$ belongs to
$\bar{\cM}_{\io_B \ot \vfi}^+$.
Because $(\rho \ot \io)(x)^* (\rho \ot \io)(x) \leq \|\rho\| \, (\rho \ot \io)(x^* x)$ (see lemma 5.3 of \cite{Lan}), this implies that $(\rho \ot \io)(x)^* (\rho \ot \io)(x)$ belongs to $\bar{\cM}_{\io_B \ot \vfi}^+$.
Therefore $(\rho \ot \io)(x)$ belongs to $\bar{\cN}_{\io_B \ot \vfi}$.

\medskip

Take $v \in H$. Because of result~\ref{cutoff} we have for all $\om \in \cG_\vfi$ that
\begin{eqnarray*}
& & \rho\bigl(\,(1 \ot \th_v^*)(1 \ot T_\om^{\frac{1}{2}})(\io_A \ot \lafi)(x)\,\bigr)
= \rho\bigl(\, (1 \ot \th_v^*)(\io_A \ot \pifi)(x)(1 \ot \th_{\xi_\om})\,\bigr) \\
& & \spat = \rho\bigl((\io_A \ot \om_{\xi_\om,v})(x)\bigr)
= (\io_B \ot \om_{\xi_\om,v})\bigl((\rho \ot \io)(x)\bigr) \\
& & \spat =   (1 \ot \th_v^*)(\io_B \ot \pifi)\bigl((\rho \ot \io)(x) \bigr)(1 \ot \th_{\xi_\om})
= (1 \ot \th_v^*)(1 \ot T_\om^{\frac{1}{2}})(\io_B \ot \lafi)\bigl((\rho \ot \io)(x) \bigr) \ .
\end{eqnarray*}
Since $(T_\om^{\frac{1}{2}})_{\om \in \cG_\vfi}$ converges strongly to 1, this implies that $\rho\bigl(\, (1 \ot \th_v^*)(\io_A \ot \lafi)(x)\, \bigr) $ \newline $= (1 \ot \th_v^*)(\io_B \ot \lafi)\bigl((\rho \ot \io)(x)\bigr)$.
\end{enumerate}
\end{demo}

\bigskip\medskip

\section{The tensor products of weights}

It this section we will say something about the tensor product of two proper weights and a partial GNS-construction for this tensor product.

Therefore we fix two \cst-algebras $A$ and $B$, a proper weight $\vfi$
on $A$ and a proper weight $\psi$ on $B$. Let us also fix a
GNS-construction $(H_\vfi,\pifi,\lafi)$ for $\vfi$ and a
GNS-construction $(H_\psi,\pips,\laps)$ for $\psi$.

\medskip

\begin{definition}
We define the tensor product weight $\vfi \ot \psi$ on $A \ot B$ in such a way that
$$(\vfi \ot \psi)(x) = \sup \, \{ \, (\om \ot \th)(x) \mid
\om \in \cF_\vfi , \th \in \cF_\psi \, \} $$
for every $x \in (A \ot B)^+$. Then $\vfi \ot \psi$ is a proper weight on $A \ot B$.
\end{definition}

It is then easy to see that the following properties hold:
\begin{itemize}
\item $\cM_\vfi \od \cM_\psi \subseteq \cM_{\vfi
\ot \psi}$ and $(\vfi \ot \psi)(a \ot b) =
\vfi(a) \, \psi(b)$ for $a \in \cM_\vfi$ and $b \in \cM_\psi$
\item $\cN_\vfi \od \cN_\psi \subseteq \cN_{\vfi
\ot \psi}$
\end{itemize}

\medskip

Notice that the family $\{\,\om \ot \th \mid \om \in \cF_\vfi , \th \in \cF_\psi\,\}$ is upwardly directed. This will imply that the map $\vfi \ot \psi$ above is additive.

\bigskip

The defining formula for $\vfi \ot \psi$ can be easily extended to the multiplier algebra of the tensor product.

\begin{proposition}
We have for all $x \in M(A \ot B)^+$ that
$$(\vfi \ot \psi)(x) = \sup \, \{ \, (\om \ot \th)(x) \mid
\om \in \cF_\vfi , \th \in \cF_\psi \, \}. $$
\end{proposition}
\begin{demo}
Define the function $\eta : M(A \ot B)^+ \rightarrow [0,\infty]$ such that
$$\eta(x) = \sup \, \{ \, (\om \ot \th)(x) \mid
\om \in \cF_\vfi , \th \in \cF_\psi \, \} \ .$$
Then $\eta$ is clearly a weight on $M(A \ot B)$ which is strictly lower semi-continuous and is equal to $\vfi \ot \psi$ on $(A \ot B)^+$.
Choose $x \in M(A \ot B)^+$. Take an approximate unit $(e_i)_{i \in I}$ for $A \ot B$. Then $(\,x^{\frac{1}{2}} e_i \, x^{\frac{1}{2}}\,)_{i \in I}$ is an increasing net in $(A \ot B)^+$ which converges strictly to $x$. Hence the strict lower semi-continuity of $\overline{\vfi \ot \psi}$ and $\eta$ implies that $(\,(\vfi \ot \psi)(x^{\frac{1}{2}} e_i \, x^{\frac{1}{2}})\,)_{i \in I}$ converges to $(\vfi \ot \psi)(x)$ and that
$(\,\eta(x^{\frac{1}{2}} e_i \, x^{\frac{1}{2}})\,)_{i \in I}$ converges to $\eta(x)$. Because $(\vfi \ot \psi)(x^{\frac{1}{2}} e_i \, x^{\frac{1}{2}}) = \eta(x^{\frac{1}{2}} e_i \, x^{\frac{1}{2}})$ for $i \in I$, this implies that $(\vfi \ot \psi)(x) = \eta(x)$.
\end{demo}

This will again imply the following results:
\begin{itemize}
\item $\bar{\cM}_\vfi \od \bar{\cM}_\psi \subseteq  \bar{\cM}_{\vfi
\ot \psi}$ and $(\vfi \ot \psi)(a \ot b) =
\vfi(a) \, \psi(b)$ for $a \in \bar{\cM}_\vfi$ and $b \in \bar{\cM}_\psi$
\item $\bar{\cN}_\vfi \od \bar{\cN}_\psi \subseteq \bar{\cN}_{\vfi
\ot \psi}$
\end{itemize}

\bigskip

Let $(H,\pi,\la)$ be a GNS-construction for $\vfi \ot \psi$.
Then it is easy to see that there exists an isometry $U : H_\vfi \ot H_\psi \rightarrow H$ such that $U (\lafi(a) \ot \laps(b))
= \la(a \ot b)$ for $a \in \bar{\cN}_\vfi$ and $b \in \bar{\cN}_\psi$.
A straightforward calculation shows that $\pi(a \ot b) \, U
= U \, (\pifi \ot \pips)(a \ot b)$ for $a \in A$ and $b \in B$. Therefore $\pi(x) \, U = U \, (\pifi \ot \pips)(x)$ for all $x \in M(A \ot B)$

\medskip

It is not clear to us whether $U$ is unitary in the general case. However, if one of the weights $\psi$ and $\vfi$ is a KMS weight, then it is. We will further comment on this later.

\bigskip

We introduce a special subset of $\bar{\cN}_{\vfi \ot \psi}$ consisting of all elements which can be properly represented in $H_\vfi \ot H_\psi$.

\medskip

\begin{definition} \label{leukedef}
We define the following objects:
\begin{itemize}
\item We define $\bar{\cN}(\vfi,\psi)$ as the set  of elements $x \in \bar{\cN}_{\vfi \ot \psi}$ such that there exists $v \in H_\vfi \ot H_\psi$ such that $\|v\|^2 = (\vfi \ot \psi)(x^* x)$ and $\langle v , \lafi(a) \ot \lafi(b) \rangle = (\vfi \ot \psi)((a^* \ot b^*)\,x)$ for all
$a \in \Nfi$ and $b  \in \Npsi$.
\item Also define $\cN(\vfi,\psi) = \bar{\cN}(\vfi,\psi) \cap (A \ot B)$.
\item We define the mapping $\lafi \ot \laps : \cN(\vfi,\psi) \rightarrow H_\vfi \ot H_\psi$ as follows. Let $x \in \cN(\vfi,\psi)$. Then we define $(\lafi \ot \laps)(x) \in H_\vfi \ot H_\psi$ such that
$$\langle (\lafi \ot \laps)(x) , \lafi(a) \ot \lafi(b) \rangle
= (\vfi \ot \psi)((a^* \ot b^*)\, x)$$ for $a \in \Nfi$ and $b  \in \Npsi$.
\end{itemize}
\end{definition}

It is clear that $(H_\vfi \ot H_\psi,\pifi \ot \pips,\lafi \ot \laps)$ is a GNS-construction for $\vfi \ot \psi$ if and only if $\cN(\vfi,\psi)= \cN_{\vfi \ot \psi}$. It is also easy to see that $\cN_\vfi \odot \cN_\psi \subseteq \cN(\vfi,\psi)$ and that $(\lafi \ot \laps)(x)  = (\lafi \od \laps)(x)$ for $x \in \cN_\vfi \odot \cN_\psi$.

\bigskip

We collect the GNS-like properties of $(H_\vfi \ot H_\psi, \pifi \ot \pips ,\lafi \ot \laps)$ in the next proposition.

\begin{proposition}
The following properties hold.
\begin{itemize}
\item The mapping $\lafi \ot \laps : \cN(\vfi,\psi) \rightarrow H_\vfi \ot H_\psi$ is a linear map which is closed with respect to the strict topology on $A \ot B$ and the norm topology on $H_\vfi \ot H_\psi$.
\item The mapping $\lafi \ot \laps : \cN(\vfi,\psi) \rightarrow H_\vfi \ot H_\psi$ is closable with respect to the strict topology on $M(A \ot B)$ and the norm topology on $H_\vfi \ot H_\psi$. Denote its closure by $\overline{\lafi \ot \laps}$. Then $D(\,\overline{\lafi \ot \laps}\,)
= \bar{\cN}(\vfi,\psi)$ and we put $(\lafi \ot \laps)(a) = (\,\overline{\lafi \ot \laps}\,)(a)$ for $a \in \bar{\cN}(\vfi,\psi)$.
\item $\langle (\lafi \ot \lapsi)(x) , (\lafi \ot \lapsi)(y) \rangle
= (\vfi \ot \psi)(y^* x)$ for $x,y \in \bar{\cN}(\vfi,\psi)$.
\item $\bar{\cN}(\vfi,\psi)$  and $\cN(\vfi,\psi)$ are left ideals in $M(A \ot B)$ and
$(\pifi \ot \pips)(x) (\lafi \ot \laps)(a) = (\lafi \ot \laps)(x \,a)$ for $x \in M(A \ot B)$ and $a \in \bar{\cN}(\vfi,\psi)$.
\end{itemize}
\end{proposition}
\begin{demo}
Define a mapping $\Ga : \bar{\cN}(\vfi,\psi) \rightarrow H_\vfi \ot H_\psi$ as follows.

Consider $x \in \bar{\cN}(\vfi,\psi)$. Then we define $\Ga(x) \in H_\vfi \ot H_\psi$ such that
$$\langle \Ga(x) , \lafi(a) \ot \laps(b) \rangle
= (\vfi \ot \psi)((a^* \ot b^*)\, x)$$ for $a \in \Nfi$ and $b  \in \Npsi$.
So $\lafi \ot \laps$ is the restriction of $\Ga$ to $\cN(\vfi,\psi)$.

\medskip

Fix a GNS-construction $(H,\pi,\la)$ for $\vfi \ot \psi$ and define the isometry $U : H_\vfi \ot H_\psi \rightarrow H$ such that
$U (\lafi(a) \ot \laps(b)) = \la(a \ot b)$ for $a \in \cN_\vfi$ and $b  \in \cN_\psi$. We have for all $x \in \bar{\cN}_{\vfi \ot \psi}$ and $a \in \Nfi$ and $b  \in \Npsi$ that
$$\langle U^* \la(x) , \lafi(a) \ot \laps(b) \rangle
= \langle \la(x) , U (\lafi(a) \ot \laps(b)) \rangle
= \langle \la(x) , \la(a \ot b) \rangle = (\vfi \ot \psi)((a^* \ot b^*) \, x) \ .$$
This implies for every $x \in \bar{\cN}_{\vfi \ot \psi}$ that
$x \in \bar{\cN}(\vfi,\psi)$ $\Leftrightarrow$ $\|U^* \la(x)\|^2 = (\vfi \ot \psi)(x^* x)$ $\Leftrightarrow$ $\|U^* \la(x)\| = \|\la(x)\|$ $\Leftrightarrow$ $\la(x) \in U(H_\vfi \ot H_\psi)$.
Moreover, $\Ga(x) = U^* \la(x)$ and $U \Ga(x) = \la(x)$ for all $x \in \bar{\cN}(\vfi,\psi)$.

\smallskip

This implies easily that $\bar{\cN}(\vfi,\psi)$ is a subspace of $M(A \ot B)$ and that the map $\Ga$ is linear. We also have for every $x,y \in \bar{\cN}(\vfi,\psi)$ that
$$\langle  \Ga(x) , \Ga(y) \rangle =
\langle U \Ga(x) , U \Ga(y) \rangle
= \langle \la(x) , \la(y) \rangle = (\vfi \ot \psi)(y^* x) \ .$$

\smallskip

Choose $x \in M(A \ot B)$ and $a \in \bar{\cN}(\vfi,\psi)$.
By the remarks before definition~\ref{leukedef} we get that
$U (\pifi \ot \pips)(x) = \pi(x) U$. Then
$x a \in \bar{\cN}_{\vfi \ot \psi}$ and
$$\la(x a) = \pi(x) \la(a) = \pi(x) U \Ga(a)
= U (\pifi \ot \pips)(x) \Ga(a) \in U(H_\vfi \ot H_\psi) \ .$$
Hence $x a \in \bar{\cN}(\vfi,\psi)$ and $\Ga(x a) = U^* \la(x a) =
(\pifi \ot \pips)(x) \Ga(a)$. \inlabel{ster3}

\smallskip

Choose a net $(x_i)_{i \in I}$ in $\bar{\cN}(\vfi,\psi)$, $x \in M(A \ot B)$ and $v \in H_\vfi \ot H_\psi$ such that
$(x_i)_{i \in I}$ converges strictly to $x$ and
$(\, \Ga (x_i)\,)_{i \in I}$ converges to $v$.
Because $U \Ga(x_i) = \la(x_i)$ for $i \in I$, this implies that $(\la(x_i))_{i \in I}$ converges to $U v$. So the strict closedness of
$\overline{\la}$ (see proposition \ref{weight.prop11}) implies that
$x \in \bar{\cN}_{\vfi \ot \psi}$ and $\la(x) = U v$.
Since $\la(x) \in U(H_\vfi \ot H_\psi)$, we get that $x \in \bar{\cN}(\vfi,\psi)$ and $\Ga(x)
= U^* \la(x) = v$. Hence $\Ga$ is strict-norm closed. Since $\lafi \ot \laps$ is the restriction of $\Ga$ to $\bar{\cN}(\vfi,\psi) \cap A$, it follows that $\lafi \ot \laps$ is also strict-norm closed.

\medskip

Using equation~\ref{ster3} and an approximate unit for $A$, we see that $\cN(\vfi,\psi)$ is a strict core for $\Ga$.
\end{demo}

\medskip

It is easy to see that $\bar{\cN}_\vfi \odot \bar{\cN}_\psi \subseteq \bar{\cN}(\vfi,\psi)$ and that $(\lafi \ot \laps)(x)  = (\lafi \od \laps)(x)$ for $x \in \bar{\cN}_\vfi \odot \bar{\cN}_\psi$.

\bigskip

In notation \ref{weight.not1}, we used the operators $T_\om$ to cut off $\lafi$ to a continuous map. In the next proposition, we see that $\lafi \ot \laps$ is also cut off properly by the tensor product of these operators.

\begin{proposition} \label{weight.prop12}
We have for all $\om \in \cF_\vfi$, $\th \in \cF_\psi$ and $x \in \bar{\cN}(\vfi,\psi)$ that
$(T_\om^{\frac{1}{2}} \ot T_\th^{\frac{1}{2}})(\lafi \ot \lapsi)(x)
= (\pifi \ot \pips)(x) (\xi_\om \ot \xi_\th)$.
\end{proposition}
\begin{demo}
We have for all $x \in \bar{\cN}(\vfi,\psi)$ that
$$\| (\pifi \ot \pips)(x) (\xi_\om \ot \xi_\th) \|^2
= (\om \ot \th)(x^* x) \leq (\vfi \ot \psi)(x^* x) = \| (\lafi \ot \lapsi)(x)\|^2$$
which implies the existence of a bounded operator $T \in B(H_\vfi \ot H_\psi)$ such that $T (\lafi \ot \lapsi)(x) =$\newline $(\pifi \ot \pips)(x) (\xi_\om \ot \xi_\th)$ for $x \in \bar{\cN}(\vfi,\psi)$.
We have for $a \in \Nfi$, $b \in \Nps$ that
$$(T_\om^{\frac{1}{2}} \ot T_\th^{\frac{1}{2}})(\lafi(a) \ot \laps(b))
= \pifi(a) \xi_\om \ot \pips(b) \xi_\th = (\pifi \ot \pips)(a \ot b)(\xi_\om \ot \xi_\th) = T (\lafi(a) \ot \laps(b)) \ .$$ Hence $T = T_\om^{\frac{1}{2}} \ot T_\th^{\frac{1}{2}}$.
\end{demo}

\bigskip

Let us quickly prove the next Fubini-type result.

\begin{result} \label{tensor.res1}
Let $x \in \bar{\cM}_{\io \ot \psi}^+$. Then $\vfi\bigl((\io \ot \psi)(x)\bigr)
= (\vfi \ot \psi)(x)$.
\end{result}
\begin{demo}
We have that
\begin{eqnarray*}
 (\vfi \ot \psi)(x) & =  & \sup_{\om \in \cF_\vfi,\th \in \cF_\psi} \ (\om \ot \th)(x)
 =  \sup_{\om \in \cF_\vfi} \,\bigl(\, \sup_{\th \in \cF_\psi} \, (\om \ot \th)(x) \, \bigr) \\
&  = &\sup_{\om \in \cF_\vfi} \,\bigl(\, \sup_{\th \in \cF_\psi} \, \om((\io \ot \th)(x)) \, \bigr)
=  \sup_{\om \in \cF_\vfi} \, \om\bigl((\io \ot \psi)(x)\bigr)
= \vfi\bigl((\io \ot \psi)(x)\bigr) \ .
\end{eqnarray*}
\end{demo}

\medskip

\begin{remark} \rm
Consider $a \in \bar{\cN}_\psi$ and $\om \in A^*$. Then we have for $\th \in \cF_\psi$ that
$$\om\bigl((\io \ot a\,\th\, a^*)(x)\bigr)
=  \th(a^*\,(\om \ot \io)(x) \, a) = \langle \pips((\om \ot \io)(x)) \, T_\th^{\frac{1}{2}} \, \laps(a) , T_\th^{\frac{1}{2}} \, \laps(a) \rangle \, $$
So we see that the net $\bigl(\,\om\bigl((\io \ot a\,\th\, a^*)(x)\bigr)\,\bigr)_{\th \in \cG_\psi}$ converges to $\psi(a^* \, (\io \ot \om)(x) \, a)$, which is equal to \newline $\om\bigl((\io \ot a\,\psi\,a^*)(x)\bigr)$.

\medskip

Now lemma \ref{weight.lem5} tells us that the net $\bigl(\,(\io \ot a\,\th\, a^*)(x)\,)_{\th \in \cG_\psi}$ converges strictly to $(\io \ot a\,\psi\,a^*)(x)$. If $(\io \ot a\,\psi\,a^*)(x)$ happens to belong to $A$, then this convergence is in norm.
\end{remark}

\medskip

The previous remark implies that result \ref{tensor.res1} implies the next one.

\begin{result} \label{tensor.res2}
Let $x \in M(A \ot B)^+$ and $a \in \bar{\cN}_{\psi}$. Then  $\vfi\bigl((\io \ot a \psi a^*)(x)\bigr)
= (\vfi \ot \psi)((1 \ot a^*)\,x \, (1 \ot a))$.
\end{result}

\bigskip

The next result  gives some interesting elements belonging to $\bar{\cN}_{\vfi \ot \psi}$. We will make use of result \ref{weight.res2}.

\begin{proposition} \label{tensor.prop1}
Consider $x \in \bar{\cN}_{\vfi \ot \psi} \cap \bar{\cN}_{\io \ot \psi}$. Then $x$ belongs to $\bar{\cN}(\vfi,\psi)$ and the following holds:
\begin{itemize}
\item Let $v \in H_\psi$ and define the element $q \in M(A)$ such that
$\mu(q) = \langle \lapsi((\mu \ot \io)(x)) , v \rangle$ for $\mu \in A^*$.
Then $q$ belongs to $\bar{\cN}_\vfi$ and $\vfi(q^* q) \leq \|v\|^2 \, (\vfi \ot \psi)(x^* x)$.
\item Let $(e_i)_{i \in I}$ be an orthonormal basis of $H_\psi$. Define for every $i \in I$ the element $q_i \in \bar{\cN}_\vfi$ such that
$\mu(q_i) = \langle \lapsi((\mu \ot \io)(x)) , e_i \rangle$ for $\mu \in A^*$.
Then $\sum_{i \in I} \|\lafi(q_i)\|^2 = (\vfi \ot \psi)(x^* x) < \infty$ and
$$(\lafi \ot \laps)(x) = \sum_{i \in I} \, \lafi(q_i) \ot e_i \ .$$
\end{itemize}
\end{proposition}
\begin{demo}
It is enough to prove the second result. But we saw in the proof of proposition \ref{weight.prop2} that for every $\mu \in A^*_+$, we have that
$\sum_{i \in I} \mu(q_i^* q_i) = \mu((\io \ot \psi)(x^* x))$.
Using definition \ref{weight1.def3} this implies that $$\sum_{i \in I} \vfi(q_i^* q_i)
= \vfi((\io \ot \psi)(x^* x)) = (\vfi \ot \psi)(x^* x) < \infty \ .$$
So we find in particular that $q_i \in \bar{\cN}_\vfi$ for all $i \in I$.

\medskip

Define $v = \sum_{i \in I} \, \lafi(q_i) \ot e_i$. Then the above equality implies that $\|v\|^2 = (\vfi \ot \psi)(x^* x)$.

\medskip

Take $a \in \Nfi$ and $b \in \Npsi$. Then we have for all $\om \in \cF_\vfi$ that
\begin{eqnarray*}
& & \langle (T_\om \ot 1) v , \lafi(a) \ot \lapsi(b) \rangle
=  \sum_{i \in I} \langle  T_\om^{\frac{1}{2}} \lafi(q_i) \ot e_i , T_\om^{\frac{1}{2}} \lafi(a) \ot \lapsi(b) \rangle \\
& & \spat = \sum_{i \in I} \langle \pifi(q_i) \xi_\om ,T_\om^{\frac{1}{2}} \lafi(a) \rangle \, \langle e_i , \lapsi(b) \rangle
= \sum_{i \in I} \langle \lapsi ((\om_{\xi_\om,T_\om^{\frac{1}{2}} \lafi(a)} \ot \io)(x)) , e_i \rangle \, \langle e_i , \lapsi(b) \rangle \\
& & \spat = \langle \lapsi\bigl((\om_{\xi_\om,T_\om^{\frac{1}{2}}\lafi(a)} \ot \io)(x)\bigr) ,\lapsi(b) \rangle
=  \psi\bigl(b^* (\om_{\xi_\om,T_\om^{\frac{1}{2}} \lafi(a)} \ot \io)(x)\bigr) \\
& & \spat = \om_{\xi_\om, T_\om^{\frac{1}{2}} \lafi(a)}\bigl((\io \ot \psi)((1 \ot b^*) \, x)\bigr)
=  \langle \pifi\bigl((\io \ot \psi)((1 \ot b^*)\, x)\bigr)\, \xi_\om , T_\om^{\frac{1}{2}} \lafi(a) \rangle \\
& & \spat = \langle \pifi\bigl((\io \ot \psi)((1 \ot b^*)\, x)\bigr) \,\xi_\om , \pifi(a) \xi_\om \rangle
=  \om(a^* (\io \ot \psi)((1 \ot b^*)\, x)) \\
& & \spat = \om\bigl((\io \ot \psi)((a^* \ot b^*)\, x)\bigr) \ .
\end{eqnarray*}
Since $(T_\om^{\frac{1}{2}})_{\om \in \cG_\vfi}$ converges to $1$,  we see that
$\langle v , \lafi(a) \ot \lapsi(b) \rangle = (\vfi \ot \psi)((a^* \ot  b^*) \,x)$. This implies by definition that $x \in \bar{\cN}(\vfi,\psi)$ and
$(\lafi \ot \laps)(x) = v$.
\end{demo}

\medskip

We have of course a similar result  if $x$ belongs to $\bar{\cN}_{\vfi \ot \psi} \cap \bar{\cN}_{\vfi \ot \io}$.

\medskip\medskip

\begin{corollary}
Consider $x \in \bar{\cN}_{\vfi \ot \io}$, $y \in \bar{\cN}_\psi$ and an orthonormal basis $(e_i)_{i \in I}$ for $H_\psi$. Then $x\,(1 \ot y)$ belongs to $\bar{\cN}_{\vfi \ot \psi} \, \cap \,
\bar{\cN}_{\vfi \ot \io} \, \cap \, \bar{\cN}_{\io \ot \psi}$,
$\sum_{i \in I} \, \|\lafi\bigl((\io \ot \om_{\laps(y),e_i})(x)\bigr)\|^2
= (\vfi \ot \psi)((1 \ot y^*)\,x^* x\,(1 \ot y)) < \infty$ and
$$(\lafi \ot \laps)(x\,(1 \ot y))
= \sum_{i \in I} \lafi\bigl((\io \ot \om_{\laps(y),e_i})(x)\bigr) \ot e_i \ .$$
\end{corollary}
\begin{demo}
It is clear that $x\,(1 \ot y)$ belongs to $\bar{\cN}_{\io \ot \psi} \, \cap \, \bar{\cN}_{\vfi \ot \io}$.

By result \ref{tensor.res2}, we have  that
$$(\vfi \ot \psi)\bigl((1 \ot y^*)\,x^* x\,(1 \ot y)\bigr)
= \vfi((\io  \ot y \psi y^*)(x^* x) )  =
 \psi(y^* \, (\vfi \ot \io)(x^* x) \, y) \ .$$
So we get that $x\,(1 \ot y)$ belongs to $\bar{\cN}_{\vfi \ot \psi}$ and we can apply the previous proposition. Thus we have for every $i \in I$ an element $q_i \in \bar{\cN}_\psi$ such that $\mu(q_i) = \langle \laps\bigl((\mu \ot \io)(x\,(1 \ot y))\bigr) , e_i \rangle$ for all $\mu \in A^*$. We also know that $\sum_{i \in I} \|\lafi(q_i)\|^2 = (\vfi \ot \psi)((1 \ot y^*)\,x^* x\, (1 \ot y))$ and
$$(\lafi \ot \laps)(x\,(1 \ot y)) = \sum_{i \in I} \lafi(q_i) \ot e_i \ .$$
But we have for all $\mu \in A^*$ and $i \in I$ that
\begin{eqnarray*}
\mu(q_i) & = & \langle \laps\bigl((\mu \ot \io)(x\,(1 \ot y))\bigr) , e_i \rangle = \langle \laps((\mu \ot \io)(x)\,y) , e_i \rangle\\
& = & \langle \pips((\mu \ot \io)(x))\, \laps(y) , e_i \rangle
= \mu\bigl((\io \ot \om_{\laps(y),e_i})(x)\bigr) \ .
\end{eqnarray*}
So we see that $q_i = (\io \ot \om_{\laps(y),e_i})(x)$ for all $i \in I$.
\end{demo}

\bigskip

Not surprisingly, relative invariance is preserved under the tensor product construction.

\begin{proposition} \label{weight.prop16}
Let $\al : A \rightarrow A$ and
$\be : B \rightarrow B$ be $^*$-isomorphisms such that there exist numbers $\lambda, \nu > 0$ such that $\vfi \, \al = \lambda \, \vfi$ and $\psi \be = \nu \, \psi$. Then $(\vfi \ot \psi)(\al \ot \be) = \lambda \nu \, (\vfi \ot \psi)$.
\end{proposition}
\begin{demo}
It is easy to see that $\cF_\vfi = \{ \, \lambda^{-1} \, \om \al \mid \om \in \cF_\vfi \, \}$ and
$\cF_\psi = \{ \, \nu^{-1} \, \th \be \mid \th \in \cF_\psi \, \}$. Hence,
\begin{eqnarray*}
(\vfi \ot \psi)(x)  & = & \sup \, \{ \, (\om \ot \th)(x) \mid
\om \in \cF_\vfi , \th \in \cF_\psi \, \}  \\
& = & \sup \, \{ \, (\lambda \nu)^{-1} \, (\om \al \ot \th \be)(x) \mid \om \in \cF_\vfi , \th \in \cF_\psi \, \} \\
& = & \sup \, \{ \, (\lambda \nu)^{-1} \, (\om \ot \th)\bigl((\al \ot \be)(x)\bigr) \mid
\om \in \cF_\vfi , \th \in \cF_\psi \, \}  \\
& = &(\lambda \nu)^{-1} \, (\vfi \ot \psi)\bigl((\al \ot \be)(x)\bigr)\ .
\end{eqnarray*}
\end{demo}

\medskip

\begin{proposition} \label{weight.prop17}
Let $\al : A \rightarrow A$ and $\be : B \rightarrow B$ be $^*$-isomorphisms such that there exist numbers $\lambda, \nu > 0$ such that $\vfi\, \al = \lambda \, \vfi$ and $\psi\, \be = \nu \, \psi$.

Define unitary operators $U \in B(H_\vfi)$ and $V \in B(H_\psi)$ such that
$U \lafi(a)= \lambda^{-\frac{1}{2}} \, \lafi(\al(a))$ for $a \in \Nfi$ and
$V \laps(b)= \nu^{-\frac{1}{2}} \, \laps(\be(a))$ for $b \in \Npsi$.
Then we have for all $x \in \bar{\cN}(\vfi,\psi)$ that
$(\al \ot \be)(x)$ belongs to $\bar{\cN}(\vfi,\psi)$ and
$(\lafi \ot \laps)\bigl((\al \ot \be)(x)\bigr) = (\lambda \nu)^{\frac{1}{2}} \, (U \ot V)(\lafi \ot \laps)(x)$.
\end{proposition}
\begin{demo}
By the previous result, we know that $(\al \ot \be)(x) \in \bar{\cN}_{\vfi \ot \psi}$ and
\begin{eqnarray*}
(\vfi \ot \psi)\bigl(\,(\al \ot \be)(x)^* (\al \ot \be)(x)\,\bigr) & = & \lambda \nu \, (\vfi \ot \psi)(x^* x) = \lambda \nu \, \|(\lafi \ot \laps)(x)\|^2 \\
& = & \lambda \nu \, \|(U \ot V)(\lafi \ot \laps)(x)\|^2 \ .
\end{eqnarray*}
We have moreover for all $a \in \Nfi$, $b \in \Npsi$ that
\begin{eqnarray*}
& &  \lambda^{\frac{1}{2}} \nu^{\frac{1}{2}}  \, \langle (U \ot V)(\lafi \ot \laps)(x) , \lafi(a) \ot \laps(b) \rangle
= \lambda^{\frac{1}{2}} \nu^{\frac{1}{2}}  \, \langle (\lafi \ot \laps)(x) , U^* \lafi(a) \ot V^* \laps(b) \rangle \\
& & \spat = \lambda \nu \, \langle (\lafi \ot \laps)(x) , \lafi(\al^{-1}(a) )  \ot \laps(\be^{-1}(b)) \rangle
= \lambda \nu \, (\vfi \ot \psi)\bigl(\,(\al^{-1}(a)^* \ot \be^{-1}(b)^*)\, x \,\bigr) \ .
\end{eqnarray*}
Using the previous proposition again, this implies that
\begin{eqnarray*}
\lambda^{\frac{1}{2}} \nu^{\frac{1}{2}}  \, \langle (U \ot V)(\lafi \ot \laps)(x) , \lafi(a) \ot \laps(b) \rangle
& = & (\vfi \ot \psi)\bigl(\,(\al \ot \be)((\al^{-1}(a)^* \ot \be^{-1}(b)^*)\, x)\,\bigr) \\
& = & (\vfi \ot \psi)((a^* \ot b^*)\,(\al \ot \be)(x))\ .
\end{eqnarray*}
By definition, this implies that $(\al \ot \be)(x) \in \bar{\cN}(\vfi,\psi)$ and $(\lafi \ot \laps)((\al \ot \be)(x))
=$ \newline $ \lambda^{\frac{1}{2}} \nu^{\frac{1}{2}} \, (U \ot V)(\lafi \ot \laps)(x)$.
\end{demo}

\bigskip\medskip

In a last part of this section, we will look a little bit closer to the case of KMS weights. In the next proposition, we combine a technique of Jan Verding (see \cite{Verd}) with proposition \ref{tensor.prop1}. We will need  lemma \ref{app.lem1}  which is due to Jan Verding (see lemma A.1.2 of \cite{Verd}).

\begin{proposition} \label{tensor.prop2}
Suppose that $\psi$ is a KMS weight with modular group $\si$ and let $J$ denote the modular conjugation of $\psi$ in the GNS-construction $(H_\psi,\pips,\laps)$. Then $\bar{\cN}(\vfi,\psi) = \bar{\cN}_{\vfi \ot \psi}$ and we have for all $x \in \bar{\cN}(\vfi,\psi)$ and $a \in D(\overline{\si}_{\frac{i}{2}})$ that
$x\,(1 \ot a)$ belongs to $\bar{\cN}(\vfi,\psi)$ and
$(\lafi \ot \laps)(x\,(1 \ot a)) = (1 \ot J \pips(\si_{\frac{i}{2}}(a))^* J )(\lafi \ot \laps)(x)$.
\end{proposition}
\begin{demo}
Let $(H,\pi,\la)$ be a GNS-construction for $\vfi \ot \psi$.
We know there exists a net $(e_j)_{j \in J}$ in $\Npsi \cap D(\si_{\frac{i}{2}})$ such that
\begin{itemize}
\item $(e_j)_{j \in J}$ is bounded and converges strictly to 1
\item $(\si_{\frac{i}{2}}(e_j))_{j \in J}$ is bounded and converges strongly to 1
\end{itemize}

\smallskip

Consider $b \in D(\overline{\si}_{\frac{i}{2}})$. Choose $y \in \bar{\cN}_{\vfi \ot \psi}$. Take $\om \in \cF_\vfi$. By result \ref{tensor.res1}, we have that
$$\psi((\om \ot \io)(y^* y)) = (\om \ot \psi)(y^* y) \leq (\vfi \ot \psi)(y^* y) < \infty \ .$$
So we get for all $\om \in \cF_\vfi$ and $\th \in \cF_\psi$ that
\begin{eqnarray*}
& & (\om \ot \th)((1 \ot b^*)\,y^* y\,(1 \ot b))
= \th(b^* (\om \ot \io)(y^* y)  b) \\
& & \spat \leq \psi(b^* (\om \ot \io)(y^* y) \, b)
= \| \laps((\om \ot \io)(y^* y)^{\frac{1}{2}} \, b) \|^2 \\
& & \spat = \| J \pips(\si_{\frac{i}{2}}(b))^* J \laps((\om \ot \io)(y^* y)^{\frac{1}{2}}) \|^2 \leq \|\si_{\frac{i}{2}}(b)\|^2 \, \|\laps((\om \ot \io)(y^* y)^{\frac{1}{2}}) \|^2 \\
& & \spat = \|\si_{\frac{i}{2}}(b)\|^2 \, \psi((\om \ot \io)(y^* y))
\leq \|\si_{\frac{i}{2}}(b)\|^2 \, (\vfi \ot \psi)(y^* y) \ .
\end{eqnarray*}
This implies that $y\,(1 \ot b)$ belongs to $\bar{\cN}_{\vfi \ot \psi}$ and
\begin{equation} \label{vgl1}
(\vfi \ot \psi)((1 \ot b^*)\,y^* y\,(1 \ot b)) \leq \|\si_{\frac{i}{2}}(b)\|^2 \, (\vfi \ot \psi)(y^* y) \ .
\end{equation}

Choose $x \in \cN_{\vfi \ot \psi}$.
We have for $j \in J$ that $x(1 \ot e_j) \in \bar{\cN}_{\io \ot \psi} \cap \cN_{\vfi \ot \psi}$ which by proposition \ref{tensor.prop1} implies that $x \,(1 \ot e_j)$ belongs to $\cN(\vfi,\psi)$. Moreover,
$$\|(\lafi \ot \lapsi)(x\,(1 \ot e_j))\|^2 = (\vfi \ot \psi)((1 \ot e_j^*)\, x^* x \, (1 \ot e_j)) \leq \|\si_{\frac{i}{2}}(e_j)\|^2 \, (\vfi \ot \psi)(x^* x) \ ,$$
so the net $\bigl(\, (\lafi \ot \laps)(x\,(1 \ot e_j))\,\bigr)_{j \in J}$ is bounded. Because $(x\,(1 \ot e_j))_{j \in J}$ converges to $x$,  lemma \ref{app.lem1}  implies the existence of a sequence $(x_n)_{n=1}^\infty$ in the convex hull of $\{ \, x \,(1 \ot e_j) \mid j \in J\,\}$, which is contained in $\bar{\cN}(\vfi,\psi)$, such that
$(x_n)_{n=1}^\infty$ converges to $x$ and
$(\,(\lafi \ot \laps)(x_n)\,)_{j \in J}$ is convergent. The closedness of $\lafi \ot \laps$ therefore gives us that $x$ belongs to $\cN(\vfi,\psi)$.

\medskip

Take a GNS-construction $(H,\pi,\la)$ for $\vfi \ot \psi$ and define the isometry $U : H_\vfi \ot H_\psi \rightarrow H$ such that
$U (\lafi \odot \laps)(p) = \la(p)$ for $p \in \bar{\cN}_\vfi \odot \bar{\cN}_\psi$.

\medskip

Now take $z \in \bar{\cN}_{\vfi \ot \psi}$. Fix an approximate unit $(u_k)_{k \in K}$ for $A \ot B$. By the above result, we have for every $k \in K$  that $u_k \, z$ belongs to $\bar{\cN}(\vfi,\psi)$ (since $u_k \, z$ belongs to $\cN_{\vfi \ot \psi}$).
We have moreover for every $k \in K$ that
$$(\lafi \ot \lapsi)(u_k \, z) = U^* \la(u_k \, z) = U^* \pi(u_k) \la(z) \ ,$$
which implies that $(\,(\lafi \ot \lapsi)(u_k \, z)\,)_{k \in K}$ converges to $U^* \la(z)$. Because $(u_k \, z)_{k \in K}$ converges strictly to $z$, the strict closedness of $\overline{\lafi \ot \lapsi}$ gives us that $z$ belongs to $\bar{\cN}(\vfi,\psi)$.

\smallskip

So we have proven that $\bar{\cN}(\vfi,\psi) = \bar{\cN}_{\vfi \ot \psi}$.

\medskip

Choose $a \in D(\overline{\si}_{\frac{i}{2}})$. By inequality~\ref{vgl1}, we know already that $y \, (1 \ot a)$ belongs to $\bar{\cN}(\vfi,\psi)$ for all $y \in \bar{\cN}_{\vfi \ot \psi}$. Inequality~\ref{vgl1} also allows to define a bounded operator $T \in B(H_\vfi \ot H_\psi)$ such that \newline $T (\lafi  \ot \laps)(y) = (\lafi \ot \laps)(y\,(1 \ot a))$ for $y \in
\bar{\cN}(\vfi,\psi)$. We have for all $p \in \Nfi$, $q \in \Npsi$ that
$$T (\lafi(p) \ot \laps(q)) = (\lafi \ot \laps)((p \ot q)(1 \ot a))
=\lafi(p) \ot \laps(q a) = \lafi(p) \ot J \pips(\si_{\frac{i}{2}}(a))^* J \laps(q) \ .$$
So $T$ must be equal to $1 \ot J \pips(\si_{\frac{i}{2}}(a))^* J $.
\end{demo}

\medskip

It goes without saying that a similar result holds if $\vfi$ is a KMS weight.

\bigskip

\begin{proposition}
Suppose that $\vfi$ is a KMS weight with modular group $\si$ and $\psi$ a KMS weight on $B$ with modular group $\si'$ and use the following notations.
\begin{itemize}
\item Let $J$ be the modular conjugation of $\vfi$ in the GNS-construction $(H_\vfi,\pifi,\lafi)$  and \newline
$J'$ the modular conjugation of $\psi$ in the GNS-construction $(H_\psi,\pips,\laps)$.
\item Let $\nab$ be the modular operator of $\vfi$ in the GNS-construction $(H_\vfi,\pifi,\lafi)$  and \newline
$\nabp$ the modular operator of $\psi$ in the GNS-construction $(H_\psi,\pips,\laps)$.
\end{itemize}
Then $\vfi \ot \psi$ is a KMS weight on $A \ot B$ with GNS-construction $(H_\vfi \ot H_\psi,\pifi \ot \pips,\lafi \ot \laps)$. We have moreover that
\begin{itemize}
\item $\Nfi \odot \Nps$ is a core for $\lafi \ot \laps$.
\item $\si \ot \si'$ is a modular group for $\vfi \ot \psi$.
\item $J \ot J'$ is the modular conjugation of $\vfi \ot \psi$ in the GNS-construction $(H_\vfi \ot H_\psi,\pifi \ot \pips,\la)$.
\item $\nab \ot \nabp$ is the modular operator of $\vfi \ot \psi$ in the GNS-construction $(H_\vfi \ot H_\psi,\pifi \ot \pips,\la)$.
\end{itemize}
\end{proposition}

\medskip

In section 7 of \cite{JK1}, we proved that $\vfi \ot \psi$ is a KMS weight on $A \ot B$ with modular group $\si \ot \si'$. We also proved the following fact. The mapping $\lafi \odot \lafi$ is closable and the closure $\la$ of $\lafi \odot \laps$ satisfies the following properties :
\begin{itemize}
\item $(H_\vfi \ot H_\psi,\pifi \ot \pips,\la)$ is a GNS-construction for $\vfi \ot \psi$.
\item $J \ot J'$ is the modular conjugation of $\vfi \ot \psi$ in the GNS-construction $(H_\vfi \ot H_\psi,\pifi \ot \pips,\la)$.
\item $\nab \ot \nabp$ is the modular operator of $\vfi \ot \psi$ in the GNS-construction $(H_\vfi \ot H_\psi,\pifi \ot \pips,\la)$.
\end{itemize}
Because the linear mapping $\lafi \ot \laps : \cN(\vfi,\psi) \rightarrow H_\vfi \ot H_\psi  $
is closed and $(\lafi \ot \laps)(x) = (\lafi \odot \laps)(x)$ for $x \in \Nfi \odot \Nps$, we have that $\cN_{\vfi \ot \psi} \subseteq \cN(\vfi,\psi)$ and $\la(x) = (\lafi \ot \laps)(x)$ for $x \in \cN_{\vfi \ot \psi}$.  But we have always that $\cN(\vfi,\psi) \subseteq \cN_{\vfi \ot \psi}$, hence $\cN(\vfi,\psi) = \cN_{\vfi \ot \psi}$.

\bigskip\bigskip

The next proposition is easy to check by using continuity arguments. Notice that in the notations of the next proposition, $(H_\om,\pi_\om,\xi_\om)$ and $(H_\th,\pi_\th,\xi_\th)$ are cyclic GNS-construction for $\om$ and $\th$ respectively and that
$\la_\om(a) = \pi_\om(a) \xi_\om$ for $a \in M(A)$ and
$\la_\th(b) = \pi_\th(b) \xi_\th$ for $b \in M(B)$.

\smallskip

\begin{proposition} \label{weight.prop13}
Consider two \cst-algebras $A$ and $B$. Let $\om$ be a positive functional  on $A$  with GNS-construction $(H_\om,\pi_\om,\la_\om)$ and $\th$ a positive functional on $B$ with GNS-construction $(H_\th,\pi_\th,\pi_\th)$. Define
$\xi_\om = \la_\om(1) \in H_\om$ and $\xi_\th = \la_\th(1) \in H_\th$.
Then $\bar{\cN}(\om,\th) = M(A \ot B)$, $\cN(\om,\th) = A \ot B$ and
$(\la_\om \ot \la_\th)(x) = (\pi_\om \ot \pi_\th)(x)(\xi_\om \ot \xi_\th)$ for $x \in M(A \ot B)$.
\end{proposition}

\bigskip

Let $A$ be a \cst-algebra, $\vfi$ a proper weight on $A$ and $\om \in \cF_\vfi$. Then we have a natural GNS-construction $(H_\om,\pi_\om,H_\om)$ for $\om$ such that
\begin{itemize}
\item $H_\om = [ \, T_\om^{\frac{1}{2}} v \mid v \in H \, ] = [\, \pi(a) \xi_\om \mid a \in A  \, ]$

\item $\pi_\om(x) v = \pi(x) v$ for $v \in H_\om$.
\item $\la_\om(a) = \pi(a) \xi_\om$ for $a \in A$
\end{itemize}
Notice that $\la_\om(a) = T_\om^{\frac{1}{2}} \la(a)$ for $a \in \Nfi$.

\bigskip

\begin{proposition}
Consider two \cst-algebras $A$ and $B$. Let $\vfi$ be a proper weight on $A$  with GNS-construction $(H_\vfi,\pifi,\lafi)$ and $\psi$ a proper weight on $B$ with GNS-construction $(H_\psi,\pips,\laps)$.
\begin{itemize}
\item Let $\om \in \cF_\vfi$ and $x \in \bar{\cN}(\vfi,\psi)$. Then
$x \in \bar{\cN}(\om,\psi)$ and $(T_\om^{\frac{1}{2}} \ot 1)(\lafi \ot \lapsi)(x) = (\la_\om \ot \lapsi)(x)$.
\item Let $\th \in \cF_\psi$ and $x \in \bar{\cN}(\vfi,\psi)$. Then
$x \in \bar{\cN}(\vfi,\th)$ and $(1 \ot T_\th^{\frac{1}{2}})(\lafi \ot \lapsi)(x) = (\la_\vfi \ot \la_\th)(x)$.
\item Let $\om \in \cF_\vfi$, $\th \in \cF_\psi$ and $x \in \bar{\cN}(\vfi,\psi)$. Then $(T_\om^{\frac{1}{2}} \ot T_\th^{\frac{1}{2}})(\lafi \ot \lapsi)(x) = (\la_\om \ot \la_\th)(x)$.
\end{itemize}
\end{proposition}
\begin{demo}
\begin{itemize}
\item Since $\om \ot \psi \leq \vfi \ot \psi$, we have certainly that
$x \in \bar{\cN}_{\om \ot \vfi}$. Because $H_\om = \overline{T_\om^{\frac{1}{2}} H_\vfi}$, it is clear that $(T_\om^{\frac{1}{2}} \ot 1)(\lafi \ot \lapsi)(x)$ belongs to $H_\om \ot H_\psi$.
By proposition \ref{weight.prop12}, we have for all $\th \in \cF_\psi$ that
$$\|(1 \ot T_\th^{\frac{1}{2}})\, (T_\om^{\frac{1}{2}} \ot 1)(\lafi \ot \lapsi)(x)\|^2 = \|(\pifi \ot \pips)(x)(\xi_\om \ot \xi_\th)\|^2
= (\om \ot \th)(x^* x) \ .$$
As $(\,T_\th^{\frac{1}{2}}\,)_{\th \in \cG_\psi}$ converges to $1$ and
$(\, (\om \ot \th)(x^* x)\,)_{\th \in \cG_\psi}$ converges to $(\om \ot \psi)(x^* x)$, we get that $\|(T_\om^{\frac{1}{2}} \ot 1)(\lafi \ot \lapsi)(x)\|^2 = (\om \ot \psi)(x^* x) \ .$

Let $a \in A$ and $b \in \cN_\psi$. Using proposition \ref{weight.prop12} again, we see that
\begin{eqnarray*}
& & \langle (1 \ot T_\th^{\frac{1}{2}})\, (T_\om^{\frac{1}{2}} \ot 1)(\lafi \ot \lapsi)(x) , (1 \ot T_\th^{\frac{1}{2}}) (\la_\om(a) \ot \laps(b)) \rangle \\
& & \spat = \langle (\pifi \ot \pips)(x) (\xi_\om \ot  \xi_\th) , \pifi(a) \xi_\om \ot \pips(b) \xi_\th \rangle
 = (\om \ot \th)((a^* \ot b^*)\, x) \ .
\end{eqnarray*}
As above, this implies that $\langle (T_\om^{\frac{1}{2}} \ot 1)(\lafi \ot \lapsi)(x) , \la_\om(a) \ot \laps(b) \rangle
= (\om \ot \psi)((a^* \ot b^*)\,x)$.
\item Analogous as above.
\item Follows from propositions \ref{weight.prop12} and \ref{weight.prop13}.
\end{itemize}
\end{demo}

\bigskip\medskip

\section{Integrals and closed linear mappings}

If we have a bounded linear map between two Banach spaces, the integral commutes with this bounded linear map. A similar result holds for closed linear mappings and the proof turns out to be very simple (see e.g. lemma 1.1 of \cite{JK1}).
By $\mu$-integrability of a Banach-valued function
we mean strong integrability as defined in e.g. {\S}3, Chapter X of
\cite{Lang}.

\begin{result}
Consider Banach spaces $E$ and $F$ and a closed linear map $\pi$ from
$D(\pi) \subseteq E$ into $F$. Let $X$ be a locally compact space with a regular Borel
measure $\mu$ and consider a function $f : X \rightarrow D(\pi)$ such
that
\begin{enumerate}
\item $f$ is $\mu$-integrable.
\item $\pi \circ f$ is $\mu$-integrable.
\end{enumerate}
Then $\int f(t) \, d\mu(t)$ belongs to $D(\pi)$ and $\pi\bigl(\,\int f(t) \, d\mu(t)\,\bigr) = \int \pi(f(t)) \, d\mu(t)$.
\end{result}

\bigskip\bigskip

If we work with integrals of strictly continuous functions, we will always work within the following setting.

\medskip

Consider a \cst-algebra $A$ and a locally compact space $X$ with regular Borel measure $\mu$. Let $f : X \rightarrow M(A)$ be a strictly continuous function such that there exists a $\mu$-integrable function $g : X \rightarrow \R^+$ such that $\|f(t)\| \leq g(t)$ for $t \in X$.

Then we have for every $b \in B$ that the functions
$X \rightarrow A : t \mapsto f(t) \, b$ and $X \rightarrow A : t \mapsto  b \, f(t)$ are norm continuous and dominated by $\|b\| \, g$.

It is then not difficult to see that there exists a unique element $\int f(t) \ d\mu(t) \in M(A)$ such that
$$\left(\, \int f(t) \, d\mu(t)\,\right) \, b = \int f(t) \, b \, d\mu(t)
\hspace{1.5cm} \text{ and } \hspace{1.5cm}  b \, \left(\, \int f(t) \, d\mu(t)\,\right)  = \int b\,f(t)  \, d\mu(t)$$
for all $b \in B$.

\medskip

Because any continuous functional on $A$ is of the form $\th \, b$ for some $\th \in A^*$ and $b \in A$, the following holds.
Consider $\om \in A^*$, then the function $X \rightarrow \C : t \mapsto
\om(f(t))$ is continuous and dominated by $\|\om\| \, g$. We have moreover that
\begin{equation} \label{vgl5}
\om \left(\, \int f(t) \, d\mu(t)\,\right) = \int \om(f(t))  d\mu(t) \ .
\end{equation}

\bigskip

Now consider two \cst-algebras $A$ and $B$ and a non-degenerate $^*$-homomorphism $\pi$ from $A$ into $M(B)$. Let $X$ be a locally compact space $X$ with regular Borel measure $\mu$ and  $f : X \rightarrow M(A)$  a bounded strictly continuous function such that there exists a $\mu$-integrable function $g : X \rightarrow \R^+$ such that $\|f(t)\| \leq g(t)$
for $t \in \R$

Then the function $X \rightarrow M(B) : t \mapsto \pi(f(t))$ is strictly continuous and dominated by $g$. We have (by using elements in $A^*$) moreover that
$$ \pi\left(\, \int f(t) \, d\mu(t)\,\right) = \int \pi(f(t)) \, d\mu(t) \ .$$

\bigskip

Using Hahn-Banach, this can be easily extended to strictly closed mappings.

\begin{result}
Consider \cst-algebras $A$ and $B$ and a strictly closed linear map $\pi$ from $D(\pi) \subseteq M(A)$ into $M(B)$. Let $X$ be a locally compact space with Borel measure $\mu$ and consider a function  $f : X \rightarrow D(\pi)$  such that there exists a $\mu$-integrable function $g : X \rightarrow \R^+$ such that
\begin{enumerate}
\item $f$ is strictly continuous and $\|f(t)\| \leq g(t)$ for $t \in X$.
\item $\pi \circ f$ is  strictly continuous and $\|\pi(f(t))\| \leq g(t)$ for $t \in X$.
\end{enumerate}
Then $\int f(t) \, d\mu(t)$ belongs to $D(\pi)$ and $\pi\bigl(\,\int f(t) \, d\mu(t)\,\bigr) = \int \pi(f(t)) \, d\mu(t)$.
\end{result}
\begin{demo}
Let $G$ denote the graph of $\pi$ in $M(A) \oplus M(B)$. By assumption, $G$ is strictly closed in $M(A) \oplus M(B)$.
Take $\om \in (M(A) \oplus M(B))'$ which is strictly continuous and such that $\om = 0$ on $G$. So there exist $\th \in M(A)'$ and $\eta \in M(B)'$ which are both strictly continuous and such that $\om(x,y) = \th(x) + \eta(y)$ for $x \in M(A)$ and $y \in M(B)$.

Notice that the functions $\th \circ f$ and $\eta \circ \pi \circ f$
are continuous and are dominated by a $\mu$-integrable function, so
they are $\mu$-integrable. By equality~\ref{vgl5} in the discussion before
this proposition, we have moreover that $\th(\int f(t) \,d\mu(t)) =
\int \th(f(t)) \, d\mu(t)$ and $\eta(\int \pi(f((t)) \, d\mu(t)) =
\int \eta(\pi(f((t))) \, d\mu(t)$. Hence
\begin{eqnarray*}
& & \om\bigl(\, \int f(t) \, d\mu(t) , \int \pi(f(t)) \, d\mu(t)\,\bigr)
= \th(\int f(t) \,d\mu(t)) + \eta(\int \pi(f((t)) \, d\mu(t)) \\
& & \spat = \int \th(f(t)) \, d\mu(t) + \int \eta(\pi(f((t))) \, d\mu(t)
= \int \bigl(\,\th(f(t)) , \eta(\pi(f(t)))\,\bigr) \,\, d\mu(t) \ .
\end{eqnarray*}
Let $t \in X$. Since $(f(t),\pi(f(t))) \in G$, we see that
$\bigl(\,\th(f(t)),\eta(\pi(f(t)))\,\bigr) = \om(f(t),\pi(f(t))) = 0$.
This implies that $\om\bigl(\, \int f(t) \, d\mu(t) , \int \pi(f(t)) \, d\mu(t)\,\bigr) = 0$. So Hahn-Banach implies that the couple \newline
$\bigl(\, \int f(t) \, d\mu(t) , \int \pi(f(t)) \, d\mu(t)\,\bigr)$ belongs to $G$.
\end{demo}

In a similar way, one proves the next result. This result is intended for use in connection with the GNS-map $\overline{\la}_\vfi$ of a lower semi-continuous weight $\vfi$ (see proposition \ref{weight.prop11}.)

\begin{result}
Consider a \cst-algebra $A$, a Hilbert space $H$  and a strictly closed linear map $\la$ from $D(\la) \subseteq M(A)$ into $H$. Let $X$ be a locally compact space  with regular Borel measure $\mu$ and consider a function $f : X \rightarrow D(\la)$  such that there exists a $\mu$-integrable function $g : X \rightarrow \R^+$ such that
\begin{enumerate}
\item $f$ is strictly continuous and $\|f(t)\| \leq g(t)$ for $t \in \R$.
\item $\la \circ f$ is  continuous and $\|\la(f(t))\| \leq g(t)$ for $t \in \R$.
\end{enumerate}
Then $\int f(t) \, d\mu(t)$ belongs to $D(\la)$ and $\la\bigl(\,\int f(t) \, d\mu(t)\,\bigr) = \int \la(f(t)) \, d\mu(t)$.
\end{result}

\bigskip

The following basic integration/strict convergence result turns out to be extremely useful when working with norm continuous one-parameter groups.

\medskip

\begin{result} \label{int.res1}
Consider a \cst-algebra $A$, a locally compact group $G$ with regular Borel measure $\mu$ and a norm continuous one-parameter group  $\al$ of $G$ on $A$. Let $f : G \rightarrow \C$ be a continuous integrable function, $(a_i)_{i \in I}$ a bounded net in $M(A)$ and $a \in M(A)$ such that $(a_i)_{i \in I}$ converges strictly to $a$.
Then $\bigl(\, \int  f(t) \, \al_t(a_i)\, d\mu(t) \,\bigr)_{i \in I}$ is a bounded net in
$M(A)$ which converges strictly to $\int  f(t) \, \al_t(a)\, d\mu(t)$.
\end{result}
\begin{demo}
Define $g = |f|$, which is a positive continuous integrable function.
Choose $b \in A$.

Let $i \in I$. Then we have that
\begin{eqnarray}
& & \| \, \int f(t) \, \al_t(a_i) \, b \, d\mu(t) - \int f(t) \, \al_t(a) \, b \, d\mu(t) \|
= \| \, \int f(t) \, \al_t(a_i -  a) \, b \, d\mu(t) \, \| \nonumber \\
& & \spat \leq  \, \int g(t) \, \| \al_t(a_i - a ) \, b \| \, d\mu(t)
= \int g(t) \, \| \al_t((a_i - a) \al_{-t}(b)) \| \, d\mu(t) \nonumber \\
& & \spat = \int g(t) \, \| (a_i - a)\al_{-t}(b) \| \, d\mu(t) \ .
\label{vgl6}
\end{eqnarray}

Take $\vep > 0$. By assumption, there exists a number $M > 0$ such that $\|a_i\| \leq M$ for $i \in I$ and $\|a\| \leq M$.
Because $g$ is integrable, there exists a nonempty compact subset $K$ of $G$ such that $\int_{G \setminus K} g(t) \, d\mu(t) \leq \frac{\vep}{4M(\|b\|+1)}$. Because $\|(a_i - a) \,\al_{-t}(b)\| \leq 2M \|b\|$ for $i \in I$ and $t \in \R$, this implies that
\begin{eqnarray*}
& & \int g(t) \, \| (a_i - a)\al_{-t}(b) \| \, d\mu(t)
\\
& & \spat = \int_{G \setminus K} g(t) \, \| (a_i - a)\al_{-t}(b) \| \, d\mu(t)
+ \int_{K} g(t) \, \| (a_i - a)\al_{-t}(b) \| \, d\mu(t) \\
& & \spat \leq 2 M \|b\| \, \int_{G \setminus K} g(t)  \, d\mu(t)
+ \int_{K} g(t) \, \| (a_i - a)\al_{-t}(b) \| \, d\mu(t)  \\
& & \spat \leq \int_{K} g(t) \, \| (a_i - a)\al_{-t}(b) \| \, d\mu(t)
+ \frac{\vep}{2}
\end{eqnarray*}
for every $i \in I$.

\medskip


The continuity of the function $G \rightarrow A : t \mapsto g(t)\,\al_t(b)$ and the usual compactness argument allows to find elements $t_1,\ldots,t_n \in K$ and open subsets $O_1,\ldots,O_n$ in $G$ such that $K \subseteq O_1 \cup \ldots \cup O_n$, $t_r \in O_r$ for all $r \in \{1,\ldots,n\}$ and $\| g(t) \, \al_t(b) - g(t_r) \, \al_{t_r}(b) \| \leq \frac{\vep}{8 M (\mu(K) + 1)}$ for all $r \in \{1,\ldots,n\}$ and $t \in O_r$.

Hence 
$$\| g(t) \, (a_i - a) \, \al_t(b) - g(t_r)\, (a_i - a) \, \al_{t_r}(b) \| \leq \frac{\vep}{4 (\mu(K) + 1)}  $$
for all $r \in \{1,\ldots,n\}$ and $t \in O_r$.

Since $(a_i)_{i \in I}$ converges strictly to $a$, we get the existence of an element $i_0 \in I$ such that \newline $\|g(t_r)\, (a_i - a)\, \al_{t_r}(b)\| \leq \frac{\vep}{4 (\mu(K) + 1)}$ for all $i \in I$ with $i \geq i_0$ and $r \in \{1,\ldots,n\}$.

Then it is clear that $\| g(t) \, (a_i - a) \, \al_t(b)\| \leq \frac{\vep}{2 (\mu(K) + 1)}$ for all $t \in K$ and $i \in I$ with $i \geq i_0$.

Hence we get
$$\int g(t) \, \| (a_i - a)\,\al_{-t}(b) \| \, d\mu(t) \leq \vep$$ for all $i \in I$ with $i \geq i_0$. Hence the calculation in inequality~\ref{vgl6} implies that
$$\| \, \int f(t) \, \al_t(a_i) \, b \, d\mu(t) - \int f(t) \, \al_t(a) \, b \, d\mu(t) \| \leq \vep$$
for all $i \in I$ with $i \geq 0$. So we have proven that
$$\bigl( \, \int f(t) \, \al_t(a_i)\, b \, d\mu(t) \,\bigr)_{i \in I} \ \ \rightarrow \ \ \int f(t) \, \al_t(a)\, b \, d\mu(t) \ .$$
Similarly, one proves that
$$\bigl( \, \int f(t) \, b\,\al_t(a_i) \, d\mu(t) \,\bigr)_{i \in I} \ \ \rightarrow \ \ \int f(t) \,b\, \al_t(a)\,  d\mu(t) \ .$$
\end{demo}

\bigskip\medskip

\section{Appendix : a technical lemma}

The next lemma is due to Jan Verding \& A. Van Daele and plays an
important role when studying weights arising from an inverse
GNS-construction.

\begin{lemma} \label{app.lem1}
Let $E$ be a normed space, $H$ a Hilbert space and $\la$ a linear mapping
from $D(\la) \subseteq E$ into $H$. Let $(x_i)_{i \in I}$ be a net in $D(\la)$ and $x$ an element in $E$ such that $(x_i)_{i\in I}$ converges to $x$ and
$(\la(x_i))_{i\in I}$ is bounded. Then there exists a sequence $(y_n)_{n=1}^\infty$ in the convex hull of $\{ x_i \mid i \in I \}$
and an element $v \in H$ such that $(y_n)_{n=1}^\infty$ converges to $y$
and $(\la(y_n))_{n=1}^\infty$ converges to $v$.
\end{lemma}
\begin{demo}
By the Banach-Alaoglu theorem, there exists a subnet $(x_{i_j})_{j\in J}$ of $(x_i)_{i\in I}$ and $v \in H$ such that
$(\la(x_{i_j}))_{j \in J}$ converges to $v$ in the weak topology on $H$. (remember that a Hilbert space is self-dual)

\medskip

Fix $n \in \N$. Then there exists $j_n \in J$ such that $\|x_{i_j} - x\| \leq \frac{1}{n}$ for all $j \in J$ with $j \geq j_n$.

Now $v$ belongs to the weak closed convex hull of the set $\{\la(x_{i_j}) \mid j \in J \text{ such that } j \geq j_n \}$, which is the same as  the norm closed convex hull.

Therefore there exist $\lambda_1,\ldots\!,\lambda_m \in \R^+$
with $\sum_{k=1}^m \lambda_k = 1$ and elements $\alpha_1,\ldots\!,\alpha_m \in J$  with $\alpha_1,\ldots\!,\alpha_m \geq j_n$  such that
$$ \| v - \sum_{k=1}^m \lambda_k \la(x_{i_{\alpha_k}}) \| \leq \frac{1}{n} .$$

Put $y_n =  \sum_{k=1}^{m} \lambda_k x_{i_{\alpha_k}}$. Then
$y_n \in D(\la)$, and $\la(y_n) = \sum_{k=1}^m \lambda_k \la(x_{i_{\alpha_k}})$.

So we have immediately that $\|v - \la(y_n)\| \leq \frac{1}{n}$.

Furthermore,
$$\|x - y_n\| = \left\|\sum_{k=1}^{m}\lambda_k(x - x_{i_{\alpha_k}})\right\| \leq
\sum_{k=1}^{m} \lambda_k  \frac{1}{n} = \frac{1}{n}$$

\medskip

Hence, we find that $(y_n)_{n=1}^\infty$ converges to $y$ and
that $(\la(y_n))_{n=1}^\infty$ converges to $v$.
\end{demo}

\medskip

By modifying the argument a little bit, we can also get a similar result in the $\si$-strong$^*$ topology on a von Neumann algebra.

\begin{lemma} \label{app.lem2}
Consider a von Neumann algebra $M$ acting on a Hilbert space $K$. Let $H$ be another Hilbert space and $\la$ a linear mapping
from $D(\la) \subseteq M$ into $H$. Let $(x_i)_{i \in I}$ be a net in $D(\la)$ and $x$ an element in $M$ such that $(x_i)_{i\in I}$ converges strongly$^*$ to $x$ and
$(\la(x_i))_{i\in I}$ is bounded. Then there exists a net $(y_j)_{j \in J}$ in the convex hull of $\{\, x_i \mid i \in I \, \}$
and an element $v \in H$ such that $(y_j)_{j \in J}$ converges strongly$^*$ to $y$
and $(\la(y_j))_{j \in J}$ converges to $v$.
\end{lemma}
\begin{demo}
Because the net $(\la(x_i))_{i \in I}$ is bounded, the Banach-Alaoglu theorem implies the existence of a subnet $(x_{i_l})_{l\in L}$ of $(x_i)_{i \in I}$ and $v \in H$ such that
$(\la(x_{i_l}))_{l \in L}$ converges to $v$ in the weak topology on $H$.

\medskip

Define $J = \{\, (E,n) \mid E \text{ a finite subset of } K \text{ and } n \in \N \, \}$ and put the natural order on $J$: If $(E,m),(F,n)$ $\in J$, then we say that $(E,m) \leq (F,n) \Leftrightarrow E \subseteq F$ and $m \leq n$.  Then $J$ becomes a directed set.

\medskip

Choose $j=(E,n) \in J$. Because $(x_{i_l})_{l \in L}$ converges strongly$^*$ to $x$, there exist an element $l_0 \in L$ such that
$$\|x_{i_l} \, w - x \, w\| \leq \frac{1}{n} \hspace{1.5cm} \text{ and } \hspace{1.5cm} \|x_{i_l}^* \, w - x^* \, w\| \leq  \frac{1}{n}$$
for $w \in E$ and $l \in L$ with $l_0 \leq l$.

Now $v$ belongs to the weak closed convex hull of the set $\{\la(x_{i_l}) \mid l \in L \text{ such that } l \geq l_0 \}$, which is the same as  the norm closed convex hull.
Therefore there exist $\lambda_1,\ldots\!,\lambda_m \in \R^+$
with $\sum_{r=1}^m \lambda_r = 1$ and elements $\alpha_1,\ldots\!,\alpha_m \in L$  with $\alpha_1,\ldots\!,\alpha_m \geq l_0$  such that
$$ \bigl\| \, v - \sum_{r=1}^m \lambda_r \lafi(x_{i_{\alpha_r}}) \,\bigr\| \leq \frac{1}{n} .$$

Put $y_j =  \sum_{r=1}^{m} \lambda_r \, x_{i_{\alpha_r}}$. Then
$y_j \in D(\la)$, and $\la(y_j) = \sum_{r=1}^m \lambda_r \, \lafi(x_{i_{\alpha_r}})$.
Therefore we have immediately that $\|v - \la(y_j)\| \leq \frac{1}{n}$.
Furthermore,
$$\|x \,w - y_j\,w \| = \left\|\,\sum_{r=1}^{m}\lambda_r\, (x \, w - x_{i_{\alpha_r}}\, w)\,\right\| \leq
\sum_{r=1}^{m} \lambda_r \, \frac{1}{n} = \frac{1}{n} \ .$$
for all $w \in E$.
Similarly, one shows that $\|x^* \, w - y_j^* \, w \| \leq \frac{1}{n}$ for all $w \in E$

\medskip

It is now easy to see that $(y_j)_{j \in J}$ converges strongly$^*$ to $y$ and that $(\la(y_j))_{j \in J}$ converges to $v$.
\end{demo}

\begin{remark}  \rm
Suppose that there exists a positive number $M$ such that $\|x_i\| \leq M$ for all $i \in I$. Then clearly $\|y_j\| \leq M$ for all $j \in J$.
\end{remark}

\medskip

\begin{corollary}
Consider a von Neumann algebra $M$ acting on a Hilbert space $K$. Let $H$ be another Hilbert space and $\la$ a linear mapping
from $D(\la) \subseteq M$ into $H$. Let $(x_i)_{i \in I}$ be a net in $D(\la)$ and $x$ an element in $M$ such that $(x_i)_{i\in I}$ converges $\si$-strongly$^*$ to $x$ and
$(\la(x_i))_{i\in I}$ is bounded. Then there exists a net $(y_j)_{j \in J}$ in the convex hull of $\{\, x_i \mid i \in I \, \}$
and an element $v \in H$ such that $(y_j)_{j \in J}$ converges $\si$-strongly$^*$ to $y$
and $(\la(y_j))_{j \in J}$ converges to $v$.
\end{corollary}
\begin{demo}
Define a new Hilbert space $K' = K \ot l^2(\N)$ and a new von Neumann algebra $M'$ on $K'$ by $M' = M \ot 1$. Define the $^*$-isomorphism $\pi : M \rightarrow M'$ such that $\pi(x) = x \ot 1$ for $x \in M$.
It is then easy to check for any net $(z_l)_{l \in L}$ in $M$ and any element $z \in M$ that $(z_l)_{l \in L}$ converges $\si$-strongly$^*$ to $z$
$\Leftrightarrow$ $(\pi(z_l))_{l \in L}$ converges strongly$^*$ to $\pi(z)$.
Define now a linear map $\la'$ from $D(\la') \subseteq M'$ into $H$ by $\la' = \la \circ \pi^{-1}$. Now apply the previous lemma to $\la'$, the net $(\pi(x_i))_{i \in I}$ and the element $\pi(x)$ and transfer everything back to $M$.
\end{demo}

\end{document}